\documentclass[11pt]{article}
\usepackage[margin=1in]{geometry}
\usepackage{graphicx}
\usepackage{enumitem}
\usepackage{amsfonts,amsmath,amsthm}
\usepackage{amssymb}
\usepackage[T1]{fontenc}
\usepackage{soul}
\usepackage{bm}
\usepackage{float}
\usepackage{xcolor}
\usepackage[noend]{algpseudocode}
\usepackage{algorithm}
\usepackage{hyperref}
\usepackage{xr}
\numberwithin{equation}{section}
\usepackage{subcaption} 
\usepackage{wrapfig}
\usepackage[onehalfspacing]{setspace}
\usepackage{booktabs}
\usepackage{threeparttable}
\usepackage{tikz}
\usetikzlibrary{arrows.meta}  
\usetikzlibrary{shadows}  
\usetikzlibrary{positioning}  
\usetikzlibrary{decorations.pathreplacing}  
\usetikzlibrary{shapes}  
\usetikzlibrary{calc}

\allowdisplaybreaks

\newtheorem{prop}{Proposition}
\newtheorem{rema}{Remark}

\newtheorem{assu}{Assumption}
\newtheorem{exam}{Example}
\newtheorem{theo}{Theorem}
\newtheorem{coro}{Corollary}

\newcommand{\n}{\texttt{n}}

\makeatletter
\renewcommand{\ALG@name}{Routine}
\makeatother

\title{Disjunctive Benders Decomposition}
\author{Kaiwen Fang\footnotemark[2] 
 \and Inho Sin\footnotemark[2] \and Geunyeong Byeon\footnotemark[2] \footnotemark[4]}
\date{}
\begin{document}

 \footnotetext[2]{School of Computing and Augmented Intelligence, Arizona State University, Tempe AZ, USA. Emails: \texttt{\{kfang11, inhosin, geunyeong.byeon\}@asu.edu}.}%
\footnotetext[4]{Corresponding author.}

\maketitle

\begin{abstract}
We propose an enhancement to Benders decomposition (BD) that generates valid inequalities for the convex hull of the Benders reformulation, addressing the limitation that classical BD cuts are typically tight only for the continuous relaxation. Our method integrates disjunctive programming with BD and introduces a routine that leverages existing cut-generating oracles as-is to construct convex hull inequalities.
For mixed-binary linear programs, the approach removes the need to solve the master problem as a mixed-integer program{\color{black}, even with separable subproblems. It builds on a unified normalization framework for cut-generating programs, encompassing norm-based, reverse polar, and right-hand-side normalization, and enabling the design of new normalization schemes with streamlined analysis of supporting cuts.}
Computational results on large-scale instances show substantial reductions in branch-and-bound nodes—often by orders of magnitude—while consistently outperforming commercial solvers on selected problem classes.

\vspace{2mm}

\end{abstract}

\section{Introduction}
Benders decomposition (BD) is a widely used method for solving mixed-integer programming problems and has gained increasing attention for its effectiveness in tackling large-scale instances. Its versatility has led to numerous applications in areas such as power systems, transportation, and supply chain management (see \cite{rahmaniani2017benders} for a comprehensive review). This paper advances BD by addressing its inherent limitations through the integration of disjunctive programming theory \cite{balas1998disjunctive,balas2018disjunctive}.

\subsection{Problem statement}
We focus on a mixed-integer linear programming (MILP) problem in the following form:
\begin{equation}
\begin{aligned}
\min_{\boldsymbol{x},\boldsymbol{y}} \ & \boldsymbol{c}^{\top} \boldsymbol{x} + \boldsymbol{d}^{\top} \boldsymbol{y}  \\
\mbox{s.t.} \ & A\boldsymbol{x} + B\boldsymbol{y} \geq \boldsymbol{b},\\
              & \boldsymbol{x} \in \mathcal{X} := \{\boldsymbol{x} \in [\underline{x}_1, \overline{x}_1] \times \cdots \times [\underline{x}_{n_x}, \overline{x}_{n_x}]: D \boldsymbol{x} \geq \boldsymbol{h}, x_i \in \mathbb{Z}, \forall i \in \mathcal{I}\},
\end{aligned}
\tag{\texttt{MILP}}
\label{prob:milp}
\end{equation}
where $n_x$ is the dimension of $\boldsymbol{x}$, and for each $j \in [n_x] := \{1, \dots, n_x\}$, $\underline{x}_j$ and $\overline{x}_j$ denote the lower and upper bounds of $x_j$. The set $\mathcal{I} \subseteq [n_x]$ represents the indices of $\boldsymbol{x}$ constrained to take integer values. Without loss of generality, we assume $\mathcal{I} = \{1, \cdots, \ell\}$, where $\ell = |\mathcal{I}|$. All matrices and vectors, i.e., $A$, $B$, $D$, $\boldsymbol c$, $\boldsymbol d$, $\boldsymbol b$, $\boldsymbol h$, are real and of appropriate dimensions. Fixing $\boldsymbol x$ simplifies the problem, as setting $\boldsymbol{x} = \boldsymbol{\hat{x}}$ reduces the problem to a linear programming (LP) problem:  
\begin{equation}
\min_{\boldsymbol y}\{ \boldsymbol{d}^{\top} \boldsymbol{y}: B \boldsymbol{y} \geq \boldsymbol{b} - A \boldsymbol{\hat{x}}\}.
\tag{\texttt{sub}(${\hat{\boldsymbol x}}$)}
\label{prob:sub}
\end{equation}
Let $\Pi$ denote the dual feasible region of \ref{prob:sub}, that is $\Pi = \{\boldsymbol{\pi} \geq 0: B^{\top} \boldsymbol{\pi} = \boldsymbol{d}\}$. We make the following assumptions:
\begin{assu}
\begin{enumerate}
    \item[(i)] For each $j \in \{1, \cdots, n_x\}$, $-\infty < \underline{x}_j < \overline{x}_j < \infty$, i.e., $\mathcal{X}$ is compact.
    \item[(ii)] $\Pi \neq \emptyset$.
\end{enumerate}\label{assum}
\end{assu}

\begin{rema}
Assumption \ref{assum}(i) can be relaxed with some modifications to the derivations presented in this paper, but we impose it here for simplicity. Assumption \ref{assum}(ii) is not restrictive, as $\Pi = \emptyset$ would imply that \ref{prob:sub} is infeasible or unbounded for any $\boldsymbol{\hat{x}} \in \mathbb{R}^{n_x}$, rendering the problem pathological.
\end{rema}

We define a value function $f$ that maps $\boldsymbol{\hat{x}}$ to the optimal objective value of \ref{prob:sub}. Its effective domain $\text{dom} f := \{\boldsymbol{x} \in \mathbb{R}^{n_x} : f(\boldsymbol{x}) < \infty\}$ corresponds to the set of $\boldsymbol{x}$'s for which \texttt{sub}($\boldsymbol x$) is feasible. Under Assumption \ref{assum}(ii), \ref{prob:sub} {\color{black}enjoys strong duality} for any $\boldsymbol{\hat{x}} \in \mathbb{R}^{n_x}$. Thus, we have:   
\begin{align}
f(\boldsymbol{\hat{x}}) &= \max_{\boldsymbol{\pi} \in \Pi} (\boldsymbol{b} - A \boldsymbol{\hat{x}})^{\top} \boldsymbol{\pi} \label{prob:sub:dual} \tag{\texttt{sub-dual}($\hat{\boldsymbol{x}}$)}\\
&= \begin{cases}
\max_{\boldsymbol{\hat{\pi}} \in \mathcal{J}} (\boldsymbol{b} - A \boldsymbol{\hat{x}})^{\top} \boldsymbol{\hat{\pi}} & \mbox{if } (\boldsymbol{b} - A \boldsymbol{\hat{x}})^{\top} \boldsymbol{\tilde{\pi}} \leq 0, \, \forall \boldsymbol{\tilde{\pi}} \in \mathcal{R} \\
\infty & \mbox{o.w.}
\end{cases}\nonumber
\end{align}
where $\mathcal{J}$ and $\mathcal{R}$ respectively represent the set of all extreme points and extreme rays of $\Pi$. This implies $\mbox{dom}f = \{\boldsymbol{x} \in \mathbb{R}^{n_x}: 0 \geq \boldsymbol{\tilde{\pi}}^{\top} (\boldsymbol{b}-A\boldsymbol{x}), \, \forall \boldsymbol{\tilde{\pi}} \in \mathcal{R}\}$
and
$\mbox{epi}f := \{(\boldsymbol{x}, t) \in \mathbb{R}^{n_x+1}: t \geq f(\boldsymbol{x})\} = \{(\boldsymbol{x},t)\in \mathbb{R}^{n_x+1}: t \geq \boldsymbol{\hat{\pi}}^{\top} (\boldsymbol{b}-A\boldsymbol{x}), \, \forall \boldsymbol{\hat{\pi}} \in \mathcal{J}\},$
i.e., $\mbox{epi}f$ and $\text{dom}f$ are polyhedral, and every facet of $\mbox{epi}f$ (and $\mbox{dom}f$) is associated with an extreme point (and an extreme ray) of the dual feasible region $\Pi$. 
We define a set
\begin{equation*}
\begin{aligned}
\mathcal{L} := \{(\boldsymbol{x}, t) \in \mathbb{R}^{n_x+1} : (\boldsymbol{x}, t) \in \text{epi}f, \, \boldsymbol{x} \in \text{dom}f\}.
\end{aligned}
\end{equation*} 
Using $\mathcal{L}$, \eqref{prob:milp} can be posed equivalently in the $(\boldsymbol{x},t)$-space, which is called \emph{Benders reformulation}:  
\vspace{-5mm}
\begin{subequations}
\begin{align}
\min_{\boldsymbol{x} \in \mathcal{X}, t \in \mathbb{R}} \{ \boldsymbol{c}^{\top} \boldsymbol{x} + t : (\boldsymbol{x}, t) \in \mathcal{L}\} = \min_{\boldsymbol{x} \in \mathcal{X}, t \in \mathbb{R}} \ & \boldsymbol{c}^{\top} \boldsymbol{x} + t  \\
\mbox{s.t.} \ & t \geq \boldsymbol{\hat{\pi}}^{\top} (\boldsymbol{b}-A\boldsymbol{x}), \ \forall \boldsymbol{\hat{\pi}} \in \mathcal{J}, \label{eq:benders-reformulation:epi}\\
                & 0 \geq \boldsymbol{\tilde{\pi}}^{\top} (\boldsymbol{b}-A\boldsymbol{x}), \ \forall \boldsymbol{\tilde{\pi}} \in \mathcal{R}. \label{eq:benders-reformulation:dom}
\end{align}
\label{prob:benders-reformulation}
\end{subequations}
\vspace{-5mm}
\begin{figure}[b!]
    \centering
    \begin{subfigure}[b]{0.45\textwidth}
        \centering
        \begin{tikzpicture}[
            >={Stealth},
            thick,
            math node/.style={font=\Large, text=black},
            typical oracle node/.style={
                text=black,
                font=\Large\bfseries,
                align=center,
            },
            hyperplane node/.style={align=left, font=\Large},
            scale=0.7, transform shape
        ]
            \node[math node] (input) at (0,0) {$(\boldsymbol{\hat{x}}, \hat{t})$};
            
            \node[typical oracle node] (oracle) at (3.3,0) {\textsc{typicalOracle}};
            
            \node[hyperplane node] (output) at (7.5,0) {A hyperplane\\  separating\\$(\boldsymbol{\hat{x}}, \hat{t})$ from $\mathcal{L}$};
            
            \draw[->, line width=1pt] (input) -- (oracle);
            \draw[->, line width=1pt] (oracle) -- (output);
        \end{tikzpicture}
        \caption{Conventional cut-generating oracles.}
        \label{fig:oracle-typical}
    \end{subfigure}
    \hfill
    \begin{subfigure}[b]{0.5\textwidth}
        \centering
        \begin{tikzpicture}[
            >={Stealth},
            thick,
            math node/.style={font=\Large, text=black},
            proposed oracle node/.style={
                text=black,
                font=\Large\bfseries,
                align=center,
            },
            hyperplane node/.style={align=left, font=\Large},
            scale=0.7, transform shape
        ]
            \node[math node] (input) at (0,0) {$(\boldsymbol{\hat{x}}, \hat{t})$};
            
            \node[proposed oracle node] (oracle) at (3.8,0) {\textsc{\color{black}disjunctiveOracle}};
            
            \node[hyperplane node] (output) at (9.4,0) {A hyperplane \\separating $(\boldsymbol{\hat{x}}, \hat{t})$ from\\
            $\overline{\operatorname{conv}}((\mathcal{X} \times \mathbb{R}) \cap \mathcal{L})$};

            \draw[->, line width=1pt] (input) --  (oracle);
            \draw[->, line width=1pt] (oracle) -- (output);
        \end{tikzpicture}
        \caption{Proposed cut-generating oracles.}
        \label{fig:oracle-proposed}
    \end{subfigure}
    \caption{Comparison between the proposed approach and conventional methods.}
    \label{fig:oracle}
\end{figure}

BD is an iterative algorithm for solving \eqref{prob:benders-reformulation} which often involves exponential number of constraints for describing $\mathcal{L}$. It begins by replacing $\mathcal{L}$ with a polyhedral relaxation in \eqref{prob:benders-reformulation}, forming what is known as the \emph{master problem}, which generates candidate solutions. Given a candidate solution $(\boldsymbol{\hat{x}}, \hat{t})$, an oracle verifies whether $(\boldsymbol{\hat{x}}, \hat{t}) \in \mathcal{L}$. If not, the oracle constructs a hyperplane that separates $(\boldsymbol{\hat{x}}, \hat{t})$ from $\mathcal{L}$ and adds it to the master problem, repeating this process iteratively.

The effectiveness of BD depends on two key factors: how the candidate solutions are generated and how the cut-generating oracle is designed. The candidate solutions can either be optimal solutions of the evolving master problem or a sequence of incumbent solutions identified during a branch-and-bound process. We refer to these approaches as \textsc{BendersSeq} and \textsc{BendersBnB}, respectively. The \textsc{BendersBnB} method is particularly effective for large-scale problems, as it avoids solving the mixed-integer master problem at every iteration \cite{fortz2009improved,fischetti2017redesigning,rahmaniani2017benders,maher2021implementing}.

The cut-generating oracle, another crucial component of BD, can take various forms. The classical oracle \cite{benders1962partitioning} solves \ref{prob:sub:dual}, computing $f(\boldsymbol{\hat{x}})$ and verifying whether the current solution $(\boldsymbol{\hat{x}}, \hat{t})$ violates the epigraph or domain conditions. If a violation is found, the oracle queries the solution or the unbounded ray of \ref{prob:sub:dual} to construct a violated inequality: either an {optimality cut} of the form \eqref{eq:benders-reformulation:epi}, separating the iterate from $\text{epi}f$, or a {feasibility cut} of the form \eqref{eq:benders-reformulation:dom}, separating it from $\text{dom}f$. {\color{black}The latter is often selected via normalization \cite{bonami2020implementing}.} A more advanced oracle may solve a modified subproblem to generate both optimality and feasibility cuts in a unified manner, often employing advanced cut-selection metrics \cite{fischetti2010note,brandenberg2021refined,hosseini2024deepest}. Alternatively, it may utilize a core point or analytical center to produce Pareto-optimal cuts \cite{magnanti1981accelerating,papadakos2008practical,naoum2013interior,seo2022closest} or solve \texttt{sub-dual}($\boldsymbol{x'}$) for a perturbed point $\boldsymbol{x'}$, often obtained as a convex combination of $\boldsymbol{\hat{x}}$ and a reference point \cite{ben2007acceleration,fischetti2017redesigning}. Recently, artificial intelligence-informed oracles have been developed to identify effective Benders cuts from the best candidates generated by the aforementioned oracles \cite{jia2021benders,deza2023machine}. All these oracle variants generate hyperplanes that separate the solution from $\mathcal{L}$, as illustrated in Figure \ref{fig:oracle-typical}. Throughout this paper, we refer to these oracles as \textsc{typicalOracle}. 

However, BD has an inherent limitation: the cuts generated by a \textsc{typicalOracle} are \emph{at best tight with respect to $\mathcal{L}$}. In other words, the oracle generating these cuts does not account for the integrality constraints on $\boldsymbol x$. Although such oracles have shown efficiency in various applications, their performance can degrade significantly when the integrality gap in \eqref{prob:benders-reformulation} is large—that is, when the convex hull of $(\mathcal{X} \times \mathbb{R}) \cap \mathcal{L}$ is much smaller than $\mathcal{L}$. This limitation impacts the efficiency of both \textsc{BendersSeq} and \textsc{BendersBnB} and we conjecture that it is a major contributor to the commonly observed tail-off convergence of BD \cite{byeon2022benders}, where substantial progress is often made during the initial iterations, but then the algorithm slows considerably thereafter, achieving only minimal reductions in the optimality gap over many iterations. This behavior is illustrated in the following example, which motivates the development of an oracle capable of generating cuts that are valid for the convex hull of $(\mathcal{X} \times \mathbb{R}) \cap \mathcal{L}$, as defined in Figure~\ref{fig:oracle-proposed}.

\subsection{Motivating example}\label{sec:exam}
We illustrate from a simple example that even with a state-of-the-art implementation, BD may need to enumerate a substantial portion of the discrete set $\mathcal{X}$ before reaching an optimal solution. 

\begin{figure*}[bt!]
  \centering
  \begin{subfigure}[t]{0.4\linewidth}  
    \includegraphics[width=\textwidth]{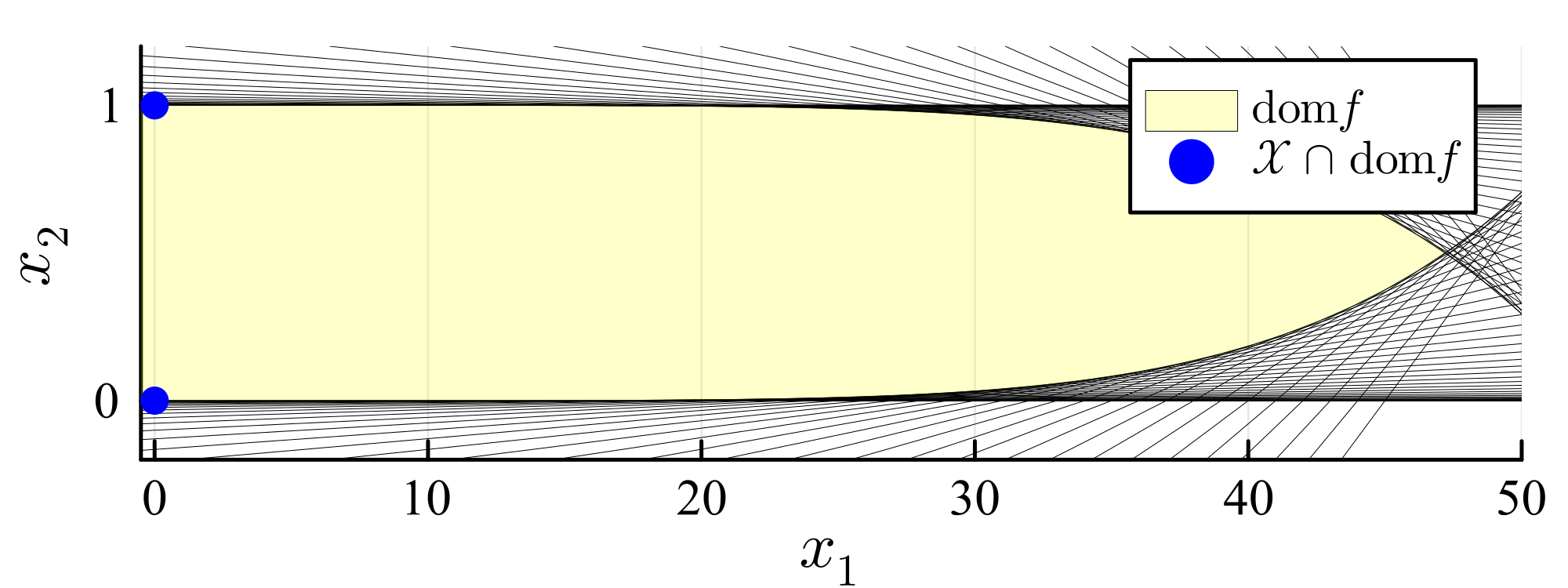}
    \caption{Instance visualization.}
    \label{fig:problem}
  \end{subfigure}
  \hfill  
  \begin{subfigure}[t]{0.55\linewidth}
    \includegraphics[width=\textwidth]{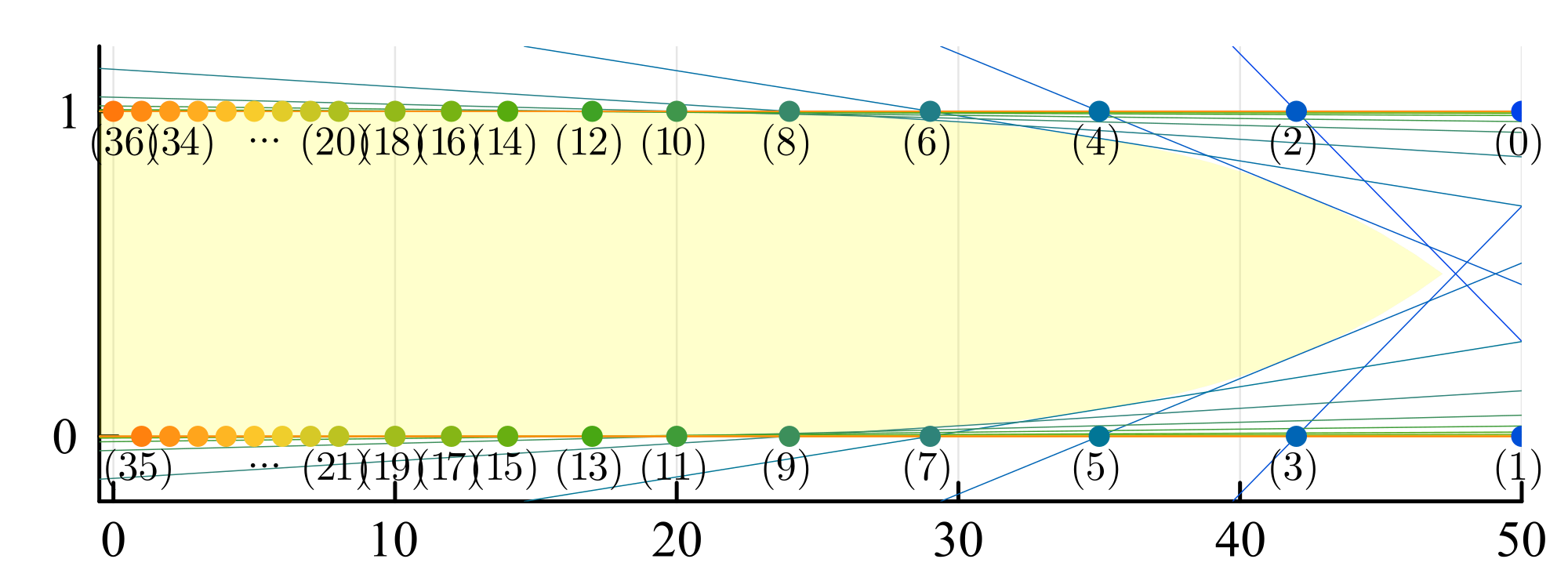}
    \caption{Iterates $\{\boldsymbol{x^{(k)}}\}$ and cuts produced by \textsc{BendersSeq}; colors associate $\boldsymbol{x^{(k)}}$ and the cuts.}
    \label{fig:vanillaBD}
  \end{subfigure}
  
  \vspace{1em}  
  
  \begin{subfigure}[t]{0.4\linewidth}
    \includegraphics[width=\textwidth]{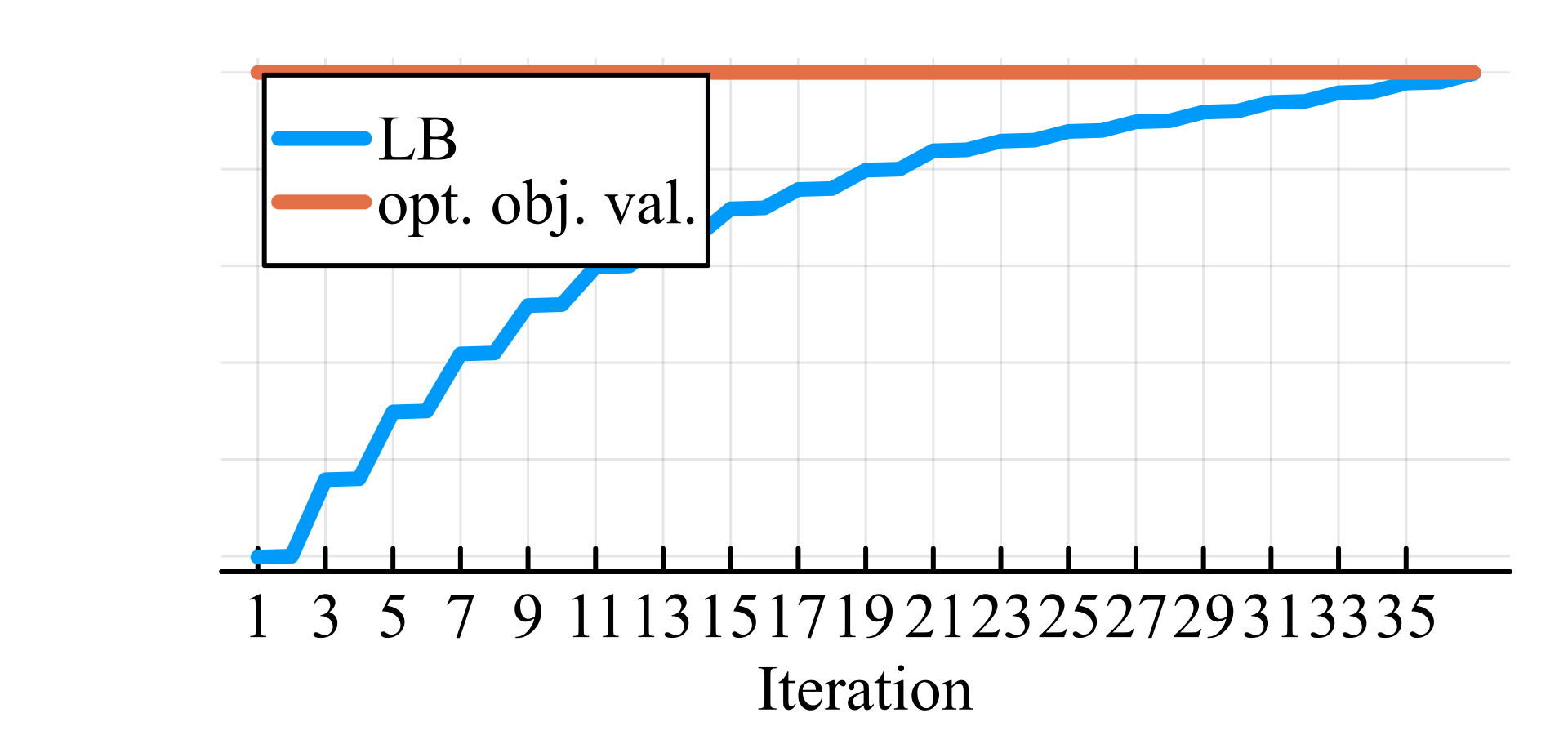}
    \caption{Tail-off convergence.}
    \label{fig:tailoff_iter}
  \end{subfigure}
  \hfill
  \begin{subfigure}[t]{0.55\linewidth}
    \includegraphics[width=\textwidth]{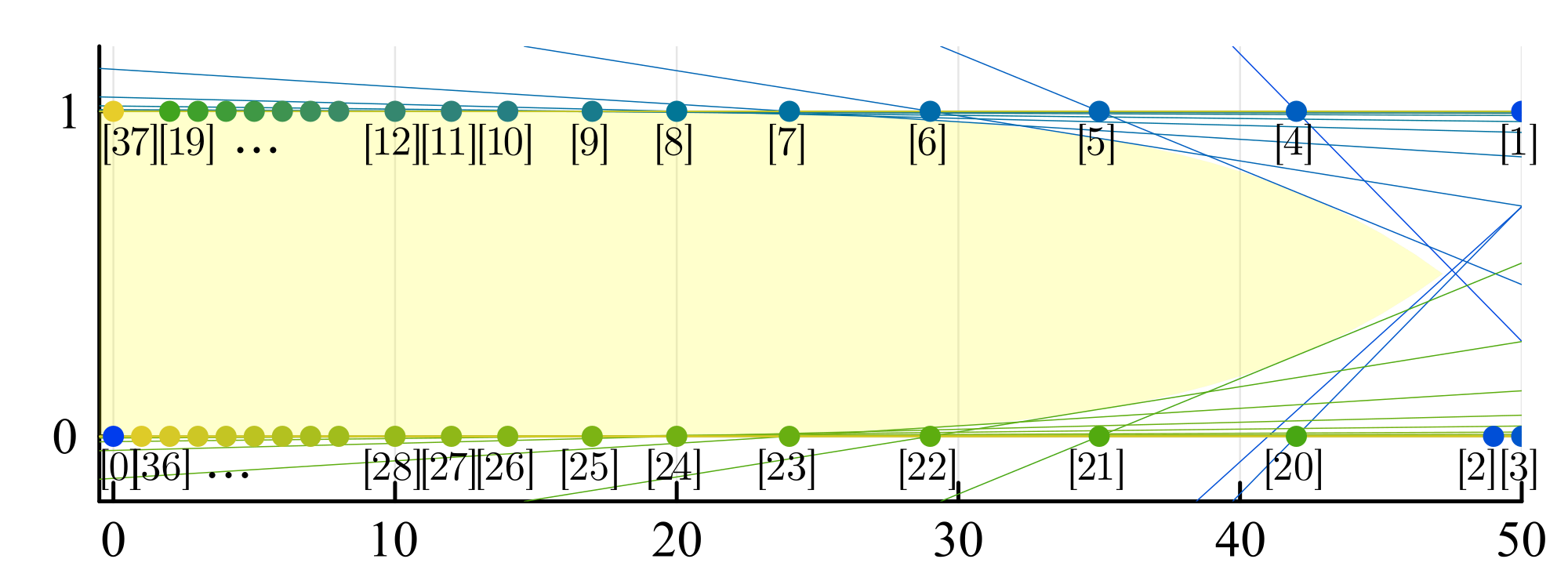}
    \caption{Candidate incumbent solutions $\{\boldsymbol{\hat{x}^{[k]}}\}$ and cuts produced by \textsc{BendersBnB}.}
    \label{fig:lazyBD}
  \end{subfigure}
  
  \caption{Visualization of BD's behavior; horizontal (or vertical) axes for $x_1$ (or $x_2$).}
  \label{fig:benders_behavior}
\end{figure*}

Consider an instance of \eqref{prob:benders-reformulation} with $\mathcal{X} = \{0,1,\cdots,50\} \times \{0,1\}$ and $\boldsymbol{c} = [-2; -1]$. Assume $f(\boldsymbol{x}) = 0$ for $\boldsymbol{x} \in \text{dom}f$, resulting in the Benders reformulation: $\min \{\boldsymbol{c}^{\top} \boldsymbol{x}: \boldsymbol{x} \in \mathcal{X} \cap \text{dom}f\}$. As depicted in Figure \ref{fig:problem}, $\text{dom}f$ is characterized by hyperplanes with gradually changing slopes so that only two points (marked in blue) are feasible, i.e., in $\mathcal{X} \cap \text{dom}f$. The optimal solution can be easily identified as $(0,1)$.

Figures \ref{fig:vanillaBD} and \ref{fig:lazyBD} depict the iterates $\{\boldsymbol{x^{(k)}}\}$ and the incumbent solutions $\{\boldsymbol{\hat{x}^{[k]}}\}$ generated by \textsc{BendersSeq} and \textsc{BendersBnB}, respectively, both equipped with an advanced \textsc{typicalOracle} proposed in \cite{fischetti2010note} and implemented via the Julia-CPLEX interface. Even in this simple example, both methods enumerate more than a quarter of $\mathcal{X}$. Figure \ref{fig:vanillaBD} highlights the trajectory of \textsc{BendersSeq}, which zigzags across more than a quarter of $\mathcal{X}$, resulting in the tailing-off convergence behavior depicted in Figure \ref{fig:tailoff_iter}. Meanwhile, although \textsc{BendersBnB} may identify a high-quality incumbent solution early on, as shown in Figure \ref{fig:lazyBD}, proving its optimality requires a substantial number of cuts, resulting in a final search tree with 16 nodes.

It is important to mention that this example is not considered pathological. For example, when dealing with bilevel programs involving convex followers, the Benders reformulations of these bilevel programs exhibit a similar pattern due to the incorporation of large coefficients to render the follower's optimality condition convex. These coefficients lead the Benders cuts to remove just one potential solution point during each iteration, a phenomenon that was observed in \cite{byeon2022benders}.  

If cutting planes are generated for the convex hull of $\mathcal{X} \cap \text{dom}f$, instead of $\text{dom}f$, the algorithm would terminate right away or in a few iterations. This example highlights why the simple in-out method \cite{ben2007acceleration,fischetti2017redesigning}, which separates $\lambda \boldsymbol{x^{(k)}} + (1-\lambda)\boldsymbol{\hat{x}^0}$ for some $\lambda \in [0,1]$ and a core point $\boldsymbol{\hat{x}^0}$, instead of $\boldsymbol{x^{(k)}}$, is often found to be effective \cite{maher2021implementing}. Nevertheless, this example suggests that the performance of the in-out method is limited and can vary significantly depending on the choice of $\boldsymbol{\hat{x}^0}$ and $\lambda$, which are challenging to determine a priori.

\subsection{Related literature}
Generating cuts for the convex hull of $(\mathcal{X} \times \mathbb{R}) \cap \mathcal{L}$ can yield significantly stronger cuts, reducing the total number of iterations in \textsc{BendersSeq} or shrinking the branch-and-bound tree in \textsc{BendersBnB}. Only recently have attempts been made to incorporate the integrality of $\boldsymbol x$ into Benders cut generation to make the cuts stronger \cite{bodur2017mixed,bodur2017strengthened,zou2019stochastic,rahmaniani2020benders}. In \cite{zou2019stochastic,rahmaniani2020benders}, the authors introduced binary variables to the subproblems to generate stronger Lagrangian-based cuts.
In \cite{bodur2017mixed}, a mixed-rounding procedure was applied to already-generated Benders cuts to further strengthen them.
Another relevant work, \cite{bodur2017strengthened}, introduced split cuts, termed project-and-cut split cuts, which are derived from a cut-generating problem formulated based on the original formulation \eqref{prob:milp}. Since the resultant cut-generating linear program can become enormous, they experimented with adding split cuts for the current master problem—i.e., for the convex hull of $(\mathcal{X} \times \mathbb{R}) \cap \widehat{\mathcal{L}}$, where $\widehat{\mathcal{L}}$ is an outer approximation of $\mathcal{L}$—as well as solving the proposed cut-generating problem separately for each scenario. They also proposed another method, termed cut-and-project split cuts, which incorporate valid constraints for the Lagrangian subproblems to the Benders subproblems, akin to approaches proposed for solving stochastic integer programming \cite{gade2014decomposition}. {\color{black}Other extensions of BD for stochastic integer programming focus on enabling cut generation in the presence of discrete subproblems \cite{sherali2002modification,sen2005c,lozano2022binary}, rather than strengthening cuts via integrality in the master variables.}

We provide a detailed comparison of the proposed approach with these existing methods in Section \ref{sec:discussion}.

\subsection{Contributions} 
We propose a new approach to enhancing BD by incorporating the generation of valid inequalities for the convex hull of the Benders reformulation. This addresses a key limitation of conventional BD, which typically generates cuts that can be at best tight with respect to the continuous relaxation. Our main contributions are summarized as follows:
\begin{itemize}[leftmargin=*] 
\item \emph{A general and scalable framework for generating disjunctive cuts for the Benders reformulation (i.e., \eqref{prob:benders-reformulation}).}
We develop a scalable routine that leverages existing cut-generating oracles as-is to produce disjunctive cuts. When fully executed, the routine constructs a supporting hyperplane for a set lying between the convex hull and the continuous relaxation. Notably, even if the routine terminates early, it still produces a valid inequality that cuts off part of the continuous relaxation, ensuring computational efficiency and practical applicability. This routine offers advantages over related approaches, particularly for large-scale problems: it neither introduces integer variables into the subproblem nor operates directly on the original formulation \eqref{prob:milp}, {\color{black}and it preserves the dual feasible region of the subproblem} (see Section~\ref{sec:discussion}).

{\color{black}\item \emph{A unified normalization framework for cut-generating programs.}
We introduce a support-function-based normalization that subsumes common schemes, including norm-based and reverse polar normalization. The resulting formulation admits a geometric dual interpretation as computing a gauge distance to the convex hull, unifying existing interpretations and facilitating the design of new normalization schemes with streamline analysis of supporting cuts.}

\item \emph{Advancements for mixed-binary linear programs.}
For mixed-binary linear programs, we show that a specialized integration of the proposed routine within BD eliminates the need to solve the master problem as a mixed-integer program {\color{black}theoretically—even for separable subproblems; to the best of our knowledge, this is the first such result.} Additionally, we extend the a posteriori strengthening and lifting procedure for lift-and-project cuts—an important class of disjunctive cuts for mixed-binary linear programs—into the BD framework. This extension provides a simple yet effective method for enhancing lift-and-project cuts produced by the proposed routine. Furthermore, it motivates the development of an approximate routine that efficiently identifies a lift-and-project cut whenever a violation exists.

\item \emph{Numerical experiments.} We conduct numerical experiments on large-scale instances of the uncapacitated facility location problem and the stochastic network interdiction problem, demonstrating the effectiveness of disjunctive cuts in improving the node efficiency of \textsc{BendersBnB}, measured through lower-bound improvement and branch-and-bound tree size reduction. The proposed approach solves instances that neither off-the-shelf solvers nor a state-of-the-art Benders implementation can solve, while reducing final tree sizes by up to an order of magnitude on difficult instances. Even in cases where no method reaches optimality, the proposed approach consistently attains the strongest lower bounds.
\end{itemize}

\subsection{Notation and preliminaries}
{\color{black}In this paper, vectors are denoted by boldface lowercase letters, e.g., $\boldsymbol{x}$, while scalars are denoted by standard lowercase letters.} $\mathbb{Z}$ and $\mathbb{R}$ respectively denote the set of integers and real numbers. For any set $\mathcal{S} \in \{ \mathbb{Z}, \mathbb{R} \}$, the notation $\mathcal{S}_{\geq 0}$ (or $\mathcal{S}_{> 0}$) refers to the set of elements in $\mathcal{S}$ that are greater than or equal to (or strictly greater than) zero. The set $[n]$ represents the set $\{1, 2, \dots, n\}$. The boldfaced $\boldsymbol{0}$ is a vector of zeros. The notation $\| \cdot \|_p$ represents the $l_p$-norm for some $p \in [1, \infty]$, and $q$ denotes the Hölder constant corresponding to $p$, satisfying $\frac{1}{p} + \frac{1}{q} = 1$. The vector $\boldsymbol{e_j} \in \mathbb{R}^{n_x}$ denotes the $j$-th unit vector, with 1 in the $j$-th element and zeros elsewhere. For a vector $\boldsymbol v$ and a scalar $a$, $(\boldsymbol{v}, a)$ denotes their concatenation, and $\text{Proj}_{x}(\mathcal{S})$ refers to the orthogonal projection of a set $\mathcal{S}$ onto the $\boldsymbol x$-space. The symbols $\texttt{!}$ and $\land$ denote the logical negation (complement) and logical conjunction (and), respectively.

{\color{black}
We recall several notions from convex analysis used throughout the paper. We use the standard notation $f^*$ to denote the \emph{conjugate} of a function $f$, defined as
\(
f^*(\boldsymbol{y}) := \sup_{\boldsymbol{x}} \{ \boldsymbol{y}^\top \boldsymbol{x} - f(\boldsymbol{x}) \}.
\)
Let $\mathcal S \subseteq \mathbb{R}^n$ be a set. The \emph{polar} of $\mathcal S$ is defined as
\(
{\mathcal S}^\circ := \{\boldsymbol{y} \in \mathbb{R}^n : \boldsymbol{y}^\top \boldsymbol{x} \le 1, \ \forall\boldsymbol{x} \in \mathcal S\}.
\)
For a convex cone $\mathcal K$, its polar is
\(
{\mathcal K}^\circ := \{\boldsymbol{y} \in \mathbb{R}^n : \boldsymbol{y}^\top \boldsymbol{x} \le 0, \ \forall \boldsymbol{x} \in \mathcal{K}\}.
\)
We denote by $\overline{\operatorname{conv}}(\mathcal{S})$ the closed convex hull of $\mathcal{S}$, respectively.
Let $\mathcal C \subseteq \mathbb{R}^n$ be a nonempty closed convex set.
The \emph{indicator function} of $\mathcal C$ is defined as
\[
\delta_{\mathcal C}(\boldsymbol{x}) :=
\begin{cases}
0, & \boldsymbol{x} \in \mathcal C,\\
+\infty, & \text{o.w.}
\end{cases}
\]
The \emph{support function} of $\mathcal C$ is
\[
\sigma_{\mathcal C}(\boldsymbol{x}) := \sup_{\boldsymbol{y} \in \mathcal C} \boldsymbol{x}^\top \boldsymbol{y},
\]
and they satisfy $\sigma_{\mathcal C}^* = \delta_{\mathcal C}$. The \emph{gauge function} associated with $\mathcal C$ is defined as
\[
\gamma_{\mathcal C}(\boldsymbol{x}) := \inf\{\tau \ge 0 : \boldsymbol{x} \in \tau \mathcal{C}\},
\]
and it satisfies
\(
\gamma_{\mathcal C}^* = \delta_{{\mathcal C}^\circ}.
\)
}

\subsection{Organization of the paper}
The paper is organized as follows. {\color{black}Section~\ref{sec:cgp} introduces a cut-generating program with a unified treatment of normalization via support functions and provides a geometric dual interpretation.} Section~\ref{sec:dbd} presents the proposed approach: Section~\ref{sec:dbd:disjunction} defines the disjunction in the Benders reformulation, Section~\ref{sec:cglp} develops a cut-generating problem, Section~\ref{sec:dbd:oracle} describes an oracle for solving this problem, and Section~\ref{sec:dbd:dbd} integrates these components into a disjunctive Benders decomposition. Section~\ref{sec:strenthening-lifting} presents extensions for mixed-binary linear programs. Section~\ref{sec:discussion} compares the proposed approach with related work, and Section~\ref{sec:experiment} reports computational results. Section~\ref{sec:conclusion} concludes.

{\color{black}
\section{Cut generating program with support-function normalization} \label{sec:cgp}
Cutting planes are typically generated by solving a cut-generating program (CGP). In this section, we introduce a CGP with a support-function-based normalization that serves as the foundation of our framework. We show that two widely used normalization schemes in CGPs can be expressed through a single support-function constraint. This unified formulation admits a geometric dual interpretation: the CGP computes a gauge distance from the separation point to the convex hull of the set, thereby casting cut generation as a distance-to-convex-hull problem under a user-specified gauge. This perspective recovers the known geometric interpretation of several popular normalization schemes while providing a simpler and more unified derivation. It also provides a systematic way to design new normalization schemes.

For a set $\mathcal S \subseteq \mathbb R^n$, let $\mathcal S^*$ denote its reverse polar cone \cite{balas2002lift}, that is, the cone of all valid inequalities for $\mathcal S$: $\mathcal S^* := \{(\boldsymbol{\alpha_x}, \alpha_0)  \in \mathbb{R}^{n+1}: \boldsymbol{\alpha_x}^\top \boldsymbol{x} \ge \alpha_0, \ \forall \boldsymbol{x} \in \mathcal S\}$. The reverse polar cone $\mathcal S^*$ and the closed convex hull $\overline{\operatorname{conv}}(\mathcal S)$ of $\mathcal S$ has the following relationship:
\begin{prop}
Let $\mathcal S \subseteq \mathbb R^n$. Then, 
\[\boldsymbol{x} \in \overline{\operatorname{conv}}(\mathcal S) \Longleftrightarrow (-\boldsymbol{x},1) \in (\mathcal S^*)^{\circ}.\]
\label{prop:reverse}
\end{prop}

To test whether a point $\boldsymbol{\hat{x}}$ lies in $\overline{\operatorname{conv}}(\mathcal S)$ and, if not, produce a separating inequality, the CGP maximizes the violation $\alpha_0 - \boldsymbol{\alpha}_x^\top \boldsymbol{\hat{x}}$ over all valid inequalities, i.e., over $\mathcal S^*$. Since the feasible region $\mathcal S^*$ is a cone, this maximization is unbounded whenever $\boldsymbol{\hat{x}} \notin \overline{\operatorname{conv}}(\mathcal S)$. To ensure boundedness and obtain meaningful inequalities, it is therefore necessary to impose a normalization constraint.

We introduce a generalized normalization constraint via the support function $\sigma_{\mathcal C}$ of a nonempty, convex, and compact set $\mathcal C \subseteq \mathbb{R}^n$:
\begin{subequations}
\makeatletter
        \def\@currentlabel{\texttt{CGP}}
        \makeatother
        \label{eq:w:cgp}
        \renewcommand{\theequation}{{\texttt{CGP}}.{\textnormal{\alph{equation}}}}
\begin{align}
\sup \ & \alpha_{0} - \boldsymbol{\alpha}_{x}^{\top} \hat{x} \\
\text{s.t. } & (\boldsymbol{\alpha_x}, \alpha_0) \in {\mathcal S}^*, \label{eq:cgp:b}\\
& \sigma_{\mathcal C}(\boldsymbol{\alpha_x}) \le 1.\label{eq:cgp:c}
\end{align}
\label{prob:cgp}
\end{subequations}
Constraint \eqref{eq:cgp:c} subsumes two common normalization schemes. The first is norm-based normalization \cite{balas2002lift,kilincc2017lift}, given by
\begin{equation}
\|\boldsymbol{\alpha}_x\|_q \le 1, \quad q \in [1,\infty],
\label{eq:normalization:coeff}
\end{equation}
which corresponds to choosing $\mathcal C$ as the dual unit ball,
\(
\mathcal C = \{\boldsymbol{x} \in \mathbb{R}^n : \|\boldsymbol{x}\|_p \le 1\},\) where $p$ is such that
\(\tfrac{1}{p} + \tfrac{1}{q} = 1.
\)
The second is reverse polar normalization \cite{balas2002lift,perregaard2001generating,cornuejols1991comparison,serra2020reformulating,brandenberg2021refined}, defined as
\begin{equation}
\boldsymbol{\alpha}_x^\top \boldsymbol{y} \le 1,
\label{eq:normalization:reverse-polar}
\end{equation}
for a given vector $\boldsymbol{y}$, which corresponds to $\mathcal C$ being a singleton, $\mathcal C = \{\boldsymbol{y}\}$.

Both normalization schemes are known to yield supporting hyperplanes for appropriate choices of $\boldsymbol{y}$, and their corresponding dual formulations admit interpretations as distance computations to $\overline{\operatorname{conv}}(\mathcal S)$. The formulation \eqref{prob:cgp} recovers these results within a unified and streamlined framework in the following section.

\begin{rema}[Related formulations]
  A normalization constraint similar to \eqref{eq:cglp:c} has appeared in \cite{fullner2024new}, but without being formalized via support functions. As a result, existing analyses treat coefficient and reverse polar normalization separately, whereas our formulation provides a unified analysis.
  \end{rema}

\subsection{Dual of \eqref{prob:cgp}}
By expressing the normalization constraint through the support function in \eqref{eq:cgp:c}, we obtain a unified dual interpretation of \eqref{prob:cgp} for different normalization schemes. This dual perspective also motivates the design of a practical oracle for generating disjunctive cuts for \eqref{prob:benders-reformulation} later. We now consider the following problem, which we denote by (\texttt{DCGP}):
\begin{subequations}
\makeatletter
        \def\@currentlabel{\texttt{DCGP}}
        \makeatother
        \label{eq:w:dcgp}
        \renewcommand{\theequation}{{\texttt{DCGP}}.{\textnormal{\alph{equation}}}}
\begin{align}
\inf_{\tau,\boldsymbol{s},\boldsymbol{x},z} \ & \tau\\
\text{s.t.} \ 
              [\boldsymbol{\alpha_x}] \quad   & \boldsymbol{s} = \boldsymbol{x}  - \hat{\boldsymbol x},\label{eq:dcgp:b}\\
              [{\alpha_0}] \quad   & z = 1, \label{eq:dcgp:c}\\
              \quad   & (-\boldsymbol{x}, z) \in (\mathcal S^*)^\circ,\label{eq:dcgp:d}\\
              & \boldsymbol{s} \in \tau \mathcal C.\label{eq:dcgp:e}
\end{align}
\label{prob:dcgp}
\end{subequations}
\begin{theo}\label{theo:dcgp}
Problems \eqref{prob:cgp} and \eqref{prob:dcgp} enjoy strong duality, and \eqref{prob:dcgp} is equivalent to
\begin{equation}
\inf \gamma_{\mathcal C}(\boldsymbol x - \boldsymbol{\hat x}): \boldsymbol x \in \overline{\operatorname{conv}}(\mathcal{S}),\label{prob:dcgp:gauge}
\end{equation}
which measures the distance from $\boldsymbol{\hat{x}}$ to $\overline{\operatorname{conv}}(\mathcal{S})$ induced by the gauge function $\gamma_{\mathcal C}$. Equivalently, it computes the smallest $\tau \ge 0$ such that the translated set $\boldsymbol{\hat{x}}+\tau\mathcal C$ intersects $\overline{\operatorname{conv}}(\mathcal{S})$.
\end{theo}

\begin{rema}[Alternative assumption]
Let $\boldsymbol{p}:= (-\boldsymbol{\hat x}, 1)$, $\boldsymbol{\alpha}:=(\boldsymbol{\alpha_x},\alpha_0)$, and $\mathcal C' := \mathcal C \times \{0\}$ such that $\sigma_{\mathcal C'}(\boldsymbol{\alpha}) = \sigma_{\mathcal C}(\boldsymbol{\alpha_x})$.
If instead $\mathcal C$ is closed and convex and satisfies a standard regularity condition such as $\operatorname{ri}((\mathcal C')^\circ)\cap \operatorname{ri}(\mathcal S^*)\neq\varnothing$, then one may write
\begin{align*}
\eqref{prob:cgp}
&= \sup_{\boldsymbol{\alpha}}
\alpha_0 - \boldsymbol{\alpha}_x^\top \boldsymbol{\hat{x}} 
- \big(\delta_{(\mathcal C')^\circ} + \delta_{\mathcal S^*}\big)(\boldsymbol{\alpha}) \\
&= \big(\delta_{(\mathcal C')^\circ} + \delta_{\mathcal S^*}\big)^*(\boldsymbol p) \\
&= \big(\gamma_{\mathcal C'} \square \delta_{(\mathcal S^*)^\circ}\big)(\boldsymbol p) \\
&= \inf_{\boldsymbol s}\Big\{\gamma_{\mathcal C'}(\boldsymbol s):
\boldsymbol p-\boldsymbol s \in (\mathcal S^*)^\circ\Big\},
\end{align*}
and this gives the same result. This argument does not require compactness, but it does require the usual qualification condition to remove the closure.
\label{rema:alternative-assump}
\end{rema}

Theorem~\ref{theo:dcgp} implies that \eqref{prob:dcgp} attains a finite optimal value if and only if $\overline{\operatorname{conv}}(\mathcal{S})$ is nonempty and the translated set $\boldsymbol{\hat{x}} + \tau \mathcal C$ intersects $\overline{\operatorname{conv}}(\mathcal{S})$ for some $\tau \ge 0$; otherwise, the problem is infeasible. Moreover, it implies that \eqref{prob:cgp} yields a supporting hyperplane of $\overline{\operatorname{conv}}(\mathcal{S})$ whenever it attains a positive optimum:
\begin{prop}\label{prop:supporting}
  Suppose \eqref{prob:dcgp} attains a finite optimum at $(\hat\tau,\boldsymbol{\hat{s}},\boldsymbol{x'}, \hat z)$, and let $\boldsymbol{\hat\alpha} = (\boldsymbol{\hat{\alpha}_x}, \hat\alpha_0)$ be an optimal dual solution associated with \eqref{eq:dcgp:b}-\eqref{eq:dcgp:c}. If $\hat\tau > 0$, then
  \(
  \hat\alpha_0 - (\boldsymbol{\hat{\alpha}_x})^\top \boldsymbol{x} \le 0
  \)
  is a valid inequality for $\overline{\operatorname{conv}}(\mathcal{S})$ that separates $\boldsymbol{\hat{x}}$ and supports $\overline{\operatorname{conv}}(\mathcal{S})$ at $\boldsymbol{x}'$. If $\hat\tau = 0$, then $\boldsymbol{\hat{x}} \in \overline{\operatorname{conv}}(\mathcal{S})$.  
\end{prop}

Theorem~\ref{theo:dcgp} and Proposition~\ref{prop:supporting} recover the known geometric interpretations of the normalization constraints \eqref{eq:normalization:coeff} and \eqref{eq:normalization:reverse-polar} as special cases:
\begin{exam}[Norm-based normalization \eqref{eq:normalization:coeff}] 
  When $\mathcal C$ is the unit ball, \eqref{prob:dcgp} reduces to a projection problem that computes the point in $\overline{\operatorname{conv}}(\mathcal{S})$ closest to $\boldsymbol{\hat{x}}$ with respect to the $\ell_p$-norm. The optimal dual solution then defines a supporting hyperplane at the resulting point. In this case, the gauge in \eqref{prob:dcgp:gauge} coincides with the $\ell_p$-distance, yielding
  \begin{equation*}
    \min \ \|\boldsymbol{x} - \boldsymbol{\hat{x}}\|_p
    \quad \text{s.t.} \quad \boldsymbol{x} \in \overline{\operatorname{conv}}(\mathcal{S}).
  \end{equation*}
  This recovers the results in \cite{boyd1995convergence,cadoux2010computing}. Constraint \eqref{eq:dcgp:e} is expressed as $\tau \ge \|\boldsymbol{s}\|_p$. 
  \end{exam}
  
  \begin{exam}[Reverse polar normalization \eqref{eq:normalization:reverse-polar}]
  When $\mathcal C$ is a singleton, \eqref{prob:dcgp} reduces to finding the smallest $\tau \ge 0$ such that
  \[
  \boldsymbol{\hat{x}} + \tau \boldsymbol{y}
  \in \overline{\operatorname{conv}}(\mathcal S),
  \]
  and the optimal dual solution defines a supporting hyperplane at the resulting point. This recovers the results in \cite{balas2002lift,cornuejols2006convex}. In this case, \eqref{eq:dcgp:e} reduces to $\boldsymbol{s} = \tau \boldsymbol{y}$. 
  \end{exam}

  \begin{rema}[Design of new normalization constraints]
    Note that \eqref{eq:cgp:c} is equivalent to $\boldsymbol{\alpha_x} \in \mathcal C^\circ$. Thus, by choosing different nonempty, convex, and compact sets $\mathcal C$, one can derive customized normalization constraints via their polars while preserving all preceding derivations.
    
    For example, if $\mathcal C$ is the convex hull of vectors $\{\boldsymbol{y^{i}}\}_{i=1}^k$, then
    \(
    \mathcal C^\circ = \{\boldsymbol{\alpha_x} : \boldsymbol{\alpha_x}^\top \boldsymbol{y^i} \le 1,\ i=1,\dots,k\},
    \)
    yielding a normalization constraint defined by multiple reverse polar constraints. This extension may be useful, for instance, when incorporating a new (potentially invalid) direction $\boldsymbol{\hat{y}}$ together with known valid core directions. In this case, \eqref{eq:dcgp:e} can be written as $\boldsymbol{s} = \sum_{i=1}^k \lambda_i\boldsymbol{y^i}$, $\sum_{i=1}^k \lambda_i= \tau$, and $\lambda_i \ge 0, \ \forall i=1,\cdots,k$. Similarly, knapsack-based sets $\mathcal C = \{\boldsymbol{x} \ge 0 : \boldsymbol{w}^\top \boldsymbol{x} \le 1\}$ for some $\boldsymbol{w} \in \mathbb{R}_{>0}^n$ lead to componentwise asymmetric bounds on $\boldsymbol{\alpha_x}$. In this case, \eqref{eq:dcgp:e} is equivalently expressed as $\boldsymbol{s} \ge 0$ and $\boldsymbol{w}^\top \boldsymbol{s} \le \tau$.
    \end{rema}

\begin{rema}[A more general setting] The normalization constraint \eqref{eq:cgp:c} can be generalized to $\sigma_{\mathcal C'}(\boldsymbol{\alpha}) \le 1$ for any nonempty, convex, and compact set $\mathcal C' \subseteq \mathbb{R}^{n+1}$, and the approach proposed in the next section extends naturally to this formulation. Following the same argument as in the proof of Theorem~\ref{theo:dcgp} or Remark~\ref{rema:alternative-assump}, one obtains the corresponding dual: 
\begin{align}\label{eq:generalized-cgp} \eqref{prob:cgp} &= \inf_{(\boldsymbol{s_x}, s_0)} \gamma_{\mathcal C'}(\boldsymbol{s_x}, s_0) :\, (-\boldsymbol{\hat x} - \boldsymbol{s_x},\, 1 - s_0) \in (\mathcal S^*)^\circ \\
&= \inf_{(\boldsymbol{x}, s_0)} \gamma_{\mathcal C'}(\boldsymbol{x}-\boldsymbol{\hat x}, s_0) :\, \boldsymbol{x} \in (1 - s_0)\overline{\operatorname{conv}}(\mathcal S),\nonumber
\end{align} 
where the second equality follows from Proposition~\ref{prop:reverse}. 
This introduces to \eqref{prob:dcgp} a nonnegative variable $s_0$, replaces the right-hand side of \eqref{eq:dcgp:c} with $1 - s_0$, and modifies \eqref{eq:dcgp:e} to $\boldsymbol{s} \in \tau \mathcal C'$. 

This includes normalization on $\alpha_0$ as a special case, the one considered in \cite{balas1996mixed, conforti2019facet}:
\begin{equation*}
|\alpha_0| \le 1,
\end{equation*} 
which corresponds to $\mathcal C'$ being $\boldsymbol{0} \times \{z : |z| \le 1\}$. 
In this case, $\gamma_{\mathcal C'}(\boldsymbol{s_x}, s_0) = |s_0|,$ if $\boldsymbol{s_x} = \boldsymbol{0}$, and $\infty$, otherwise.
Hence, 
\begin{align*}
\eqref{eq:generalized-cgp}=\inf_{s_0} |s_0| : (-\boldsymbol{\hat x}, 1 - s_0) \in (\mathcal S^*)^\circ
= \inf_{s_0} |s_0| : \boldsymbol{\hat x} \in (1 - s_0)\overline{\operatorname{conv}}(\mathcal S).
\end{align*}
This coincides with the dual studied in \cite{conforti2019facet}.
\end{rema}
}

\section{Proposed approach: disjunctive Benders decomposition}\label{sec:dbd}

Since Balas's seminal work \cite{balas1979disjunctive}, disjunctive cuts have been extensively used to solve mixed-integer programs in both linear and nonlinear convex settings \cite{balas2018disjunctive}. Notably, disjunctive cuts—particularly split cuts, which are derived from two-term disjunctions—subsume several prominent classes of cutting planes, including Chvátal-Gomory cuts, Gomory Mixed-Integer cuts, and Mixed-Integer Rounding cuts \cite{balas2003precise,conforti2014integer}. 
Disjunctions are useful for conveying integrality information in the cut-generation process. In this section, we propose a framework that strengthens BD by leveraging disjunctions {\color{black} and \eqref{prob:cgp}.}

In what follows, $\mathcal{P}$ denotes the continuous relaxation of the feasible region of \eqref{prob:benders-reformulation}, i.e., $\mathcal{P} = (\mathcal{X}_{LP} \times \mathbb{R}) \cap \mathcal{L}$ where $\mathcal{X}_{LP} = \{\boldsymbol{x} \in \mathbb{R}^{n_x}: \boldsymbol{\underline x} \leq \boldsymbol{x} \leq \boldsymbol{\overline x}, D \boldsymbol{x} \geq \boldsymbol{h}\}$.

\subsection{Disjunction of the Benders reformulation }\label{sec:dbd:disjunction}
{\color{black}
Consider a disjunction of the form 
\begin{equation}
    \bigvee_{r \in \mathcal R} U^r \boldsymbol{x} \ge \boldsymbol{u^r},
    \label{eq:disj}
\end{equation}
where $\mathcal R$ is a finite index set, and for each $r \in \mathcal R$, $U^r$ and $\boldsymbol{u^r}$ are a real matrix and a real vector of appropriate dimensions. The disjunction is constructed by exploiting the integrality of $\boldsymbol x$ while ensuring that no integer solution is excluded. As an example, a two-term disjunction (i.e., $|\mathcal R| = 2$) can be derived from any integral vector $(\boldsymbol{\phi}, \phi_0) \in \mathbb Z^{n_x}\times \mathbb Z$ satisfying $\phi_j = 0$ for all $j \in [n_x] \setminus \mathcal{I}$, which ensures that $\boldsymbol{\phi}^\top \boldsymbol{x} \in \mathbb{Z}$ for all $\boldsymbol{x} \in \mathcal{X}$. Then, the disjunction $\boldsymbol{\phi}^{\top} \boldsymbol{x} \leq \phi_0 \lor \boldsymbol{\phi}^{\top} \boldsymbol{x} \geq \phi_0 + 1$, defined by $U^1 = -\boldsymbol{\phi}^\top$, $u^1 = -\phi_0$, $U^2 = \boldsymbol{\phi}^\top$, and $u^2 = \phi_0 + 1$, does not exclude any integral point.}
The integral vector $(\boldsymbol{\phi}, \phi_0)$ is called a \textit{split} for \eqref{prob:benders-reformulation}. 

{\color{black}
Using the disjunction \eqref{eq:disj}, a union of polyhedra $\mathcal{P}_{\mathcal R} := \bigcup_{r \in \mathcal R} \mathcal{P}_{r}$ can be constructed, where
\begin{align*}
\mathcal{P}_r & := \mathcal{P} \cap \{(\boldsymbol{x}, t)\in \mathbb R^{n_x} \times \mathbb R : U^r \boldsymbol{x} \geq \boldsymbol{u^r}\}\\
    & 
    {\begin{array}{rl}
    =\{(\boldsymbol{x},t) \in \mathbb R^{n_x +1 }: [\lambda_{\hat{\pi}}^r] \quad & t \geq \boldsymbol{\hat{\pi}}^{\top} (\boldsymbol{b}-A\boldsymbol{x}), \ \forall \boldsymbol{\hat{\pi}} \in \mathcal{J},\\
    \quad [\mu_{{\tilde{\pi}}}^r]  \quad & 0 \geq \boldsymbol{\tilde{\pi}}^{\top} (\boldsymbol{b}-A\boldsymbol{x}), \ \forall \boldsymbol{\tilde{\pi}} \in \mathcal{R},\\
    \quad [\boldsymbol{\theta^r}] \quad & D \boldsymbol{x} \geq \boldsymbol{h},\\
    \quad [\boldsymbol{\kappa^r}]  \quad & U^r \boldsymbol{x} \ge \boldsymbol{u^r},\\
    \quad [\boldsymbol{\eta^r}] \quad & -\boldsymbol{x} \geq -\boldsymbol{\overline{x}}, \\
    \quad [\boldsymbol{\nu^r}] \quad & \boldsymbol{x} \geq \boldsymbol{\underline{x}}\},
    \end{array}}
    \end{align*}}
where the variables in square brackets denote the Farkas multipliers associated with each constraint. 

Note that $\overline{\operatorname{conv}}\left((\mathcal{X} \times \mathbb{R}) \cap \mathcal{L}\right) \subseteq \overline{\operatorname{conv}}{\color{black}(\mathcal{P}_{\mathcal R})} \subseteq \mathcal{P} \subseteq \mathcal{L}$, as the {\color{black}disjunction} conveys partial integrality information. 
Moreover, when all integer variables are binaries (i.e., $[\underline{x}_j, \overline{x}_j] = [0, 1]$ for all $j \in \mathcal{I}=\{1,\cdots, \ell\}$), the following holds from the sequential convexification theorem \cite{balas1998disjunctive}:  
\begin{equation*}
\overline{\operatorname{conv}}\left((\mathcal{X} \times \mathbb{R}) \cap \mathcal{L}\right) = \overline{\operatorname{conv}}(\mathcal{P}_{\mathcal R_\ell} \cap \overline{\operatorname{conv}}(\cdots \cap \overline{\operatorname{conv}}(\mathcal{P}_{\mathcal R_2} \cap \overline{\operatorname{conv}}(\mathcal{P}_{\mathcal R_1})) \cdots )),
\end{equation*} 
where $\mathcal R_j$ denote the disjunction $x_j \le 0 \lor x_j \ge 1$ for $j \in \mathcal I$.
This implies that $\overline{\operatorname{conv}}((\mathcal{X} \times \mathbb{R}) \cap \mathcal{L})$ can be obtained by iteratively adding valid inequalities for simple disjunctive systems, $\overline{\operatorname{conv}}(\mathcal{P}_{\mathcal R_j})$ for $j \in \mathcal{I}$, in a specialized augmented manner \cite{balas1998disjunctive}. {\color{black}The sequential convexification result forms the foundation of a specific implementation of the proposed framework that eliminates the need to solve the master problem as a MILP in BD.}

\subsection{Cut-generating linear program}\label{sec:cglp}

{\color{black}Typically, disjunctive cuts are generated by solving \eqref{prob:cgp} where $\mathcal S = \mathcal P_{\mathcal R}$, originally introduced by Balas \cite{balas1979disjunctive}.}

\paragraph{Valid inequalities for \(\mathcal P_{\mathcal R}\)}\label{sec:cglp:valid}
{\color{black}The reverse polar cone $\mathcal{P}_{\mathcal{R}}^*$ of $\mathcal{P}_{\mathcal{R}}$ can be expressed using the Farkas multipliers in an extended space.
Assume $\mathcal{P}_r \neq \emptyset, \ \forall r \in \mathcal R$. An inequality $\hat\alpha_0 - (\boldsymbol{\hat\alpha_x})^{\top} \boldsymbol{x} - \hat\alpha_t t \leq 0$ is a valid inequality for $\mathcal P_{\mathcal R} = \bigcup_{r \in \mathcal R}\mathcal{P}_r$}, if and only if $(\boldsymbol{\hat\alpha_x}, \hat\alpha_t,\hat\alpha_0)$ meets the following constraints (see, e.g., \cite[Theorem 1.2]{balas2018disjunctive}):
{\color{black}
\begin{subequations}
\begin{align}
\ [\omega_t^r] \quad   & \alpha_{t} = \sum_{\boldsymbol{\hat{\pi}} \in \mathcal{J}} \lambda_{{\hat{\pi}}}^r, & \forall r \in \mathcal R,\\
              [\boldsymbol{\omega_x^r}] \quad   & \boldsymbol{\alpha_{x}} = \sum_{\boldsymbol{\tilde{\pi}} \in \mathcal{R}}  (A^\top \boldsymbol{\tilde{\pi}}) \mu_{{\tilde{\pi}}}^r 
              + \sum_{\boldsymbol{\hat{\pi}} \in \mathcal{J}}  (A^\top  \boldsymbol{\hat{\pi}})\lambda_{{\hat{\pi}}}^r+ D^{\top} \boldsymbol{\theta^r} + (U^r)^\top\boldsymbol{\kappa^r} - \boldsymbol{\eta^r} + \boldsymbol{\nu^r}, & \forall r \in \mathcal R,\label{eq:cglp:c}\\
              [\omega_0^r] \quad   & \alpha_{0} \leq \sum_{\boldsymbol{\tilde{\pi}} \in \mathcal{R}}  (\boldsymbol{b}^\top \boldsymbol{\tilde{\pi}})\mu_{{\tilde{\pi}}}^r
              + \sum_{\boldsymbol{\hat{\pi}} \in \mathcal{J}}  (\boldsymbol{b}^\top  \boldsymbol{\hat{\pi}})\lambda_{{\hat{\pi}}}^r 
              + \boldsymbol{h}^{\top} \boldsymbol{\theta^r}
              + (\boldsymbol{u^r})^\top \boldsymbol{\kappa^r}
              - \boldsymbol{\overline x}^{\top} \boldsymbol{\eta^r}
              + \boldsymbol{\underline x}^{\top} \boldsymbol{\nu^r}, & \forall r \in \mathcal R,\\
              & (\lambda^r_{{\hat{\pi}}})_{\boldsymbol{\hat{\pi}} \in \mathcal{J}}, (\mu_{{\tilde{\pi}}}^r)_{\boldsymbol{\tilde{\pi}} \in \mathcal{R}}, \boldsymbol{\theta^r}, \boldsymbol{\kappa^r}, \boldsymbol{\eta^r}, \boldsymbol{\nu^r} \geq 0,  & \forall r \in \mathcal R,
    \end{align}\label{eqs:valid-cut}
\end{subequations}
}
where the variables in square brackets are the associated dual variables. We let $\boldsymbol{\alpha} := (\boldsymbol{\alpha_x}, \alpha_t, \alpha_0)$ and let $\boldsymbol{\chi} := ((\lambda^r_{{\hat{\pi}}})_{\boldsymbol{\hat{\pi}} \in \mathcal{J}}, (\mu_{{\tilde{\pi}}}^r)_{\boldsymbol{\tilde{\pi}} \in \mathcal{R}}, \boldsymbol{\theta^r}, \boldsymbol{\kappa^r}, \boldsymbol{\eta^r},\boldsymbol{\nu^r})_{\color{black}r \in \mathcal R}$. 
Then, 
\begin{equation*}
\mathcal{P}_{\mathcal R}^* := \{ \boldsymbol{\alpha} : \exists \boldsymbol{\chi} \text{ s.t. } (\boldsymbol{\alpha}, \boldsymbol{\chi}) \text{ meets } \eqref{eqs:valid-cut}\}.
\end{equation*}

\paragraph{Normalization constraint} 
The choice of normalization constraint plays a critical role in generating effective disjunctive cuts, as widely discussed in the literature \cite{balas2002lift,fischetti2011separation,lodi2023disjunctive}. A commonly used normalization approach is referred to as standard normalization, which imposes $\|\boldsymbol\chi\|_1 = 1$. However, as noted in \cite{fischetti2011separation}, this method is susceptible to issues such as redundancy in constraints and sensitivity to constraint scaling. Specifically, even a weak cut may appear to be the most violated one in the presence of redundant constraints or when constraints are improperly scaled; a cut generated with standard normalization is not necessarily a supporting hyperplane to $\overline{\operatorname{conv}}(\mathcal{P}_{\mathcal{R}})$. This issue is particularly problematic in the context of BD, where numerous constraints may remain unknown, making it challenging to manage redundancy as the problem is augmented with disjunctive cuts. Additionally, feasibility cuts, defined by rays, can be arbitrarily scaled. 

{\color{black}
A more robust approach is to normalize in the projected $\boldsymbol{\alpha}$-space. Two common choices are the norm-based normalization \eqref{eq:normalization:coeff} and the reverse polar normalization \eqref{eq:normalization:reverse-polar}. As discussed in \cite{balas2002lift}, when $\overline{\operatorname{conv}}(\mathcal{P}_{\mathcal R})$ is full-dimensional, facet-defining inequalities correspond to extreme rays of $\mathcal{P}_{\mathcal R}^*$. For $q \in \{1,\infty\}$, coefficient normalization can be interpreted as intersecting the cone $\mathcal{P}_{\mathcal R}^*$ with multiple hyperplanes, yielding an extreme point that does not necessarily correspond to an extreme ray. In contrast, reverse polar normalization uses a single hyperplane. Consequently, it guarantees the existence of an optimal extreme point corresponding to an extreme ray of $\mathcal{P}_{\mathcal R}^*$ \cite{balas2002lift}, although not every extreme point of the extended cut-generating program corresponds to an extreme ray. Conditions under which reverse polar normalization yields facet-defining inequalities are given in \cite{cornuejols2006convex,brandenberg2021refined}.

}

\paragraph{Dual of the \eqref{prob:cgp}}
{\color{black}
First, note that
\begin{align*}
(\mathcal{P}_{\mathcal R}^*)^\circ
&:= \{\boldsymbol{v} : \boldsymbol{v}^\top \boldsymbol{\alpha} \le 0, \ \forall \boldsymbol{\alpha} \in \mathcal{P}_{\mathcal R}^*\} \\
&= \left\{\left(\textstyle-\sum_{r \in \mathcal R}\boldsymbol{\omega_x^r},\ -\sum_{r \in \mathcal R}\omega_t^r,\ \sum_{r \in \mathcal R}\omega_0^r\right)
:\ (\boldsymbol{\omega_x^r}, \omega_t^r, \omega_0^r) \in (\mathcal{P}_r)^\#, \ \omega_0^r \ge 0,\ \forall r \in \mathcal R \right\},
\end{align*}
where the last equality follows from strong duality of the feasible linear program $\sup_{\boldsymbol{\alpha} \in \mathcal{P}_{\mathcal R}^*} \boldsymbol{v}^\top \boldsymbol{\alpha}$. Here, $\mathcal{Q}^\#$ denotes a {\color{black}homogenized} representation of a polyhedron $\mathcal{Q} := \{\boldsymbol{u} \in \mathbb{R}^n : A\boldsymbol{u} \leq \boldsymbol{b}\}$ in a space of one higher dimension, obtained by scaling all the right-hand side expressions with a new nonnegative variable; specifically,  
$\mathcal{Q}^\# = \{(\boldsymbol{u}, v) \in \mathbb{R}^{n} \times \mathbb{R}_{\geq 0} : A\boldsymbol{u} \leq \boldsymbol{b}v\}$.
For instance, $\mathcal{L}^{\#} = \{(\boldsymbol{x},t,z) \in \mathbb{R}^{n_x+1} \times \mathbb{R}_{\geq 0}: t \geq \boldsymbol{\hat{\pi}}^{\top} (\boldsymbol{b}z - A\boldsymbol{x}), \forall \boldsymbol{\hat{\pi}} \in \mathcal{J}, 0 \geq \boldsymbol{\tilde{\pi}}^{\top} (\boldsymbol{b}z - A\boldsymbol{x}), \forall \boldsymbol{\tilde{\pi}} \in \mathcal{R}\}$ and $\mathcal{X}_{LP}^{\#}=\{(\boldsymbol{x},z) \in \mathbb{R}^{n_x} \times \mathbb{R}_{\geq 0}: D \boldsymbol{x} \geq \boldsymbol{h}z, \boldsymbol{x} \geq \boldsymbol{\underline x}z, -\boldsymbol{x} \geq -\boldsymbol{\overline x} z\}$. {\color{black}Therefore, $(\mathcal P_r)^\# = \{(\boldsymbol{x},t,z): (\boldsymbol{x},t,z) \in \mathcal L^\#, (\boldsymbol{x},z) \in \mathcal X_{LP}^\#, U^r \boldsymbol{x} \ge \boldsymbol{u^r} z\}$.}

Using this, we derive \eqref{prob:dcgp} for $\mathcal S = \mathcal P_{\mathcal R}$: 
\begin{subequations}
\begin{align}
\inf \ & \tau  \\
\text{s.t.} \ 
& (\boldsymbol{\omega_x^r}, \omega_t^r, \omega_0^r) \in {(\mathcal{P}_r)^\#},\quad \forall r \in \mathcal R,\label{eq:dcglp:b}\\
              & \omega^r_0 \geq 0, \;\;\;\qquad\qquad\qquad \forall r \in \mathcal R,\label{eq:dcglp:c}\\
              [\alpha_0] \quad   & \sum_{r \in \mathcal R}\omega_0^r = 1, \label{eq:dcglp:d}\\
              [\boldsymbol{\alpha_x}] \quad   & \boldsymbol{s_x} = \sum_{r \in \mathcal R}\boldsymbol{\omega_x^r}  - \hat{\boldsymbol x},\label{eq:dcglp:e}\\
              [\alpha_t] \quad   & s_t = \sum_{r \in \mathcal R}\omega_t^r - \hat{t},\label{eq:dcglp:f}\\
              & \boldsymbol{s} \in \tau \mathcal C,\label{eq:dcglp:g}\\
              & \tau \ge 0.\label{eq:dcglp:h}
\end{align}
\label{prob:dcglp}
\end{subequations}}
\setcounter{equation}{3} 
{\color{black}
The following corollary immediately follows from Theorem \ref{theo:dcgp} and Proposition \ref{prop:supporting}:
\begin{coro}\label{coro:dcglp}
Problem \eqref{prob:dcglp} is equivalent to
\begin{equation*}
\inf \gamma_{\mathcal C}((\boldsymbol x, t) - (\boldsymbol{\hat x}, \hat{t})): (\boldsymbol x, t) \in \overline{\operatorname{conv}}(\mathcal{P}_{\color{black}\mathcal{R}}),
\end{equation*}
i.e., the problem attains a finite optimum if and only if the translated set $(\boldsymbol{\hat{x}}, \hat{t})+\tau\mathcal C$ intersects $\overline{\operatorname{conv}}(\mathcal{P}_{\mathcal R})$ for some $\tau > 0$; otherwise, it is infeasible. Suppose Problem \eqref{prob:dcglp} attains a finite optimal value at $(\hat\tau, (\boldsymbol{\hat\omega^r})_{r \in \mathcal R}, \boldsymbol{\hat s})$. Let $\boldsymbol{\hat\alpha} = (\boldsymbol{\hat\alpha_x}, \hat\alpha_t, \hat\alpha_0)$ be the optimal dual solution associated with \eqref{eq:dcglp:d}--\eqref{eq:dcglp:f}. If $\hat\tau > 0$, 
\(
\hat\alpha_0 - \boldsymbol{\hat\alpha_x}^\top \boldsymbol{x} - \hat\alpha_t t \le 0
\)
defines a valid inequality for $\overline{\operatorname{conv}}(\mathcal{P}_{\mathcal R})$ that separates $(\boldsymbol{\hat x}, \hat t)$ and supports $\overline{\operatorname{conv}}(\mathcal{P}_{\mathcal R})$ at 
\[\textstyle
(\boldsymbol{x}', t') := \left(\sum_{r \in \mathcal R} \boldsymbol{\hat\omega_x^r}, \sum_{r \in \mathcal R} \hat\omega_t^{r}\right) \in \overline{\operatorname{conv}}(\mathcal{P}_{\mathcal R}).
\]
If $\hat\tau = 0$, then $(\boldsymbol{\hat{x}}, \hat{t}) \in \overline{\operatorname{conv}}(\mathcal{P}_{\mathcal R})$.
\end{coro}

\begin{rema}[Interpretation of the optimal solution of \eqref{prob:dcglp}]\label{rema:dcglp-sol-interpretation}
Let $(\hat\tau, (\boldsymbol{\hat\omega^r})_{r \in \mathcal R},\boldsymbol{\hat{s}})$ (and $\boldsymbol{\hat\alpha}$) be the optimal primal (and dual) solution of \eqref{prob:dcglp}. 
For each $r \in \mathcal R$, define
$$(\boldsymbol{x^{r\prime}},t^{r\prime}):=\begin{cases}
  (\boldsymbol{\hat\omega^r_x}/\hat\omega_0^r,\hat\omega_t^{r}/\hat\omega_0^r) \in \mathcal{P}_r, & \mbox{ if } \hat\omega_0^r>0,\\
  (\boldsymbol{\hat\omega^r_x},\hat\omega_t^{r}), & \mbox{ o.w.}
\end{cases}$$
and $(\boldsymbol{x'}, t') := \sum_{r \in \mathcal R} \hat\omega_0^r(\boldsymbol{x^{r\prime}},t^{r\prime}) = (\sum_{r \in \mathcal R} \boldsymbol{\hat\omega_x^r}, \sum_{r \in \mathcal R} \hat\omega_t^{r})$.

Note $(\boldsymbol{x'}, t')$ is a convex combination of $(\boldsymbol{x^{r\prime}}, t^{r\prime}) \in \mathcal P_r$ due to \eqref{eq:dcglp:c}-\eqref{eq:dcglp:d}, and from Corollary \ref{coro:dcglp} it is the projection of $(\boldsymbol{\hat{x}},\hat{t})$ onto $\overline{\operatorname{conv}}(\mathcal{P}_{\mathcal{R}})$. The projection point $(\boldsymbol{x'}, t')$ is the intersection of the valid inequality $\boldsymbol{\hat\alpha_x}^{\top} \boldsymbol{x} + \hat\alpha_t t \geq \hat\alpha_0$ with $\overline{\operatorname{conv}}(\mathcal{P}_{\mathcal{R}})$.
\end{rema}
}

\subsection{\textsc{{\color{black}disjunctiveOracle}}}\label{sec:dbd:oracle}
Corollary \ref{coro:dcglp} and Remark \ref{rema:dcglp-sol-interpretation} offer valuable insights and properties of \eqref{prob:dcglp}. Building on these results, we propose an oracle for generating disjunctive cuts for the Benders reformulation \eqref{prob:benders-reformulation}, which we refer to as \emph{disjunctive Benders cuts}. Within this oracle, \eqref{prob:dcglp} is solved via a cutting plane algorithm, and then the simple retrieval of (optimal) dual solutions associated with \eqref{eq:dcglp:d}–\eqref{eq:dcglp:f} generates (the deepest) disjunctive Benders cuts.

Note that \eqref{prob:dcglp} involves $\mathcal{L}^\#$, which is often defined by an exponential number of constraints associated with $\mathcal{J}$ and $\mathcal{R}$. However, $\mathcal{L}^\#$ has a key property, outlined in the following proposition, which allows us to generate separating hyperplanes of $\mathcal{L}^\#$ using typical Benders oracles (see Figure \ref{fig:oracle-typical} for the definition of typical Benders oracles), without the need to have all the constraints upfront:
\begin{prop}
Let $(\boldsymbol{\hat\omega_x}, \hat\omega_t, \hat\omega_0) \in \mathbb{R}^{n_x+1} \times \mathbb{R}_{> 0}$. Suppose there exists a hyperplane of the form $\boldsymbol{\beta_x}^{\top} \boldsymbol{x} + \beta_t t \geq \beta_0$ that separates $(\boldsymbol{x'}, t'):=(\boldsymbol{\hat\omega_x}/\hat\omega_0, \hat\omega_t/\hat\omega_0)$ from $\mathcal{L}$. Then, the hyperplane $\boldsymbol{\beta_x}^{\top} \boldsymbol{\omega_x} + \beta_t \omega_t \geq \beta_0 \omega_0$ separates $(\boldsymbol{\hat\omega_x}, \hat\omega_t, \hat\omega_0)$ from $\mathcal{L}^\#$. Conversely, if $(\boldsymbol{x'},t') \in \mathcal{L}$, meaning no such separating hyperplane exists, then $(\boldsymbol{\hat\omega_x},\hat\omega_t, \hat\omega_0) \in \mathcal{L}^\#$. 
\label{prop:algorithm-cut}
\end{prop}
Proposition \ref{prop:algorithm-cut} motivates Routine \ref{oracle} for generating disjunctive Benders cuts. The routine iteratively refines a relaxation \texttt{R} of \eqref{prob:dcglp}, starting with an initial approximation of $\mathcal{L}$ that includes at least one optimality cut or a trivial lower bound on $t$ (e.g., $t \geq -10^{99}$). This guarantees that the relaxation of $\mathcal{L}^\#$ contains at least one constraint of the form $t \geq \boldsymbol{\hat{\pi}}^{\top} (\boldsymbol{b}z - A\boldsymbol{x})$ for some $\boldsymbol{\hat{\pi}} \in \mathcal{J}$ or a trivial lower bound like $t \geq -10^{99}z$. 
The purpose of this construction is to guarantee that $\boldsymbol{\hat\omega^r} = (\boldsymbol{\hat\omega_x^r},\hat\omega_t^{r},\hat\omega_0^r)\in \mathcal{L}^{\#}$ whenever $\hat\omega^r_0=0$ at Line \ref{algo:master-sol} of Routine \ref{oracle}. Specifically, consider the case where $\hat\omega_0^r=0$. 
Then, it follows immediately from $\mathcal{X}_{LP}^{\#}$ in \texttt{R} that $\boldsymbol{\hat\omega_x^r}= \boldsymbol{0}$. Furthermore, the structure of the relaxation of $\mathcal{L}^{\#}$ ensures that $\hat\omega_t^{r} \geq 0$. Hence, we conclude that $\boldsymbol{\hat\omega^r} \in \mathcal{L}^{\#}$, validating Line \ref{algo:w0}.
When $ \hat\omega_0^r > 0$, based on Proposition \ref{prop:algorithm-cut}, the routine refines \texttt{R} by iteratively adding separating hyperplanes found via a \textsc{typicalOracle}, which returns a boolean variable $\mathbb{I}$ and a vector $\boldsymbol{\beta} = (\boldsymbol{\beta_x}, \beta_t, \beta_0)$ in Line \ref{algo:typical-oracle}. If $\mathbb{I}=\texttt{true}$, indicating that $(\boldsymbol{x'},t')$ belongs to $\mathcal{L}$, then $\boldsymbol{\beta} = \boldsymbol{0}$, representing a trivial valid inequality. Otherwise, $\boldsymbol{\beta}$ defines $\boldsymbol{\beta_x}^{\top} x + \beta_t t \geq \beta_0$ that separates $(\boldsymbol{x'},t')$ from $\mathcal{L}$.

The convergence behavior of Routine \ref{oracle} depends on the choice of \textsc{typicalOracle} in Line \ref{algo:typical-oracle}. However, as long as the chosen \textsc{typicalOracle} guarantees finite convergence of \textsc{BendersSeq} for solving \eqref{prob:benders-reformulation}, it likewise ensures the finite convergence of Routine \ref{oracle}. In addition, when $\textsc{typicalOracle}(\boldsymbol{x'},t')$ provides $f(\boldsymbol{x'})$ as a byproduct—such as in the case of the conventional Benders cut generation oracle, the termination can be based on optimality gap (Remark \ref{rema:opt-gap-termination}). More interestingly, even if Routine \ref{oracle} is terminated prematurely, the dual solution will still yield a valid disjunctive Benders cut that supports a relaxation of $\overline{\operatorname{conv}}(\mathcal{P}_{\mathcal{R}})$ (Proposition \ref{prop:inexact}). Moreover, Routine \ref{oracle} generates a series of Benders cuts as a byproduct (Remark \ref{rema:byproduct-benders}).
Finally, Routine \ref{oracle} can easily incorporate previously identified disjunctive cuts to generate {\color{black}a higher rank cut }(Remark \ref{rema:dcglp-augmentation}).

\begin{algorithm}[t!]
  \caption{\textsc{{\color{black}disjunctiveOracle}}: disjunctive Benders cut separation oracle}\label{oracle}
  \begin{algorithmic}[1]
    \Require tolerance $\epsilon > 0$; {\color{black}disjunction $\bigvee_{r \in \mathcal R} U^r \boldsymbol{x} \ge \boldsymbol{u^r}$}; point $(\boldsymbol{\hat{x}}, \hat{t})$ to separate
    \State $\texttt{R}$ $\gets$ a relaxation of \eqref{prob:dcglp} with $\mathcal{L}$ replaced by a relaxation of it, including at least one optimality cut or a trivial lower bound on $t$; 
    $\mathbb{I}^r_{\mathcal L^\#} \gets \texttt{false}, \ \forall r \in \mathcal R$\label{algo:initialization}
    \While{$\texttt{!}(\land_{r \in \mathcal R}\mathbb{I}^r_{\mathcal{L}^\#}$)}\label{algo:while:begin}
      \State Solve $\texttt{R}$;  $(\hat\tau, (\boldsymbol{\hat\omega^r})_{r \in \mathcal R}, \boldsymbol{\hat s}) \gets$ its optimal solution \label{algo:master-sol}
      \For{$r$ in $\mathcal R$} in parallel
      \If{$\hat\omega_0^r > \epsilon$}
      \State $\boldsymbol{x'} \gets \boldsymbol{\hat\omega_x^r} / \hat\omega^r_0;$ $t' \gets \hat\omega^{r}_t / \hat\omega^r_0$
      \State ($\mathbb{I}^r_{\mathcal{L}^\#}$, ($\boldsymbol{\beta_x}, \beta_t, \beta_0$)) $\gets$ \textsc{typicalOracle}($\boldsymbol{x'}, t'$) \label{algo:typical-oracle}
      \If{\texttt{!}$\mathbb{I}_{\mathcal{L}^\#}^r$}
      \State Add $\beta_0\omega_0 - \boldsymbol{\beta_x}^{\top} \boldsymbol{\omega^r_x} - \beta_t \omega_t^r \leq 0$ for all $r \in \mathcal R$ to $\texttt{R}$
      \EndIf
      \Else
      \State $\mathbb{I}_{\mathcal{L}^\#}^r$ $\gets \texttt{true}$\label{algo:w0}
      \EndIf
      \EndFor
      \EndWhile\label{algo:while:end}
      \State $\hat\alpha_{0}, \boldsymbol{\hat\alpha_{x}}, \hat\alpha_{t} \gets$ dual values of \eqref{eq:dcglp:d}, \eqref{eq:dcglp:e} and \eqref{eq:dcglp:f}\\
      \Return $(\hat\tau, (\boldsymbol{\hat\alpha_{x}},\hat\alpha_{t},\hat\alpha_{0}))$\label{oracle:return}
  \end{algorithmic}
\end{algorithm}

\begin{rema}[Optimality gap-based termination of Routine \ref{oracle}]\label{rema:opt-gap-termination}
  Since \texttt{R} is a relaxation of \eqref{prob:dcglp}, its optimal objective value at each iteration provides a lower bound on the optimal objective value of \eqref{prob:dcglp}. Additionally, if $\textsc{typicalOracle}(\boldsymbol{x'},t')$ provides $f(\boldsymbol{x'})$ as a byproduct—such as in the case of a conventional Benders cut generation oracle—an upper bound on the optimal objective can be computed at each iteration, allowing the routine to terminate based on the optimality gap. Specifically, after Line 11, we define
  $$\omega_t^{r\prime} := \begin{cases}
    \hat\omega_t^{r}, & \mbox{ if } \mathbb{I}^r_{\mathcal{L}^\#}=\texttt{true}\\
    \hat\omega_0^r f\left(\boldsymbol{\hat\omega^r_x}/\hat\omega^r_0\right) & \mbox{ o.w.}
  \end{cases}$$
where we use multiplication with $\infty$ as $\infty$. Then, $(\boldsymbol{\hat\omega^r_x}, \omega_t^{r\prime})$ is feasible to \eqref{prob:dcglp}, so we can compute the upper bound by: 
  {\color{black}$
  \texttt{UB} \gets \min\{\texttt{UB},\gamma_{\mathcal C}(\boldsymbol{\hat s_x}, \sum_{r \in \mathcal R}\omega_t^{r\prime} - \hat{t}) \}.
  $}
\end{rema}
  
Let $\rho \in \mathbb{Z}_{>0}$ be the iteration counter of Routine \ref{oracle}, and let $\mathcal{L}^{(\rho)} \subseteq \mathcal{L}$ denote the refined relaxation of $\mathcal{L}$ after $\rho$ iterations, that is
\begin{equation*}
\mathcal{L}^{(\rho)} = \mathcal{L}^{(0)} \cap \left\{(\boldsymbol{x},t) \in \mathbb{R}^{n_x+1} : (\boldsymbol{\beta_x^{(\rho')}})^{\top} x + \beta_t^{(\rho')} t \geq \beta_0^{(\rho')}, \quad \rho' \in \{0,1,\ldots, |\mathcal R|\rho-1\} \right\},
\end{equation*}
where $\mathcal{L}^{(0)}$ is the initial relaxation. For each $\rho' \in [\rho]$, the vectors $\boldsymbol{\beta}^{(|\mathcal R|(\rho'-1))}, \cdots, \boldsymbol{\beta}^{(|\mathcal R|\rho'-1)}$ correspond to the separating hyperplanes for $\mathcal{L}$ identified in Line \ref{algo:typical-oracle} during iteration $\rho'$; if Line \ref{algo:typical-oracle} is not executed for some $r \in \mathcal R$, the respective $\boldsymbol{\beta}$ is simply a zero vector. Let $\texttt{R}^{(\rho)}$ denote $\texttt{R}$ after $\rho$ iterations, defined as \eqref{prob:dcglp} with $\mathcal{L}^\#$ replaced by $(\mathcal{L}^{(\rho)})^\#$, a homogenized version of $\mathcal{L}^{(\rho)}$. Let $(\tau^{(\rho)}, (\boldsymbol{\omega^{r(\rho)}})_{r \in \mathcal R}, \boldsymbol{s^{(\rho)}})$ be the optimal solution of $\texttt{R}^{(\rho)}$, and let $\boldsymbol{\alpha}^{(\rho)}$ denote the optimal dual solution of $\texttt{R}^{(\rho)}$ associated with \eqref{eq:dcglp:d}-\eqref{eq:dcglp:f}. 
  
Since $\texttt{R}^{(\rho)}$ is a relaxation of \eqref{prob:dcglp}, its dual serves as a restriction of \eqref{prob:cgp}. Consequently, the intermediate value $\boldsymbol{\alpha}^{(\rho)}$ still provides a valid inequality for \eqref{prob:benders-reformulation}. This is formally stated in the following proposition:

\begin{prop}\label{prop:inexact}
$\texttt{R}^{(\rho)}$ corresponds to a projection problem:
{\color{black}\begin{equation*}
\tau^{(\rho)} = \min \gamma_{\mathcal C}(( \boldsymbol{x},t) - (\boldsymbol{\hat{x}}, \hat{t}) ) : (\boldsymbol{x},t) \in \overline{\operatorname{conv}}(\mathcal{P}^{(\rho)}_{\mathcal R}),
\end{equation*}}
where $\mathcal{P}$ in \eqref{prob:dcglp} is replaced by $\mathcal{P}^{(\rho)} = (\mathcal{X}_{LP} \times \mathbb{R}) \cap \mathcal{L}^{(\rho)}$. Therefore, at each intermediate iteration $\rho \in \mathbb{Z}_{\geq 0}$ with $\tau^{(\rho)} > 0$, the point  
$
(\sum_{r \in \mathcal R}\boldsymbol{\omega_x^{r(\rho)}}, \sum_{r \in \mathcal R}{\omega_t^{r(\rho)}})
$
is the projection of $(\boldsymbol{\hat{x}}, \hat{t})$ onto $\overline{\operatorname{conv}}(\mathcal{P}^{(\rho)}_{\mathcal R})$. The inequality
$
\alpha_0^{(\rho)} - (\boldsymbol{\alpha_x^{(\rho)}})^{\top} \boldsymbol{x} - \alpha_t^{(\rho)} t \leq 0
$
defines a supporting hyperplane to $\overline{\operatorname{conv}}(\mathcal{P}^{(\rho)}_{\mathcal R}) \supseteq \overline{\operatorname{conv}}(\mathcal{P}_{\color{black}\mathcal{R}})$ that separates $(\boldsymbol{\hat{x}}, \hat{t})$. 
\end{prop}

Proposition \ref{prop:inexact} suggests that solving \eqref{prob:dcglp} inexactly—by terminating Routine \ref{oracle} early and adding the cut derived from the intermediate solution—remains valid. The following example illustrates Proposition \ref{prop:inexact} while showing how disjunctive Benders cuts are found via Routine \ref{oracle}:

\begin{exam}
\textnormal{Consider an example illustrated in Figure \ref{fig:iter0}, where $\mathcal{X} = \mathbb{Z}$ and $\text{dom}f=\mathbb R$. The solid thin black line represents the graph of the value function $f$, and $(\hat{x}, \hat{t})$ is the point to be separated. The set $\mathcal{L}$ corresponds to the epigraph of $f$, and its initial relaxation $\mathcal{L}^{(0)}$ contains a single optimality cut, shown as the thick green line in Figure \ref{fig:iter0}. Suppose we use the disjunction $x \le 0 \lor x \ge 1$. Figure \ref{fig:iter} illustrates the progress of Routine \ref{oracle} when it uses the $\|\cdot\|_2$-norm-based normalization in \eqref{prob:cgp} and employs a \textsc{typicalOracle} that generates the deepest optimality cut based on the $\|\cdot\|_2$-norm in Line \ref{algo:typical-oracle} (see \cite{hosseini2024deepest} for details on the \textsc{typicalOracle}). Let $\rho$ denote the iteration counter of Routine \ref{oracle}, and define the following sets: $\mathcal{P}^{(\rho)} = (\mathcal{X}_{LP} \times \mathbb{R}) \cap \mathcal{L}^{(\rho)}$, $\mathcal{P}^{(\rho)}_1 = \mathcal{P}^{(\rho)} \cap \{(x,t): x \le 0\}$, $\mathcal{P}^{(\rho)}_2 = \mathcal{P}^{(\rho)} \cap \{(x,t): x \ge 1\}$, $\mathcal{P}^{(\rho)}_{\mathcal R} = \mathcal{P}^{(\rho)}_1\cup \mathcal{P}^{(\rho)}_2$ and vectors $(x^{r(\rho)},t^{r(\rho)}) := (\omega_x^{r(\rho)}/ \omega_0^{r(\rho)}, \omega_t^{r(\rho)}/ \omega_0^{r(\rho)}) \in \mathcal{P}^{(\rho)}_r,$ $r=1,2.$
{\color{black}In the figure, the shaded region represents the closed convex hull of $\mathcal{P}^{(\rho)}_{\mathcal R}$.} }

\textnormal{Initially, $(\hat{x}, \hat{t}) \in \mathcal{L}^{(0)}$ and is a convex combination of $(x^{r(0)},t^{r(0)})$, $r=1,2$, as depicted in Figure \ref{fig:iter0}.  
Line \ref{algo:typical-oracle} then generates the two red lines illustrated in Figure \ref{fig:iter0} to separate $(x^{r(0)},t^{r(0)})$ from $\mathcal{L}$ for $r=1,2$. The two red {\color{black}dash dotted} lines refine $\mathcal{L}^{(0)}$ into $\mathcal{L}^{(1)}$, as shown in Figure \ref{fig:iter1}. At this point, $(\hat{x}, \hat{t}) \notin \overline{\operatorname{conv}}(\mathcal{P}^{(1)}_{\mathcal R})$, and the projection of $(\hat{x}, \hat{t})$ onto $\overline{\operatorname{conv}}(\mathcal{P}^{(1)}_{\mathcal R})$ is $(\omega_x^{1(1)} + \omega_x^{2(1)}, \omega_t^{1(1)} + \omega_t^{2(1)})$, a convex combination of $(x^{r(1)},t^{r(1)})$ $r=1,2$. The blue {\color{black}thick dotted} line in Figure \ref{fig:iter1} represents the disjunctive cut generated by the intermediate dual solution $\boldsymbol{\alpha}^{(1)}$. \emph{This intermediate cut still eliminates a considerable portion of $\mathcal{L}$, something that typical oracles cannot achieve.} After one additional iteration, the deepest disjunctive Benders cut for $\mathcal{P}_{\mathcal R}$ is identified as shown in Figure \ref{fig:iter2}.}

\begin{figure}[t!]
    \centering
    \begin{subfigure}[b]{0.32\textwidth}  
        \includegraphics[clip,width=\textwidth]{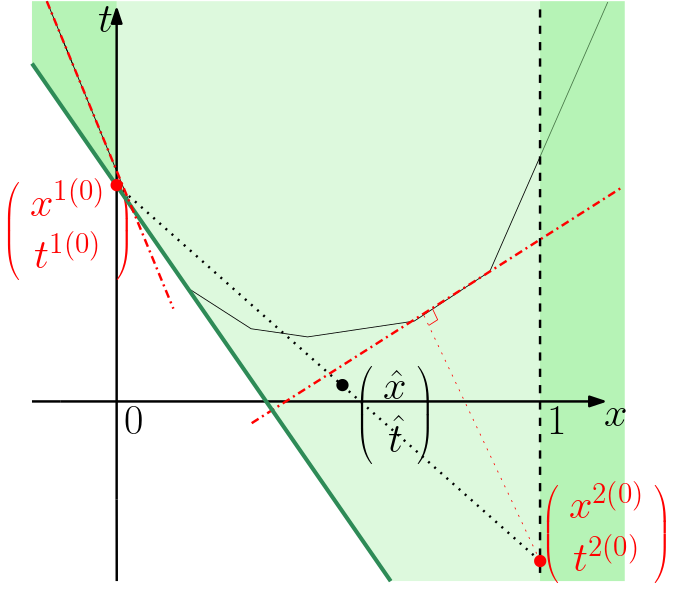}
        \caption{$\rho=0$}
        \label{fig:iter0}
    \end{subfigure}
    \hfill  
    \begin{subfigure}[b]{0.32\textwidth}
        \includegraphics[clip,width=\textwidth]{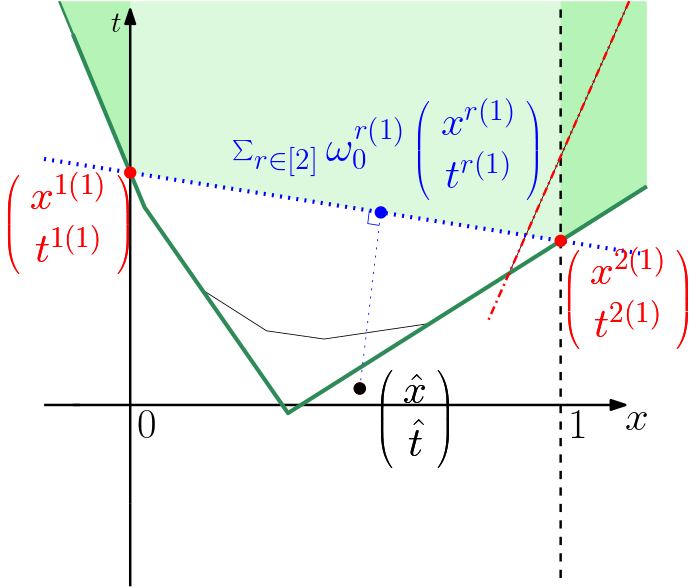}
        \caption{$\rho=1$}
        \label{fig:iter1}
    \end{subfigure}
    \hfill
    \begin{subfigure}[b]{0.32\textwidth}
        \includegraphics[clip,width=\textwidth]{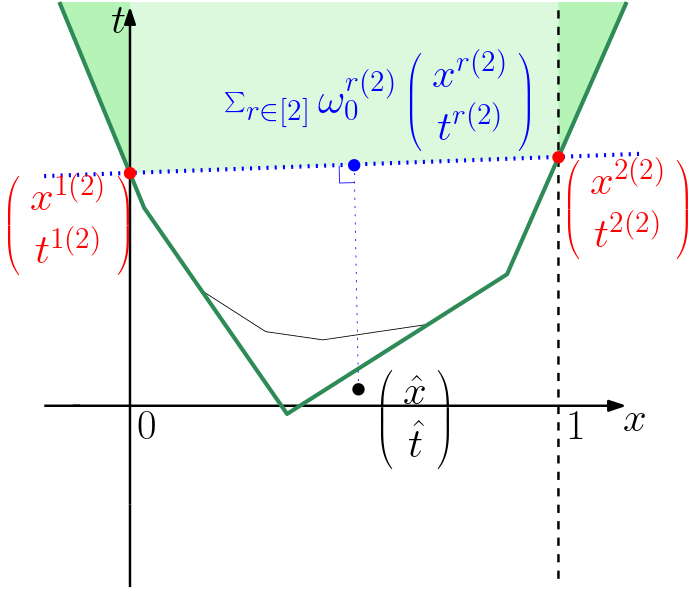}
        \caption{$\rho=2$}
        \label{fig:iter2}
    \end{subfigure}

    \caption{Visualization of Routine \ref{oracle} over iterations}
    \label{fig:iter}
\end{figure}

\end{exam}
\begin{rema}[Benders cuts obtained via automatic perturbation]  \label{rema:byproduct-benders}
  While identifying a supporting hyperplane to $\overline{\operatorname{conv}}(\mathcal{P}_{\mathcal{R}})$ that separates $(\boldsymbol{\hat{x}}, \hat{t})$, Routine \ref{oracle} also produces Benders cuts for separating $(\boldsymbol{x'},t')$ from $\mathcal{L}$ (see Line \ref{algo:typical-oracle}). These byproduct cuts can be integrated into the master problem alongside the disjunctive Benders cut. They can be interpreted as Benders cuts derived by perturbing the separation point $(\boldsymbol{\hat{x}}, \hat{t})$, similar to the in-out method \cite{ben2007acceleration,fischetti2017redesigning}.
   However, unlike the in-out method, which requires a carefully chosen core point, Routine \ref{oracle} performs this perturbation automatically.  
\end{rema}
  
\begin{rema}[Continual refinement of $\mathcal{P}$ in \eqref{prob:dcglp}]\label{rema:dcglp-augmentation}
  Let $\mathcal{P}'$ denote the set $\mathcal{P}$ (the continuous relaxation of the feasible region of \eqref{prob:benders-reformulation}) augmented with a set of $k$ disjunctive Benders cuts $\alpha_0^{(k')} - (\boldsymbol{\alpha_x^{(k')}})^{\top} \boldsymbol{x} - \alpha_t^{(k')} t \geq 0, k'=1,\ldots,k$. Then, Routine \ref{oracle} can generate a {\color{black}higher rank} disjunctive cut for $\overline{\operatorname{conv}}(\mathcal{P}'_{\mathcal R}) \subseteq \overline{\operatorname{conv}}(\mathcal{P}_{\mathcal{R}})$ by incorporating the homogenized version of the disjunctive cuts, i.e., $\alpha_0^{(k')} \omega_0^r - (\boldsymbol{\alpha_x^{(k')}})^{\top} \boldsymbol{\omega_x^r} - \alpha_t^{(k')} \omega_t^r \geq 0$ for $r \in \mathcal R$, into the relaxation \texttt{R}, for $k' \in [k]$. 
\end{rema}
{\color{black}\begin{rema}[Relation to a prior work] 
Routine \ref{oracle} closely resembles the method proposed in \cite{kilincc2017lift}, which focuses on generating lift-and-project cuts, i.e., cuts based on two-term disjunctions, for convex mixed-integer nonlinear programs (MINLPs) using a norm-based normalization constraint. Routine \ref{oracle} can be viewed as a generalization of that approach, as it accommodates nondifferentiable value functions, allows support-function-based normalization (which subsumes norm-based normalization as a special case), and handles multi-term disjunctions. 
Through this connection, one can also view our approach as applicable to the generation of disjunctive cuts for general convex MINLPs, of which the Benders reformulation considered here is a special case.
\end{rema}}

\subsubsection{Decomposable \ref{prob:sub}}
\makeatletter
\renewcommand{\ALG@name}{Routine}
\makeatother
\begin{algorithm}[t!]
  \caption{\textsc{{\color{black}disjunctiveOracle}} for separable subproblems}\label{oracle:separable}
  \begin{algorithmic}[1]
    \Require tolerance $\epsilon > 0$; {\color{black}disjunction $\bigvee_{r \in \mathcal R} U^r \boldsymbol{x} \ge \boldsymbol{u^r}$}; point $(\boldsymbol{\hat{x}}, \boldsymbol{\hat{t}})$ to separate
    \State $\texttt{R}$ $\gets$ a relaxation of \eqref{prob:dcglp:separable} with $\mathcal{L}_j$ replaced by a relaxation of it, including at least one optimality cut or a trivial lower bound on $t_j$ for all $j \in [N]$; 
    $\mathbb{I}_{\mathcal{L}_j^\#}^r \gets \texttt{false}, \ \forall j \in [N], r \in \mathcal R$; 
    \While{$\texttt{!}\land_{j \in [N]}(\land_{r \in \mathcal R}\mathbb{I}_{\mathcal{L}_j^\#}^r)$}
      \State Solve $\texttt{R}$;  $(\hat\tau, (\boldsymbol{\hat\omega^r})_{r \in \mathcal R}, \boldsymbol{\hat s}) \gets$ its optimal solution
      \For{$r \in \mathcal R$} in parallel
      \If{${\hat\omega}_0^{r} > \epsilon$}
      \For{$j \in [N]$} in parallel \label{line:oracle:separable:paralle:start}
      \State $\boldsymbol{x'} \gets \boldsymbol{\hat\omega_x^r} / \hat\omega_0^r;$ $t'_j \gets (\boldsymbol{\hat\omega_{t}^r})_j / {\hat\omega}_0^{r}$
      \State ($\mathbb{I}_{\mathcal{L}_j^\#}^r$, ($\boldsymbol{\beta_x}, \boldsymbol{\beta_{t}}, \beta_0$)) $\gets$ \textsc{typicalOracle}$_{\mathcal{L}_j}$($\boldsymbol{x'}, t'_j$) \label{line:oracle:separable:typical-oracle}
      \If{\texttt{!}$\mathbb{I}_{\mathcal{L}_j^\#}^r$}
      \State Add $\beta_0\omega_0^r - \boldsymbol{\beta_x}^{\top} \boldsymbol{\omega_x^r} - \boldsymbol{\beta_{t}}^{\top} \boldsymbol{\omega_t^r} \leq 0$ for all $r \in \mathcal R$ to $\texttt{R}$ \label{line:oracle:separable:paralle:end}
      \EndIf 
      \EndFor
      \Else
      \State $\mathbb{I}_{\mathcal{L}_j^\#}^r$ $\gets \texttt{true}$ for all $j \in [N]$
      \EndIf
      \EndFor
      \EndWhile
      \State $\hat\alpha_{0}, \boldsymbol{\hat\alpha_{x}}, \boldsymbol{\hat\alpha_t} \gets$ dual values of \eqref{eq:dcglp-gamma:separable:gamma:0}, \eqref{eq:dcglp-gamma:separable:gamma:x} and \eqref{eq:dcglp-gamma:separable:gamma:t}\\
      \Return $(\hat\tau, (\boldsymbol{\hat\alpha_{x}},\boldsymbol{\hat\alpha_{t}},\hat\alpha_{0}))$
  \end{algorithmic}
\end{algorithm}

\textsc{{\color{black}disjunctiveOracle}} can leverage distributed computing resources more effectively when \ref{prob:sub} decomposes into $N$ smaller problems, i.e., $f(\boldsymbol{\hat{x}}) = \sum_{j \in [N]}f_j (\boldsymbol{\hat{x}})$, where $f_j(\boldsymbol{\hat{x}}):= \min_{\boldsymbol{y_j}} \{ \boldsymbol{d_j}^\top\boldsymbol{y_j}: B_j \boldsymbol{y_j} \geq \boldsymbol{b_j} - A_j \boldsymbol{\hat{x}}\}$. 
Let $\mathcal{L}_j = \{(\boldsymbol{x},t_j) \in \mathbb{R}^{n_x+1}: t_j \ge f_j(\boldsymbol{x}), \boldsymbol{x} \in \text{dom}f_j\}$. It is easy to see that $(\boldsymbol{x},t) \in \mathcal{L}$ if and only if there exist $t_1,\cdots, t_N$ such that $t = \sum_{j \in [N]} t_j, (\boldsymbol{x},t_j) \in \mathcal{L}_j, \forall j \in [N]$. Thus, the Benders reformulation \eqref{prob:benders-reformulation} can be equivalently rewritten by replacing $t$ with $\sum_{j \in [N]} t_j$:  
\begin{equation}
\begin{aligned}
\min_{\boldsymbol{x} \in \mathcal{X},\boldsymbol{t} \in \mathbb{R}^N} \{ \boldsymbol{c}^\top \boldsymbol{x} + \sum_{j \in [N]}t_j : (\boldsymbol{x}, t_j) \in \mathcal{L}_j, \forall j \in [N]\}
\end{aligned}
\label{prob:benders-reformulation:separable}
\end{equation}
Similarly, for a candidate solution $(\boldsymbol{\hat{x}}, \boldsymbol{\hat{t}}) \in \mathbb{R}^{n_x} \times \mathbb{R}^N$ obtained from a relaxation of \eqref{prob:benders-reformulation:separable}, the corresponding \eqref{prob:dcglp} can be equivalently expressed as:
\begin{subequations}
\begin{align}
\inf \ & \tau  \\
\mbox{s.t.} \ & (\boldsymbol{\omega_x^r}, (\boldsymbol{\omega_{t}^r})_j, \omega_0^r) \in \mathcal{L}_j^\#, \ \forall j \in [N],r \in \mathcal{R},\\
            & (\boldsymbol{\omega_x^r}, \omega_0^r) \in \mathcal{X}_{LP}^\#, \ \forall r \in \mathcal R,\\
            & U^r \boldsymbol{\omega_x^r} \geq \boldsymbol{u^r}\omega^r_0, \ \forall r \in \mathcal R,\\
            & \omega^r_0 \geq 0,  \ \forall r \in \mathcal R,\\
            [\alpha_0] \quad   & \sum_{r \in \mathcal R}\omega_0^r = 1, \label{eq:dcglp-gamma:separable:gamma:0}\\
            [\boldsymbol{\alpha_x}] \quad   & \boldsymbol{s_x} = \sum_{r \in \mathcal R}\boldsymbol{\omega_x^r} - \boldsymbol{\hat{x}},\label{eq:dcglp-gamma:separable:gamma:x}\\
            [\boldsymbol{\alpha_{t}}] \quad   & \boldsymbol{s_t} = \sum_{r \in \mathcal R}\boldsymbol{\omega_{t}^r} - \boldsymbol{\hat{t}}, \label{eq:dcglp-gamma:separable:gamma:t}\\
            & \boldsymbol{s} \in \tau \mathcal C,\\
            & \tau \ge 0.
\end{align}
\label{prob:dcglp:separable}
\end{subequations}
This naturally extends Routine \ref{oracle} to Routine \ref{oracle:separable}, which enables additional parallel computation in Lines \ref{line:oracle:separable:paralle:start}-\ref{line:oracle:separable:paralle:end}. In Line \ref{line:oracle:separable:typical-oracle}, \textsc{typicalOracle}$_{\mathcal{L}_j}$ is any typical Benders oracle that generates a separating hyperplane for $\mathcal{L}_j$. All properties of Routine \ref{oracle} remain valid for this oracle, including optimality-based termination, provided that \textsc{typicalOracle}$_{\mathcal{L}_j}$ also returns $f_{j}(\boldsymbol{x'})$ for all $j \in [N]$ (Remark \ref{rema:opt-gap-termination}). Additionally, intermediate solutions still yield valid disjunctive cuts for \eqref{prob:benders-reformulation:separable} (Proposition \ref{prop:inexact}). Benders cuts obtained as a byproduct of Routine \ref{oracle:separable} for each $j \in [N]$ can be used to enrich the master problem (Remark \ref{rema:byproduct-benders}). Additionally, \eqref{prob:dcglp:separable} can be continually refined by incorporating raised versions of the disjunctive cuts found so far (Remark \ref{rema:dcglp-augmentation}).

\subsection{Disjunctive Benders decomposition}\label{sec:dbd:dbd}
For the remainder of this paper, we say that a point $(\boldsymbol{\hat{x}}, \hat{t})$ is integral if $\hat{x}_j \in \mathbb{Z}$ for all $j \in \mathcal{I}$; otherwise, we call it fractional. \textsc{{\color{black}disjunctiveOracle}} can be used to separate both integral and fractional candidate solutions, provided that the disjunction is chosen such that $(\boldsymbol{\hat{x}}, \hat{t}) \not \in \overline{\operatorname{conv}}(\mathcal{P}_{\color{black}\mathcal{R}})$. However, selecting such a disjunction for an integral $(\boldsymbol{\hat{x}}, \hat{t})$ is not always straightforward. As a result, we primarily focus on adding disjunctive Benders cuts to separate fractional $(\boldsymbol{\hat{x}}, \hat{t})$.   

We refer to BD equipped with \textsc{{\color{black}disjunctiveOracle}}, i.e., Routine \ref{oracle} or \ref{oracle:separable}, as \emph{disjunctive Benders decomposition}. There are multiple ways to integrate disjunctive Benders cuts into \textsc{BendersSeq} and \textsc{BendersBnB}. For instance, in \textsc{BendersSeq}, one can solve the continuous relaxation of the master problem, and if the resulting candidate solution $(\boldsymbol{\hat{x}}, \hat{t})$ is fractional, use \textsc{{\color{black}disjunctiveOracle}} to generate a disjunctive Benders cut for some disjunction that excludes $\boldsymbol{\hat{x}}$. If the separation is unsuccessful 
or if the candidate solution is integral, a \textsc{typicalOracle} can be called for proper termination.
 Similarly, in \textsc{BendersBnB}, disjunctive Benders cuts can be added at fractional nodes, i.e., nodes where the associated LP subproblems yield fractional solutions, of the branch-and-bound tree. Since \textsc{{\color{black}disjunctiveOracle}} incorporates integrality information that was previously kept within the master problem, the generated cuts can be significantly stronger (see Figure \ref{fig:iter}). 

More interestingly, for certain special cases of \eqref{prob:milp}, such as mixed-binary linear programs, a specialized strategy for constructing splits and leveraging previously found disjunctive cuts (see Algorithm \ref{algo:special}) can entirely eliminate the need to solve the master problem as a MILP in \textsc{BendersSeq}. In other words, in \textsc{BendersBnB}, it can prevent the branch-and-bound tree from expanding beyond the root node, enabling termination at the root.
\makeatletter
\renewcommand{\ALG@name}{Algorithm}
\makeatother
\begin{algorithm}[t!]
  \caption{A specialized disjunctive \textsc{BendersSeq} for mixed-binary linear programs}\label{algo:special}
  \begin{algorithmic}[1]
        \State $k \gets 0$; \texttt{M} $\gets$ a \emph{{continuous}} relaxation of \eqref{prob:benders-reformulation} with a subset of \eqref{eq:benders-reformulation:dom} and \eqref{eq:benders-reformulation:epi}; $\mathbb{I}_{\mathcal{L}} \gets \texttt{false}$ 
        \While{\texttt{true}}
        \State Find an optimal basic solution $(\boldsymbol{\hat x}, \hat t)$ of \texttt{M}
        \State $(\mathbb{I}_{\mathcal{L}}, (\boldsymbol{\beta_x}, \beta_t, \beta_0))\gets$\textsc{typicalOracle}$(\boldsymbol{\hat x}, \hat t)$
        \If{\texttt{!}$\mathbb{I}_{\mathcal{L}}$} Add $\boldsymbol{\beta_x}^\top \boldsymbol{x} + \beta_t t \gets \beta_0$ to \texttt{M}.
        \Else 
        \State $(\boldsymbol{x^{(k)}}, t^{(k)}) \gets (\boldsymbol{\hat x}, \hat t)$ 
        \If{$x_j^{(k)} \in \{0,1\}, \forall j \in \mathcal{I}$} \texttt{break}\label{line:special:disjunctive:begin}
        \EndIf
        \State (Disjunction selection) $x_i \le 0 \lor x_i \ge 1$, where $i$ is the largest index in $\mathcal{I} = \{1,\cdots,\ell\}$ such that $0 < {x}_{i}^{(k)} < 1$, i.e., ${x}_j^{(k)} \in \{0,1\}$ for all $j=i+1,\cdots,\ell$. \label{line:special:disjunctive:split}
        \State (Initialization of \eqref{prob:dcglp}) Construct \eqref{prob:dcglp} for the disjunction with an initial relaxation of $\mathcal{L}$; include the homogenized version of all previously found disjunctive cuts for disjunctions $x_j \le 0 \lor x_j \ge 1$ where $j = 1, \dots, {i-1}$ (see Remark \ref{rema:dcglp-augmentation}).    
        \State $(\hat\tau, (\boldsymbol{\hat\alpha_x},\hat\alpha_t, \hat\alpha_0)) \gets \textsc{{\color{black}disjunctiveOracle}}(\boldsymbol{x^{(k)}}, t^{(k)})$ \label{line:special:disjunctive:end}
        \State Add $\boldsymbol{\hat\alpha_x}^\top \boldsymbol{x} + \hat\alpha_t t \geq \hat\alpha_0$ to \texttt{M} \label{line:special:disjunctive:add}
        \State $k \gets k+1$
        \EndIf
        \EndWhile
        \Return $x^{(k)}$
    \end{algorithmic}
\end{algorithm}
\begin{theo}\label{theo:finite-convergence}
  Suppose \eqref{prob:milp} is a mixed-binary linear program, i.e., $[\underline{x}_j, \overline{x}_j] = [0,1]$ for all $j \in \mathcal{I}$, and {\color{black}the normalization constraint \eqref{eq:cgp:c} is expressed as linear constraints.} Then, Algorithm~\ref{algo:special} terminates in a finite number of iterations with an optimal solution to \eqref{prob:benders-reformulation}.
  \end{theo}
\begin{rema}
  {\color{black}Although primarily theoretical and not expected to have direct practical impact through Algorithm~\ref{algo:special}, Theorem~\ref{theo:finite-convergence} has important implications.} It shows that disjunctive Benders decomposition can eliminate the need to solve the master problem as a MILP, thereby removing the primary computational bottleneck in \textsc{BendersSeq}. It suggests that \textsc{BendersBnB} can terminate at the root node without further branching, highlighting its potential to substantially reduce the size of the branch-and-bound tree. {\color{black}A related but more limited result is established in \cite{rahmaniani2020benders}, where the elimination of the MILP master is shown only for nonseparable subproblems.}
  \end{rema}
\begin{exam}
\begin{figure}[t!]
\centering
\includegraphics[width=0.5\textwidth]{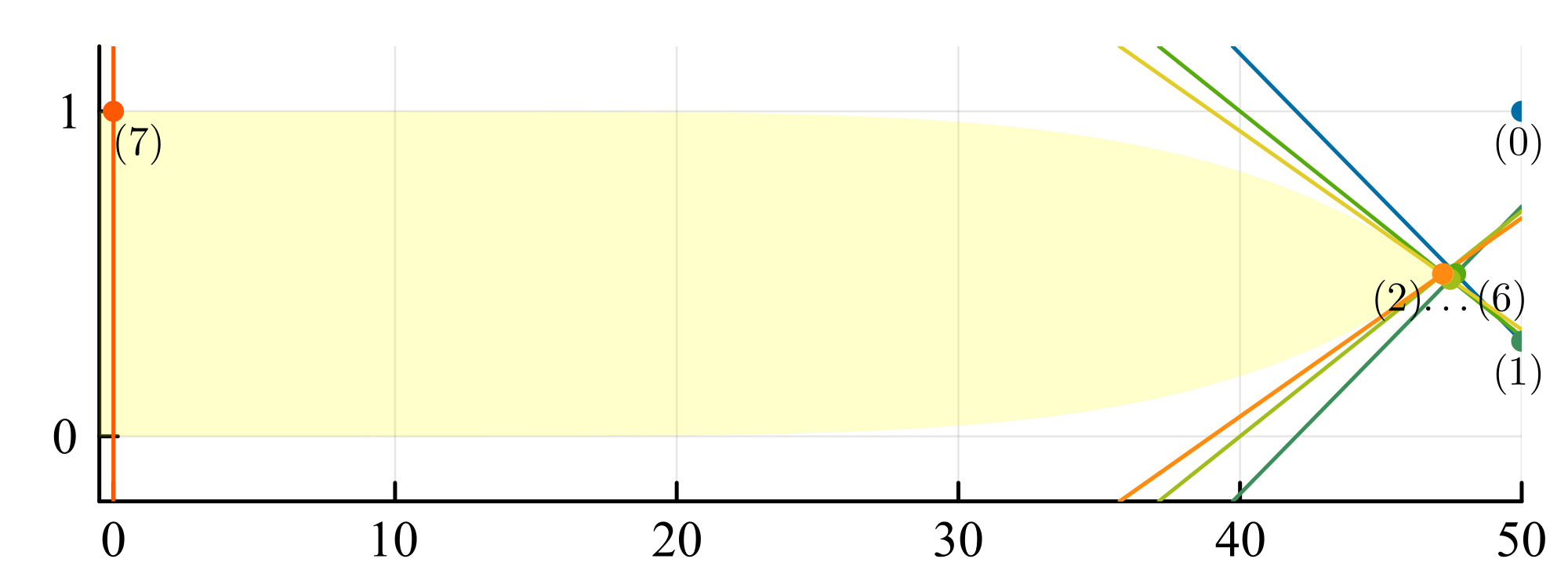}
\caption{Illustration of Algorithm \ref{algo:special} on the motivating example; colors associate the iterates $\{\boldsymbol{x^{(k)}}\}$ and cuts. }
\label{fig:algo:lp}
\end{figure}
\textnormal{Consider applying Algorithm \ref{algo:special} to the motivating example in Section \ref{sec:exam}, despite the lack of a finite convergence guarantee due to $x_2$ being a general integer variable. As shown in Figure \ref{fig:algo:lp}, the algorithm first adds seven Benders cuts until it identifies the optimal solution $\boldsymbol{x^{(6)}}=(47.2,0.5)$ to the LP relaxation of \eqref{prob:benders-reformulation}. At this point, it executes Line \ref{line:special:disjunctive:begin}, constructs the disjunction $x_2 \le 0 \lor x_2 \ge 1$ in Line \ref{line:special:disjunctive:split}, and initiate the approximation of $\mathcal{L}$ with the seven Benders cuts in \eqref{prob:dcglp}. Routine \ref{oracle} then generates the disjunctive Benders cut $x_1 \geq 0$, immediately leading to the identification of the optimal solution after incorporating just one disjunctive cut.}
\end{exam}

\section{Further improvements for mixed-binary linear programs}\label{sec:strenthening-lifting}
In a mixed-binary linear program, where $[\underline{x}_j,\overline{x}_j]=[0,1]$ for $j \in \mathcal{I}$, a disjunctive cut derived from a simple split disjunction $x_i \le 0 \lor x_i \ge 1$ for some $i \in \mathcal{I}$ is known as a lift-and-project cut \cite{balas2002lift}. This type of cut can be strengthened \emph{a posteriori} by exploiting the integrality of other binary variables. Additionally, a lift-and-project cut derived from a subproblem associated with a specific branch-and-bound node can be lifted to ensure its validity across the entire branch-and-bound tree \cite{balas1993lift}. {\color{black}We show how these ideas can be implemented within disjunctive Benders decomposition through a simple post-processing step. In particular, strengthened and lifted lift-and-project cuts can be derived directly from the dual solution of the final relaxation $\texttt{R}$ of \eqref{prob:dcglp} returned by Routine \eqref{oracle}, followed only by a simple elementary operation.
}

Throughout this section, we consider \eqref{prob:dcglp} constructed based on the simple split disjunction where
\[
U^1 \boldsymbol{x} \ge \boldsymbol{u}^1 \;\Leftrightarrow\; -x_i \ge 0,
\qquad
U^2 \boldsymbol{x} \ge \boldsymbol{u}^2 \;\Leftrightarrow\; x_i \ge 1.
\]
and assume that the generated disjunctive Benders cut $\hat\alpha_0 - \boldsymbol{\hat\alpha_x}^\top \boldsymbol{x} - \hat\alpha_t t \leq 0$ cuts off the current iterate $(\boldsymbol{\hat{x}}, \hat{t})$.

\subsection{Cut strengthening}\label{sec:strengthen}
Let $\boldsymbol{\hat\alpha}=(\boldsymbol{\hat\alpha_x}, \hat\alpha_t, \hat\alpha_0)$, $\hat\kappa^{1}$, and $\hat\kappa^{2}$ be the dual solutions associated with \eqref{eq:dcglp:d}-\eqref{eq:dcglp:f}, $-(\boldsymbol{\omega_x^1})_i \ge 0$, and $(\boldsymbol{\omega_x^2})_i \ge \omega_0^2$ in \eqref{prob:dcglp}, resp.
If $\hat\kappa^{1} + \hat\kappa^{2} > 0$, we can strengthen the cut $\hat\alpha_0 - \boldsymbol{\hat\alpha_x}^\top x - \hat\alpha_t t \leq 0$ as 
\begin{equation}
\hat\alpha_0' - (\boldsymbol{\hat\alpha_x'})^\top \boldsymbol{x}- \hat\alpha_t' t \leq 0
\label{eq:strengthened-cut}
\end{equation}
with
\begin{equation*}
\begin{cases}
\hat\alpha_0'  = \hat\alpha_0,\\
(\boldsymbol{\hat\alpha_{x}'})_j  = \min\left\{a^1_j + \hat\kappa^{1} \left\lfloor \frac{a_j^2 - a_j^1}{\hat\kappa^{1} + \hat\kappa^{2}} \right\rfloor, a^2_j -\hat\kappa^{2} \left\lceil  \frac{a_j^2 - a_j^1}{\hat\kappa^{1} + \hat\kappa^{2}} \right\rceil\right\}, \ \forall j \in \mathcal{I},\\
(\boldsymbol{\hat\alpha_{x}'})_j  = (\boldsymbol{\hat\alpha_x})_j, \ \forall j \in [n_x]\setminus \mathcal{I},\\
\hat\alpha_t'  = \hat\alpha_t,
\end{cases}
\end{equation*}
where $a^{r}_j := (\boldsymbol{\hat\alpha_x})_j - \hat\nu^{r}_j$ for $r=1,2, j \in \mathcal{I}$ with $\hat\nu^{r}_j$ for $r=1,2$ being the dual solutions associated with the nonnegativity constraints $(\boldsymbol{\omega_x^r})_j \geq 0$ respectively. 
\begin{prop}\label{prop:strengthen}
If $\hat\kappa^{1} + \hat\kappa^{2} > 0$, \eqref{eq:strengthened-cut} is valid for \eqref{prob:benders-reformulation} and tighter than $\hat\alpha_0 - \boldsymbol{\hat\alpha_x}^\top \boldsymbol{x} - \hat\alpha_t t \leq 0$.
\end{prop}
\begin{rema}
The proof implies that this strengthening procedure remains valid for suboptimal solutions of \eqref{prob:dcglp}, i.e., an intermediate point in Routine \ref{oracle} or \ref{oracle:separable} with ${\tau}^{(\rho)}> 0$. One can simply retrieve the dual solutions for \texttt{R}$^{(\rho)}$ and apply the same strengthening procedure.
\end{rema}

\subsection{Cut lifting} 
\textsc{{\color{black}disjunctiveOracle}}, i.e., Routine \ref{oracle} or \ref{oracle:separable}, can be called at any fractional node of a branch-and-bound tree when solving \eqref{prob:benders-reformulation} in the \textsc{BendersBnB} framework. Each node $\n$ of a branch-and-bound tree is associated with a pair of disjoint index sets $(\mathcal{I}_{\n,0}, \mathcal{I}_{\n,1})$, where $\mathcal{I}_{\n,l}$ denotes the set of indices of binary variables which values are fixed at $l$ at the node. A cut that is generated at $\n$ can be extended to be valid for \eqref{prob:benders-reformulation} by adjusting coefficients for the variables $x_j$ where $j \in \mathcal{I}_\n := \mathcal{I}_{\n,0} \cup \mathcal{I}_{\n,1}$ and the right-hand side. This process is referred to as \emph{lifting} the cut. 

We denote the feasible region of the LP subproblem at a fractional node $\n$ as 
$\mathcal{P}^{\n} := \mathcal{P} \cap \{(\boldsymbol{x},t) \in \mathbb{R}^{n_x+1}: -x_j \geq 0, \forall j \in \mathcal{I}_{\n,0}, x_j \geq 1, \forall j \in \mathcal{I}_{\n,1}\}$. Let $(\boldsymbol{\hat{x}}, \hat{t})$ be the fractional local solution. To simplify notation, let $\mathcal{P}_{i}^{\n}$ (and $\mathcal{P}_i$) denote the intersection of $\mathcal P^\n$ (and $\mathcal P$) with simple split disjunction $x_i \le 0 \lor x_i \ge 1$, respectively. A disjunctive Benders cut for separating $(\boldsymbol{\hat{x}}, \hat{t})$ from $\mathcal{P}_{i}^{\n}$, i.e., the local problem, can be derived by solving \eqref{prob:dcglp} with the addition of the following constraints:
\begin{subequations}
    \begin{align}
[\zeta^r_j] \quad &0 \geq (\boldsymbol{\omega_{x}^r})_{j} , \ \forall j \in \mathcal{I}_{n,0}, r=1,2,\label{eq:CGLP:gamma:lift:1:zeta}\\
 [\xi^r_j] \quad & 0\geq \omega_0^r -(\boldsymbol{\omega_{x}^r})_{j}, \ \forall j \in \mathcal{I}_{\n,1}, r=1,2.\label{eq:CGLP:gamma:lift:1:xi}
\end{align}\label{eq:CGLP:gamma:lift}
\end{subequations}
Suppose $(\boldsymbol{\hat\alpha}, \boldsymbol{\hat\chi}, (\boldsymbol{\hat\zeta^{r}},\boldsymbol{\hat\xi^{r}})_{r =1,2})$ {\color{black}is} a dual solution of the resulting augmented problem. The corresponding cut $\hat\alpha_0 -  \boldsymbol{\hat\alpha_x}^\top \boldsymbol{x} - \hat\alpha_t t \leq 0$ valid for $\mathcal{P}_{i}^\n$ can be lifted to be valid for $\mathcal{P}_{i}$ as follows:
\begin{equation}
\widehat{\alpha_0} - \boldsymbol{\widehat{\alpha_x}}^\top \boldsymbol{x} - \widehat{\alpha_t}t \leq 0,\label{eq:lifted-cut}
\end{equation}
where 
\begin{align*}
 \widehat{\alpha_t} &=  \hat\alpha_t, \\
 \widehat{\alpha_0} &= \hat\alpha_0 - \sum_{j \in \mathcal{I}_{\n,1}}\max\{\hat\xi^1_j, \hat\xi^2_j\}, \\
 (\boldsymbol{\widehat{\alpha_x}})_j &= \begin{cases}
    (\boldsymbol{\hat\alpha_x})_j,  & \forall j \not\in \mathcal{I}_\n,\\
    (\boldsymbol{\hat\alpha_x})_j +\max\{ \hat\zeta^1_j, \hat\zeta^2_j\}, & \forall j \in \mathcal{I}_{\n,0},\\
    (\boldsymbol{\hat\alpha_x})_j -\max\{\hat\xi^1_j, \hat\xi^2_j\}, & \forall j \in \mathcal{I}_{\n,1}.
\end{cases}
\end{align*}

\begin{prop}
The lifted inequality \eqref{eq:lifted-cut} is valid for $\mathcal{P}_i$, i.e., it is valid for any nodes of the branch-and-bound tree. In addition, it still cuts off the local solution $(\boldsymbol{\hat{x}}, \hat{t})$, specifically $\widehat{\alpha_0} - \boldsymbol{\widehat{\alpha_x}}^\top \boldsymbol{\hat{x}} - \widehat{\alpha_t}\hat{t} =  {\hat\alpha_0} - \boldsymbol{\hat\alpha_x}^\top \boldsymbol{\hat{x}} - {\hat\alpha_t}\hat{t}> 0$. 
Specifically, the lifted $\boldsymbol{\widehat{\alpha}}$ satisfies \eqref{eqs:valid-cut} with $\boldsymbol{\widehat{\chi}}$, given as 
\begin{subequations}
\begin{align} 
\forall r \in \{1,2\}, \ &((\widehat{\lambda_{\hat{\pi}}^r})_{\boldsymbol{\hat{\pi}} \in \mathcal{J}}, (\widehat{\mu_{\tilde{\pi}}^r})_{\boldsymbol{\tilde{\pi}} \in \mathcal{R}}, \boldsymbol{\widehat{\theta^r}}, \boldsymbol{\widehat{\kappa^r}}) = ((\hat\lambda_{\hat{\pi}}^r)_{\boldsymbol{\hat{\pi}} \in \mathcal{J}}, (\hat\mu_{\tilde{\pi}}^r)_{\boldsymbol{\tilde{\pi}} \in \mathcal{R}}, \boldsymbol{\hat\theta^r}, \boldsymbol{\hat\kappa^r}), \\ 
&\widehat{\nu^r_j} = \hat\nu^{r}_j - \hat\zeta^r_j + \max\{\hat\zeta^1_j,\hat\zeta^2_j\}, \ \forall j \in \mathcal{I}_{\n,0}, \\ 
&\widehat{\eta^r_j} = \hat\eta^r_j - \hat \xi^r_j + \max\{\hat\xi^1_j,\hat\xi^2_j\}, \ \forall j \in \mathcal{I}_{\n,1}. 
\end{align}
\label{eq:lifting:farkas}
\end{subequations}
\label{prop:lifting}
\end{prop}

\begin{rema}\label{rema:strenthening-lifted-cut} 
A lifted inequality can be further strengthened by applying the strengthening technique presented in Section \ref{sec:strengthen}. 
One can simply apply \eqref{eq:strengthened-cut} to $(\boldsymbol{\widehat{\alpha}}, \boldsymbol{\widehat{\chi}})$, as given in \eqref{eq:lifted-cut} and \eqref{eq:lifting:farkas}, to obtain the strengthened lifted cut. 
\end{rema}

\subsection{\textsc{ApxDisjunctiveOracle}: an approximate version of \textsc{{\color{black}disjunctiveOracle}}}
Proposition \ref{prop:lifting} motivates an approximate version of \textsc{{\color{black}disjunctiveOracle}} for generating lift-and-project cuts, which offers computational efficiency while guaranteeing the generation of a hyperplane that separates $(\boldsymbol{\hat{x}}, \hat{t})$ from $\overline{\operatorname{conv}}(\mathcal{P}_{\color{black}\mathcal{R}})$ if and only if $(\boldsymbol{\hat{x}}, \hat{t})\not\in \overline{\operatorname{conv}}(\mathcal{P}_{\color{black}\mathcal{R}})$.
\begin{algorithm}[t!]
\makeatletter
\renewcommand{\ALG@name}{Routine}
\makeatother
  \caption{\textsc{ApxDisjunctiveOracle} for generating a lift-and-project cut}\label{oracle-approx}
  \begin{algorithmic}[1]
    \Require tolerance $\epsilon > 0$; point $(\boldsymbol{\hat{x}}, \hat{t})$ to separate; split disjunction $x_i \le 0 \lor x_i \ge 1$ for some $i \in \mathcal I$
    \State Perform Line \ref{algo:initialization} in Routine \ref{oracle}
    \State $\mathcal{I}_z \gets \{j \in \mathcal{I}: \hat{x}_j = z\}$ for $z \in \{0,1\}$
    \State Add to \texttt{R} the inequalities \eqref{eq:CGLP:gamma:lift}, replacing $\mathcal{I}_{\n,z}$ with $\mathcal{I}_z$ for $z \in \{0,1\}$\label{oracle:approx:const}
    \State Execute Lines \ref{algo:while:begin}–\ref{algo:while:end} of Routine \ref{oracle}  \label{oracle-approx:while}
    \State $\hat\tau \gets$ optimal objective value of \texttt{R}; $\boldsymbol{\hat\alpha},(\boldsymbol{\hat\xi},\boldsymbol{\hat\zeta}) \gets$ optimal dual solutions associated with \eqref{eq:dcglp:d}-\eqref{eq:dcglp:f} and \eqref{eq:CGLP:gamma:lift} respectively in \texttt{R}  
    \State Compute $\boldsymbol{\widehat{\alpha}}$ using \eqref{eq:lifted-cut} applied to $(\boldsymbol{\hat\alpha},\boldsymbol{\hat\xi},\boldsymbol{\hat\zeta})$  \\
    \Return $(\hat\tau, (\boldsymbol{\widehat{\alpha_x}}, \widehat{\alpha_t}, \widehat{\alpha_0}))$
  \end{algorithmic}
\end{algorithm}
\begin{prop}
Let $(\hat\tau, \boldsymbol{\widehat{\alpha}})$ be the output of Routine \ref{oracle-approx}. 
Consider \eqref{prob:cgp} constructed for the split disjunction $x_i \le 0 \lor x_i \ge 1$ and the point $(\boldsymbol{\hat{x}}, \hat{t})$, and let $v_{\eqref{prob:cgp}}$ be its optimal objective value. Then, $\frac{\hat\tau}{\max\{\sigma_{\mathcal C}(\boldsymbol{\widehat{\alpha_x}}), 1\}} \leq v_{\eqref{prob:cgp}} \leq \hat\tau$, and $\frac{1}{\max\{\sigma_{\mathcal C}(\boldsymbol{\widehat{\alpha_x}}), 1\}}\boldsymbol{\widehat{\alpha}}$ is $\left(1-\frac{1}{\max\{\sigma_{\mathcal C}(\boldsymbol{\widehat{\alpha_x}}), 1\}}\right)$-optimal to \eqref{prob:cgp}.
 \label{prop:approx-oracle}
\end{prop}

\begin{rema}
When a fractional $\boldsymbol{\hat{x}}$ has many integral components, \textsc{ApxDisjunctiveOracle} can be more computationally efficient than \textsc{{\color{black}disjunctiveOracle}}, as it restricts the feasible region substantially. Nevertheless, it guarantees the generation of a separating hyperplane if and only if $(\boldsymbol{\hat{x}}, \hat{t}) \not \in \overline{\operatorname{conv}}(\mathcal{P}_{\color{black}\mathcal{R}})$, ensuring proper termination; note that $\hat\tau$ can be zero if and only if $(\boldsymbol{\hat{x}}, \hat{t}) \in \overline{\operatorname{conv}}(\mathcal{P}_{\color{black}\mathcal{R}})$ due to Proposition \ref{prop:approx-oracle}. Furthermore, the cut generated by \textsc{ApxDisjunctiveOracle} can be further strengthened following the approach in Remark \ref{rema:strenthening-lifted-cut}.
\end{rema}

\section{Discussion}\label{sec:discussion}
\textsc{{\color{black}disjunctiveOracle}} provides a recipe for leveraging existing typical Benders oracles—those that separate a point from $\mathcal{L}$—to discover valid inequalities for $\overline{\operatorname{conv}}((\mathcal{X} \times \mathbb{R}) \cap \mathcal{L})$.
It actively identifies previously unknown constraints that are crucial for finding a supporting hyperplane to a set that lies between $\overline{\operatorname{conv}}((\mathcal{X} \times \mathbb{R}) \cap \mathcal{L})$ and $\mathcal{L}$. 
This approach contrasts with other methods that rely on pre-generated Benders cuts to construct valid inequalities (e.g., mixed-integer rounding enhancement \cite{bodur2017mixed}); while such methods can be handled internally by off-the-shelf solvers to some extent, they can also be combined with the proposed approach for enhanced effectiveness.
In this section, we compare the proposed method with existing oracles that fall into the category shown in Figure \ref{fig:oracle-proposed}, which do not rely solely on pre-generated Benders cuts.
For simplicity, we assume for the rest of this section that $\mathcal{I} = [n_x]$ and $\text{dom}f_j = \mathbb{R}^{n_x}$ for all $j \in [N]$, i.e., the subproblems have complete recourse. The discussion that will follow will remain true in general cases without the assumption.

The project-and-cut approach proposed in \cite{bodur2017strengthened} also aims to generate a disjunctive cut for \eqref{prob:benders-reformulation}. To address the incomplete characterization of $\mathcal{L}$ in the evolving master problem—a relaxation of \eqref{prob:benders-reformulation}—during BD, they propose solving a cut-generating LP constructed from the original formulation \eqref{prob:milp}. However, this approach can be computationally demanding for large-scale problems, as the CGLP typically doubles the size of the original problem or is at least as large as the original formulation. Possibly due to this computational burden, their experiment focuses on generating disjunctive cuts for the current master problem, only using the current relaxation of $\mathcal{L}$. 

In the context of \textsc{{\color{black}disjunctiveOracle}}, this corresponds to the initial iteration of Routine \ref{oracle} or \ref{oracle:separable} applied to a relaxation of \eqref{prob:dcglp}, initialized with the added Benders cuts to the master problem. 
Consequently, the resulting cut is weaker than the one produced by \textsc{{\color{black}disjunctiveOracle}}, which, \emph{beyond the initial iteration, iteratively discovers unknown Benders cuts to generate a stronger disjunctive cut}. They also experiment with solving the CGLP for the original formulation on a per-scenario basis, possibly to mitigate the problem size explosion as more scenarios are incorporated. In contrast, Routine \ref{oracle:separable} can harness parallel computing resources effectively to generate disjunctive cuts with all the scenarios taken into account simultaneously.

The cut-and-project approach in \cite{bodur2017strengthened} and the Lagrangian cut in \cite{rahmaniani2020benders,zou2019stochastic} are inspired by the following Lagrangian relaxation of \eqref{prob:benders-reformulation} typically used in dual decomposition: 
\begin{subequations}
\begin{align}
\min \ & \boldsymbol{c}^{\top} \boldsymbol{x} + \sum_{j \in [N]}t_j\\
\mbox{s.t.} \ & (\boldsymbol{x},t_j) \in \overline{\operatorname{conv}}((\mathcal{X} \times \mathbb{R}) \cap \mathcal{L}_j), \ \forall j \in [N],
\end{align}
\label{prob:dd}    
\end{subequations}
where $(\mathcal{X} \times \mathbb{R}) \cap \mathcal{L}_j = \{(\boldsymbol{x},t_j): t_j \geq f_j(\boldsymbol{x}), \ \boldsymbol{x} \in \mathcal{X}\}$. Both the cut-and-project approach and the Lagrangian cut aim to separate a candidate solution $(\boldsymbol{\hat{x}}, \hat{t}_j)$ of \eqref{prob:dd} from $\overline{\operatorname{conv}}((\mathcal{X} \times \mathbb{R}) \cap \mathcal{L}_j))$ whenever $(\boldsymbol{\hat{x}}, \hat{t}_j) \not\in \overline{\operatorname{conv}}((\mathcal{X} \times \mathbb{R}) \cap \mathcal{L}_j))$ for each $j \in [N]$. 
They generate cuts based on two different descriptions of $\overline{\operatorname{conv}}((\mathcal{X} \times \mathbb{R}) \cap \mathcal{L}_j))$.
Specifically, the cut-and-project approach uses the following description:
\begin{align}
    \overline{\operatorname{conv}}((\mathcal{X} \times \mathbb{R}) \cap \mathcal{L}_j))  &= \overline{\operatorname{conv}}(\{(\boldsymbol{x},t_j): t_j \geq f_j(\boldsymbol{x}), \ \boldsymbol{x} \in \mathcal{X}\})\nonumber\\
    & = \overline{\operatorname{conv}}(\{(\boldsymbol{x},t_j): \exists \boldsymbol{y_j} \mbox{ such that } (\boldsymbol{x},t_j,\boldsymbol{y_j}) \in \mathcal{Q}_j\}),\nonumber\\
    & = \overline{\operatorname{conv}}(\text{Proj}_{(\boldsymbol{x},t_j)}\mathcal{Q}_j)\nonumber\\
    & = \text{Proj}_{(\boldsymbol{x},t_j)}(\overline{\operatorname{conv}}(\mathcal{Q}_j)),\label{eq:cut-and-project}
\end{align}
where $\mathcal{Q}_j := \{(\boldsymbol{x},t_j,\boldsymbol{y_j}): t_j \geq \boldsymbol{d_j}^{\top} \boldsymbol{y_j}, A_j \boldsymbol{x} + B_j \boldsymbol{y_j} \geq \boldsymbol{b_j}, \boldsymbol{x} \in \mathcal{X}\}$, and the last equality is from \cite[Exercise 3.34]{conforti2014integer}. On the other hand, the Lagrangian cut leverages the following alternative description of $\overline{\operatorname{conv}}((\mathcal{X} \times \mathbb{R}) \cap \mathcal{L}_j)$: 
\begin{equation*}
\begin{aligned}
\overline{\operatorname{conv}}((\mathcal{X} \times \mathbb{R}) \cap \mathcal{L}_j) & = \overline{\operatorname{conv}}(\{(x,t_j): t_j \geq f_j(x), \ x \in \mathcal{X}\}) = \overline{\operatorname{conv}}(\text{epi}f'_j) = \text{epi}({f_j'}^{**}),
\end{aligned}
\end{equation*}
where $f_j'$ is a function that equals $f_j(x)$ for $x \in \mathcal{X}$ and $\infty$ otherwise; ${f_j'}^{**}$ is the biconjugate of $f'_j$, whose epigraph is the closure of the convex hull of the epigraph of $f'_j$ \cite[Proposition 1.6.1]{bertsekas2009convex}. 
Since we have
\begin{equation*}
\begin{aligned}
{f_j'}^{**}(\boldsymbol{x}) & = \max_{\boldsymbol{\pi_j}}\boldsymbol{\pi_j}^{\top} \boldsymbol{x} - {f_j'}^{*}(\boldsymbol{\pi_j}) = \max_{\boldsymbol{\pi_j}}\boldsymbol{\pi_j}^{\top} \boldsymbol{x} + \min_{\boldsymbol{z_j}}f'_j(\boldsymbol{z_j}) - \boldsymbol{\pi_j}^{\top} \boldsymbol{z_j} = \max_{\boldsymbol{\pi_j}}\boldsymbol{\pi_j}^{\top} \boldsymbol{x} + \min_{\boldsymbol{z_j} \in \mathcal{X}}f_j(\boldsymbol{z_j}) - \boldsymbol{\pi_j}^{\top} \boldsymbol{z_j} \\
& = \max_{\boldsymbol{\pi_j}}\boldsymbol{\pi_j}^{\top} \boldsymbol{x} + \min_{\boldsymbol{z_j} \in \mathcal{X},\boldsymbol{y_j}}\boldsymbol{d_j}^{\top} \boldsymbol{y_j} - \boldsymbol{\pi_j}^{\top} \boldsymbol{z_j} : A_j \boldsymbol{z_j} + B_j \boldsymbol{y_j} \geq \boldsymbol{b_j},
\end{aligned}
\end{equation*}
it follows that
\begin{equation}
\overline{\operatorname{conv}}((\mathcal{X} \times \mathbb{R}) \cap \mathcal{L}_j) = \{(\boldsymbol{x},t_j):t_j \geq \max_{\boldsymbol{\pi_j}}\boldsymbol{\pi_j}^{\top} \boldsymbol{x} + \min_{\boldsymbol{z_j} \in \mathcal{X},\boldsymbol{y_j}}\{\boldsymbol{d_j}^{\top} \boldsymbol{y_j} - \boldsymbol{\pi_j}^{\top} \boldsymbol{z_j} : A_j \boldsymbol{z_j} + B_j \boldsymbol{y_j} \geq \boldsymbol{b_j}\}\}.
\label{eq:lagrangian-cut}
\end{equation}

Given a point $(\boldsymbol{\hat{x}}, (\hat{t}_j)_{j \in [N]})$ to separate, the cut-and-project approach, based on \eqref{eq:cut-and-project}, repeatedly (i) solves \eqref{prob:sub} to obtain the optimal $\boldsymbol{\hat{y}_j}$ with the optimal objective value $f_j(\boldsymbol{\hat{x}})$, (ii) generates an inequality that separates the point $(\boldsymbol{\hat{x}}, \hat{f}_j(\boldsymbol{\hat{x}}), \boldsymbol{\hat{y}_j})$ from $\overline{\operatorname{conv}}\mathcal{Q}_j$, if such an inequality exists, and (iii) incorporates the inequality into \eqref{prob:sub}. Upon termination, the Benders cut is derived from the final modified subproblem. The Lagrangian cut, based on \eqref{eq:lagrangian-cut}, is constructed by first determining a (sub)optimal $\boldsymbol{\hat{\pi}_j}$ for the maximization problem in \eqref{eq:lagrangian-cut} given $\boldsymbol{\hat{x}}$. Then, the inner MILP subproblem is solved to obtain an optimal solution $(\boldsymbol{\hat{y}_j}, \boldsymbol{\hat{z}_j})$ for the given $\boldsymbol{\hat{\pi}_j}$ and $\boldsymbol{\hat{x}}$. Finally, the cut $t_j \geq \boldsymbol{\hat{\pi}_j}^{\top} (\boldsymbol{x} - \boldsymbol{\hat{z}_j}) + \boldsymbol{d_j}^{\top} \boldsymbol{\hat{y}_j}$ is generated.

It is well known that \eqref{prob:dd} is equivalent to \eqref{prob:benders-reformulation} only when $N = 1$; for $N>1$ it yields a relaxation. {\color{black}This is illustrated in the following example:
\begin{exam}[Nonzero duality gap when $N>1$]
Consider an instance of \eqref{prob:benders-reformulation:separable}, where $N=2, \mathcal X = \{0,1\}^2$, $\boldsymbol c = \boldsymbol 0$, $f_1(\boldsymbol x) = \max\{x_1 + x_2 -1, -x_1-x_2+1\}$, and $f_2(\boldsymbol x) = \max\{x_1 - x_2, x_2-x_1\}$. It is easy to verify that the optimal objective value of \eqref{prob:benders-reformulation:separable} is 1 via enumeration. Its Lagrangian relaxation reduces to $\min (f_1^{'**}+f_2^{'**})(\boldsymbol{x})= \min (f_1^{'}+f_2^{'})(\boldsymbol{x}) = \min (f_1+f_2)(\boldsymbol{x}): \boldsymbol{x}\in [0,1]^2$, which attains 0 at $(0.5, 0.5)$. Hence, a nonzero duality gap arises.
\label{exam:lag}
\end{exam}}
Due to this, by design, oracles generating the Lagrangian cuts and the cut-and-project cuts do not eliminate the need to solve the master problem as an MILP unless $N = 1$. 
In contrast, the proposed disjunctive Benders decomposition can eliminate the need to solve the master problem as an MILP for any $N \in \mathbb{Z}_{>0}$ for mixed-binary linear programs (see, e.g., Theorem \ref{theo:finite-convergence}).  
Additionally, \textsc{{\color{black}disjunctiveOracle}} identifies effective Benders cuts for each scenario as a byproduct of finding the deepest disjunctive cut for \eqref{prob:benders-reformulation}. These cuts can be further incorporated into the master problem.

{\color{black}This observation does not contradict Lemma~2 in \cite{bodur2017strengthened}, which shows that generating rank-1 split cuts in an extended space (i.e., in \eqref{prob:milp}) and projecting them can yield stronger relaxations than generating such cuts directly in the projected space (i.e., in \eqref{prob:benders-reformulation}). This difference arises because Algorithm~\ref{algo:special} attains finite convergence by generating higher-rank cuts in a principled way.}

Furthermore, even in the nonseparable case, that is $N=1$, \textsc{{\color{black}disjunctiveOracle}} can be more effective at managing problem scale. Generating Lagrangian cuts for $\overline{\operatorname{conv}}((\mathcal{X} \times \mathbb{R}) \cap \mathcal{L}_j)$ requires solving MILP subproblems, which may not scale well for large problems. Similarly, in the cut-and-project scheme, when $N=1$, the process reduces to generating valid cuts for the original formulation \eqref{prob:milp}, which can also face scalability limitations. In contrast, \textsc{{\color{black}disjunctiveOracle}} requires only iterations of typical oracles, while solving a continuous problem \texttt{R}, whose size can be effectively controlled through proper constraint management between iterations. 

{\color{black}Lastly, an additional key advantage of \textsc{disjunctiveOracle} is that it preserves the subproblem's dual feasible region, i.e., \textsc{typicalOracle} is reused as-is, unlike in the aforementioned approaches. This property is crucial in BD, as it enables warm-starting and the exploitation of special subproblem structure through an unchanged dual feasible region—features that are often key to the success of BD. For example, in the uncapacitated facility location problem, BD can significantly outperform general-purpose solvers \cite{fischetti2017redesigning}, since the subproblem decomposes into a collection of knapsack problems that admit efficient closed-form solutions. This advantage is lost if the subproblem's feasible region is modified by introducing additional binary variables or adding inequalities.}

To summarize, for problems of manageable size, where solving the MILP master/sub problem or generating valid inequalities for the original formulation on a per-scenario basis remains computationally efficient, the aforementioned oracles can be more effective than disjunctive Benders decomposition. However, \emph{the key advantage of disjunctive Benders decomposition lies in its scalability for very large problems, achieved by strategically leveraging simple, typical oracles.}

\section{Numerical experiments} \label{sec:experiment}
As discussed in Section \ref{sec:discussion}, the key advantage of disjunctive BD lies in its ability to tackle large problems by strategically leveraging simple, typical oracles. 
To demonstrate this, we evaluate the proposed method on large-scale instances of the uncapacitated facility location problem (\textbf{UFLP}) 
and the stochastic network interdiction problem (\textbf{SNIP})—each presenting distinct structural features, which we describe in detail later. 
Formal problem formulations are provided in the electronic companion.
We assess computational performance through a comparative analysis with conventional BD, augmented with structure-exploiting \textsc{typicalOracle}s, as well as off-the-shelf solvers.

\subsection{Datasets}
\paragraph{UFLP} We utilize the \texttt{KG} instances from \cite{korkel1989exact} and \cite{ghosh2003neighborhood}. The benchmark comprises 90 problem instances, each with an equal number of customers and potential facilities from $\{250, 500, 750\}$ under varying cost configurations. 

\paragraph{SNIP}
We use the \texttt{snipno3} and \texttt{snipno4} test sets, each consisting of five base instances, originally introduced by \cite{pan2008minimizing} and subsequently used in \cite{bodur2017strengthened, hosseini2024deepest}. These instances are obtained from the online companion of \cite{hosseini2024deepest}. Each instance comprises 456 scenarios and 320 binary first-stage decision variables, and all share a common network topology with 783 nodes and 2,586 arcs.
The \texttt{snipno3} and \texttt{snipno4} instances differ in the probability of undetected traversal when a sensor is not installed. In \texttt{snipno3}, the probability of traversing a link undetected without a sensor is 0.1 times that with a sensor, whereas in \texttt{snipno4}, it is 0. For each base instance, the sensor installation budget $b$ varies from 30.0 to 90.0 in increments of 10, resulting in 70 total instances.

\subsection{Baseline Algorithms}
\begin{itemize}[leftmargin=*]
    \item \texttt{EXT}: For all problem instances, we include the extensive formulation \eqref{prob:milp} as a baseline, solved using an off-the-shelf solver. 
    \item \texttt{CBD}: As another baseline, we implement conventional \textsc{BendersBnB} for each problem, where a \textsc{typicalOracle} is invoked via a lazy constraint callback whenever an incumbent solution is identified. For \textbf{UFLP}, we use the solver-free \textsc{typicalOracle} proposed in \cite{fischetti2017redesigning}, which exploits the separable structure of the subproblem: once the facility opening decision \( x \) is fixed as $\hat x$, \ref{prob:sub} decomposes into independent continuous knapsack problems—one per customer—with closed-form solutions. We adopt the fat version of the method, where individual Benders cuts are generated for each violated customer constraint, which outperforms the slim (aggregated cut) version. 
    For \textbf{SNIP}, we use the conventional \textsc{typicalOracle}, in which the subproblem \eqref{prob:sub} is decomposed by scenario and solved via an off-the-shelf solver.
    \item \texttt{DBD:} For each problem, we implement disjunctive \textsc{BendersBnB}, incorporating the same lazy constraint callback used in \texttt{CBD}, along with an additional user callback.
     This user callback is equipped with \textsc{{\color{black}disjunctiveOracle}} (Routine~\ref{oracle}) or its approximate version (Routine~\ref{oracle-approx}), both of which internally leverage the same \textsc{typicalOracle} used in the lazy callback. The user callback invokes \textsc{{\color{black}disjunctiveOracle}} (or its approximation) periodically—every 250 fractional nodes—for both \textbf{UFLP} and \textbf{SNIP}.
\end{itemize}

For both \texttt{CBD} and \texttt{DBD}, one could opt to add Benders cuts as user cuts via \textsc{typicalOracle}s at each fractional node, up to a specified number of visits, as done in \cite{fischetti2017redesigning}. However, in our experiment, we omit such user cuts in both \texttt{CBD} and \texttt{DBD}, as we observed that they significantly increase computation time on certain easy \textbf{UFLP} instances. Furthermore, since the primary objective of this paper is to evaluate the strength of cutting planes, we intentionally refrain from using advanced heuristics or solver strategies beyond those implemented through custom cut callback functions, even if such enhancements were considered in prior studies.

\subsubsection{Implementation details}\label{sec:implementation}
 All implementations were developed in Julia, using JuMP's solver-independent callback framework, and all optimization models were solved using IBM ILOG CPLEX 22.1.1.0 through the JuMP modeling interface. All computational experiments were conducted on a Linux system equipped with two Intel Xeon Gold 6330 CPUs @ 2.00~GHz, each with 56 cores. We executed 8 experiments in parallel, allocating 60~GB of memory and 7 CPU cores to each, with a time limit of 4 hours per experiment.
 The source code is publicly available at \url{https://github.com/asu-opt-lab/BendersX.git}. {\color{black}In addition, for the \textbf{UFLP} instances, we conducted supplementary experiments using Gurobi 13.0.0 for comparison; these results are reported in Appendix~\ref{apx:gurobi}.}

\paragraph{Solver settings}
For \texttt{EXT}, we primarily use the default parameter settings of both CPLEX and Gurobi, with several modifications to control parallelism and improve numerical precision. Specifically, we set the number of threads to $7$, the integer feasibility tolerance to $10^{-9}$, the feasibility tolerance to $10^{-9}$, and the MIP gap tolerance to $10^{-6}$.

For both \texttt{CBD} and \texttt{DBD}, we apply the same parameter settings to the master problem as in \texttt{EXT}, except for the thread parameter when using CPLEX. Since CPLEX restricts callback-based executions to a single thread, the master problems in \texttt{CBD} and \texttt{DBD} are solved in single-threaded mode. In addition, for the \textbf{UFLP} instances, we set the branch direction parameter to prioritize upward branching, following \cite{fischetti2017redesigning}. This setting is applied consistently to \texttt{EXT} and to the master problems of both \texttt{CBD} and \texttt{DBD}.

For solving \texttt{R} (a relaxation of \eqref{prob:dcglp}) and \eqref{prob:sub} within solver-based \textsc{typicalOracle}s, we use the default solver settings except for several adjustments aimed at improving numerical stability and solution accuracy. In particular, we set the number of threads to $7$, the feasibility tolerance to $10^{-9}$, the optimality tolerance to $10^{-9}$, and enable the solver’s numerical emphasis option.


\paragraph{Root node processing} 
For all \textsc{BendersBnB} variants, we begin by solving the LP relaxation of \eqref{prob:benders-reformulation} at the root node using \textsc{BendersSeq} 
equipped with the \textsc{typicalOracle}s described earlier. For \textbf{UFLP} and \textbf{SNIP}, the LP relaxations of all instances are solved to a tolerance of $10^{-9}$ within 10 seconds and 80 seconds, respectively.

\paragraph{\textsc{(approx)DisjunctiveOracle} configurations} For all problems, we use simple split of the form \( x_i \le 0 \lor x_i \ge 1\) for some \( i \in [\ell] \), constructed using the variable with the largest fractional value. The generated lift-and-project cuts are strengthened via \eqref{eq:strengthened-cut} before being added to the master problem, and all previously generated disjunctive cuts are included in \eqref{prob:dcglp} (see Remark \ref{rema:dcglp-augmentation}). We 
terminate \textsc{(approx)DisjunctiveOracle} either when the optimality gap falls below $10^{-3}$, or earlier if the lower bound fails to improve over three consecutive iterations. 
{\color{black}For \textbf{UFLP}, we use the reverse polar normalization with $\boldsymbol y = (\boldsymbol{y_x}, \boldsymbol{y_t}) = (\boldsymbol{0}, \boldsymbol{1})$, and for \textbf{SNIP}, we use the norm-based normalization with $p=\infty$.} For \textbf{UFLP}, the solver-free \textsc{typicalOracle} is computationally cheap but often yields a large number of cuts, which can hinder the performance of the master problem's branch-and-bound process. To prioritize cut quality and incorporate cuts selectively, we add only the top 5\% of the most violated byproduct Benders cuts to the master problem, and we do not reuse the relaxation \texttt{R}—that is, all byproduct Benders cuts generated in earlier disjunctive cut generation are discarded from \texttt{R}. Lifting is disabled, meaning \textsc{{\color{black}disjunctiveOracle}} is used in its non-approximated form. 
For \textbf{SNIP}, the \textsc{typicalOracle} is more computationally intensive as it solves 456 solver-based subproblems. 
To reduce the frequency of solving \ref{prob:sub}, we reuse the relaxation \texttt{R}, add all violated byproduct Benders cuts to the master problem, and employ \textsc{ApxDisjunctiveOracle} as described in Routine \ref{oracle-approx}. 

{\color{black}
\begin{figure}[t!]
    \centering
     \begin{subfigure}[b]{0.45\textwidth}
        \includegraphics[width=\textwidth]{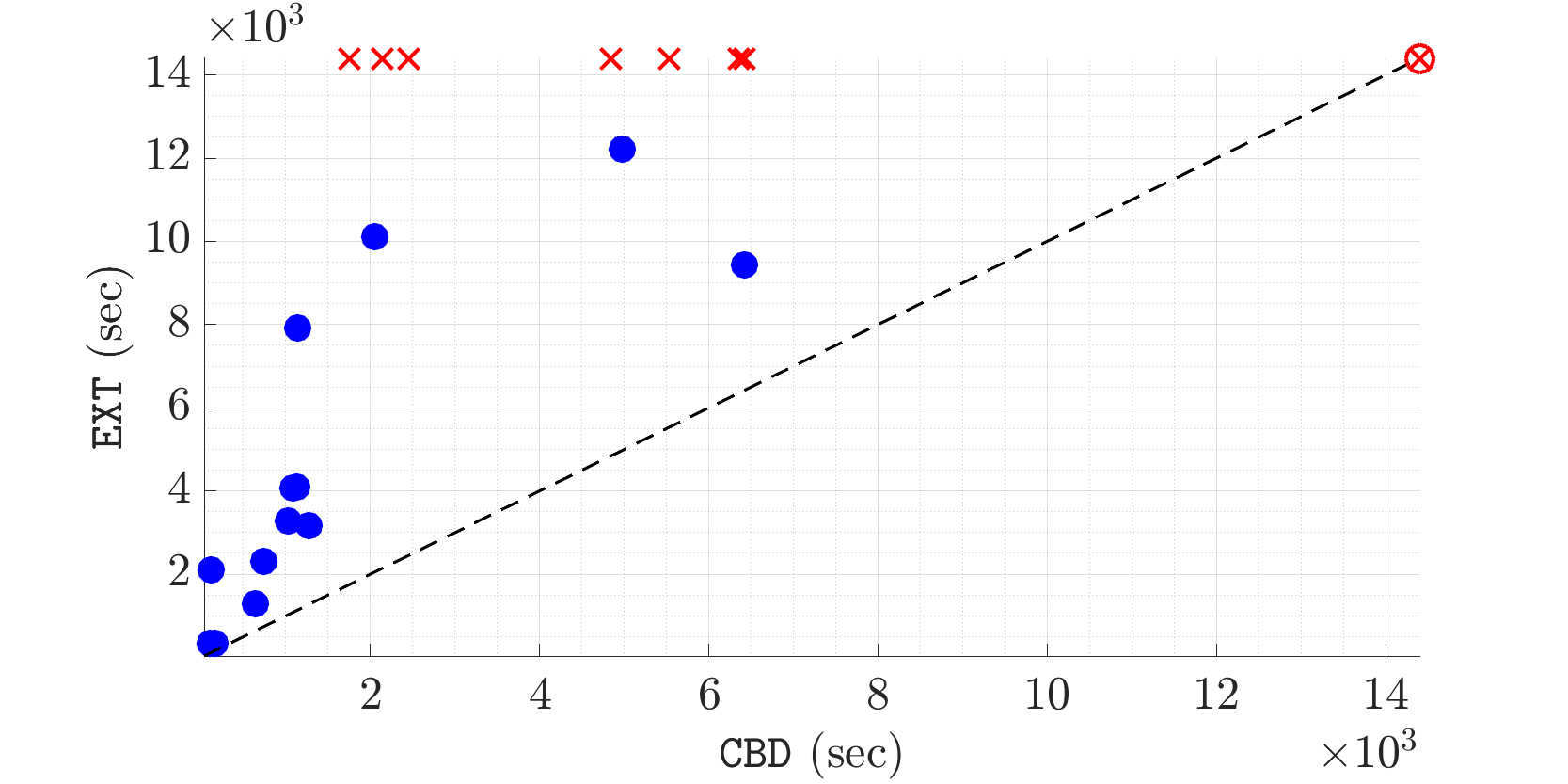}
        \caption{Time, \texttt{CBD} and \texttt{EXT}}
        \label{fig:time-ext-250-cplex}
    \end{subfigure}
    \hspace{4mm}
 \begin{subfigure}[b]{0.45\textwidth}
        \includegraphics[width=\textwidth]{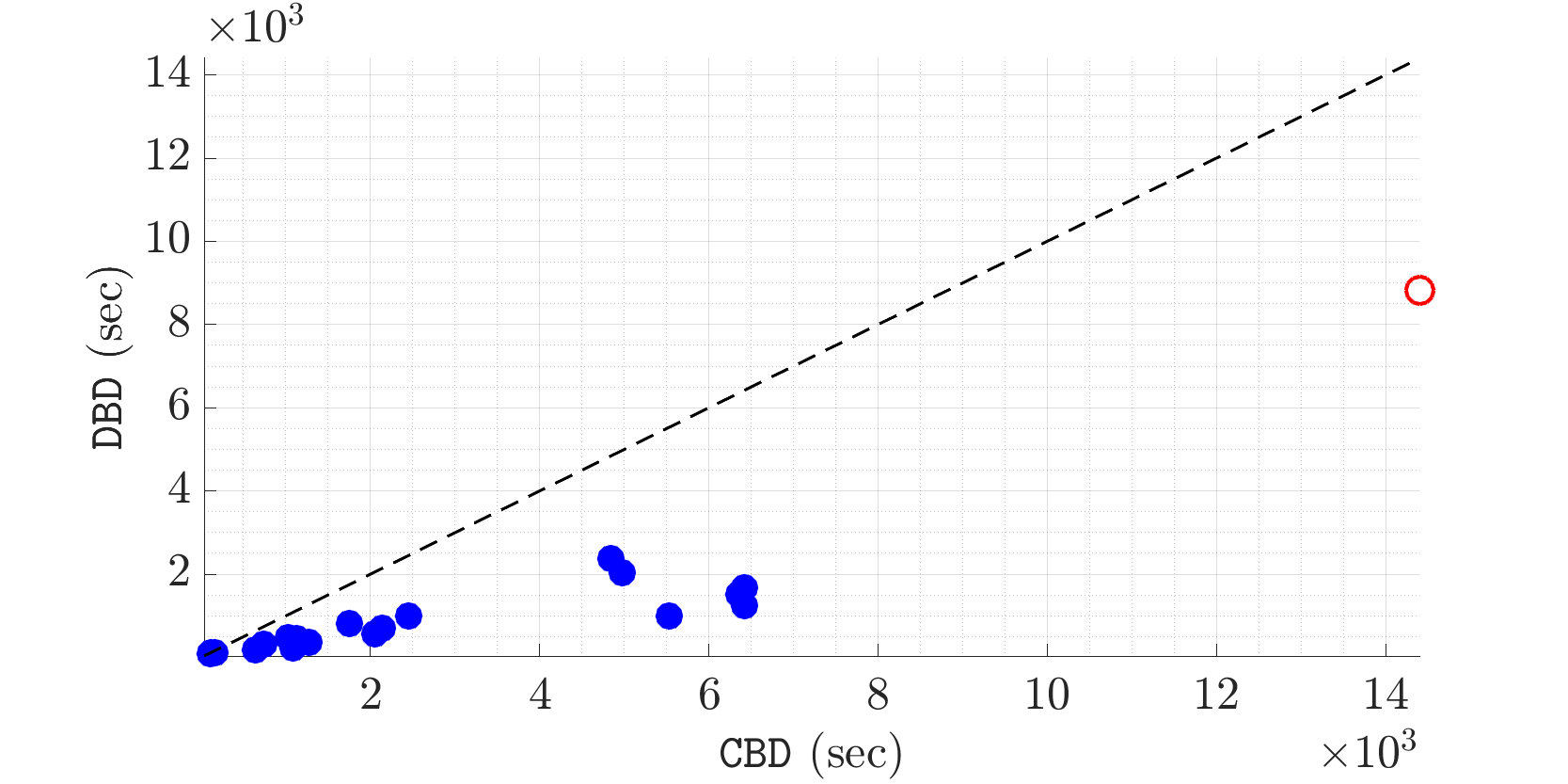}
        \caption{Time, \texttt{CBD} and \texttt{DBD}}
        \label{fig:time-dbd-250-cplex}
    \end{subfigure}
    \begin{subfigure}[b]{0.45\textwidth}  
        \includegraphics[width=\textwidth]{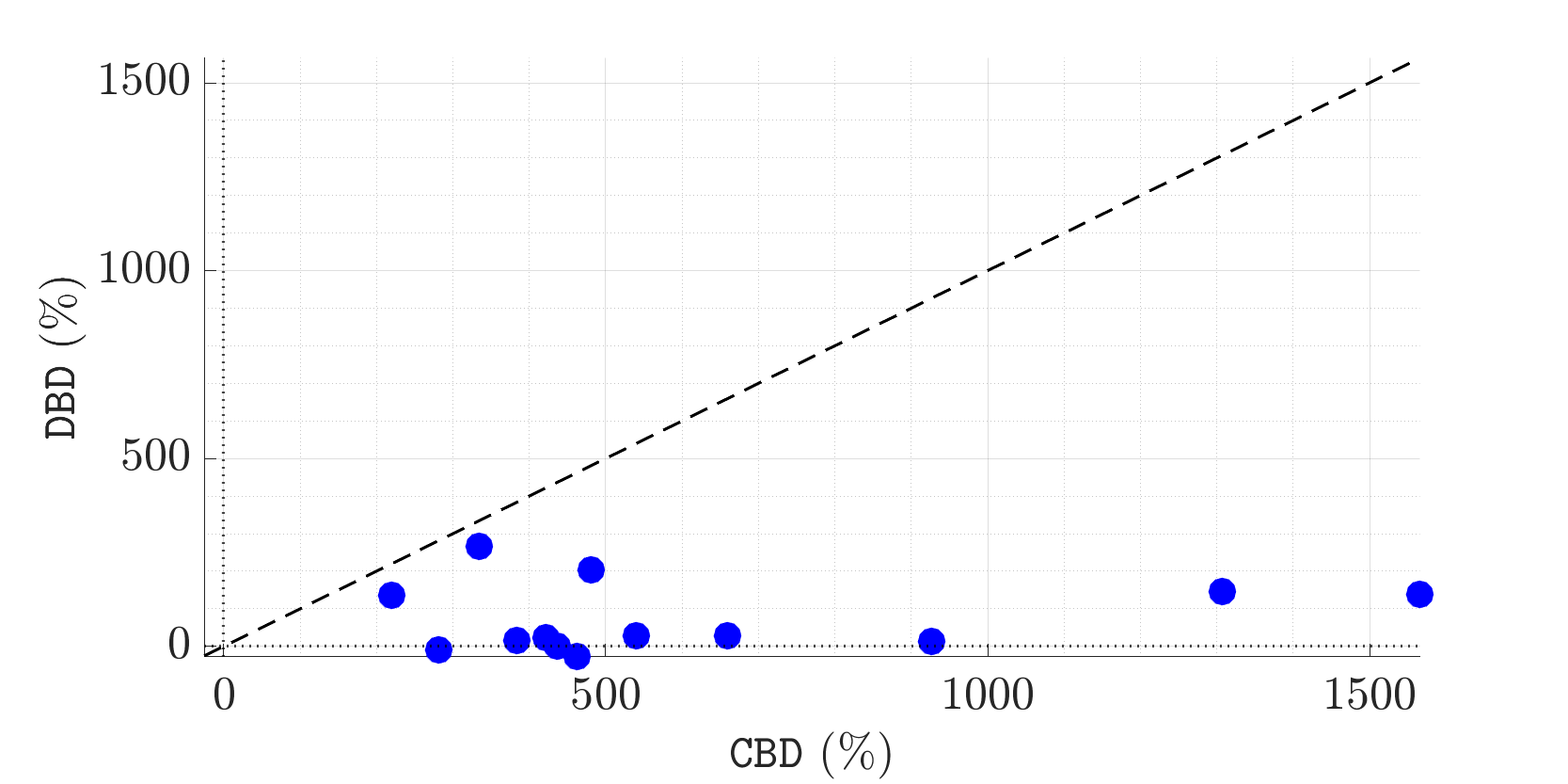}
        \caption{Relative increase in NC with respect to \texttt{EXT} on instances solved to optimality by \texttt{EXT}}
        \label{fig:nc-rel-250-cplex}
    \end{subfigure}
    \hspace{4mm}
    \begin{subfigure}[b]{0.45\textwidth}
        \includegraphics[width=\textwidth]{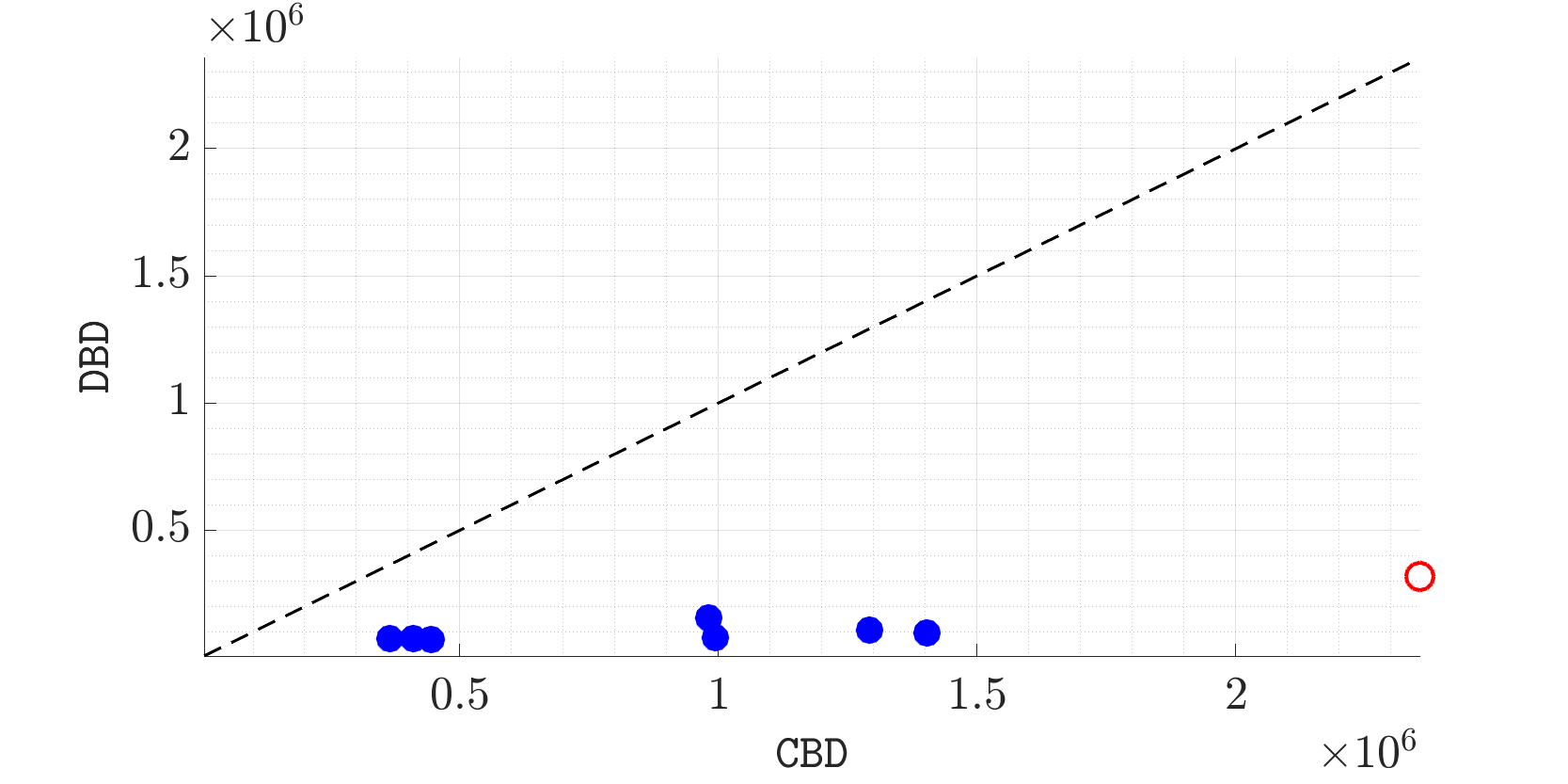}
        \caption{NC on instances for which \texttt{EXT} reached the time limit.}
        \label{fig:nc-abs-250-cplex}
    \end{subfigure}
    \caption{UFLP-250 instances; the unfilled circle and the \texttt{x} marker denote the instances for which \texttt{CBD} and \texttt{EXT}, respectively, reached the 4-hour time limit.}
    \label{fig:250-cplex}
\end{figure}
\subsection{UFLP}\label{sec:experiment:uflp}
Figure~\ref{fig:250-cplex} presents benchmark results for the 21 instances, out of the 30 instances of size 250, excluding 9 easy instances for which both \texttt{CBD} and \texttt{DBD} require less than 100 seconds. The comparison focuses on runtime and node count (NC) of the final branch-and-bound tree. In terms of runtime, Figures~\ref{fig:time-ext-250-cplex} and~\ref{fig:time-dbd-250-cplex} show that \texttt{CBD} outperforms \texttt{EXT}, while \texttt{DBD} yields a further improvement. Notably, \texttt{EXT} fails to solve 8 of the 21 instances within the time limit. Although \texttt{CBD} performs substantially better—solving 20 instances—it still remains slower than \texttt{DBD}. In contrast, \texttt{DBD} solves all 21 instances within the time limit and achieves the lowest runtime.

Figures~\ref{fig:nc-rel-250-cplex} and~\ref{fig:nc-abs-250-cplex} compare the NCs. Figure~\ref{fig:nc-rel-250-cplex} reports the relative increase in NC with respect to \texttt{EXT} for the 13 instances solved optimally by \texttt{EXT}. On average, \texttt{CBD} increases the node count by 617.1\% relative to \texttt{EXT}, whereas \texttt{DBD} increases it by only 74.2\%. \emph{Remarkably, in many instances, \texttt{DBD} achieves branch-and-bound trees comparable in size to those of \texttt{EXT}, despite operating on a decomposed formulation.} This behavior is particularly noteworthy because one would naturally expect \texttt{EXT}, which retains the complete formulation, to produce the smallest search trees. Indeed, \texttt{CBD} often requires an order of magnitude more nodes than \texttt{EXT}, reflecting the information loss induced by decomposition; its runtime advantage stems primarily from the computational savings associated with solving smaller modular problems. In contrast, \texttt{DBD} recovers much of the strength of the extensive formulation while still benefiting from the computational advantages of decomposition. 
Figure~\ref{fig:nc-abs-250-cplex} further highlights the node efficiency of \texttt{DBD} relative to \texttt{CBD}. 
The performance advantage of \texttt{DBD} becomes even more pronounced on the seven instances for which \texttt{CBD} requires more than one hour to reach optimality (see Figure~\ref{fig:time-dbd-250-cplex}). For these instances, \texttt{DBD} reduces the average runtime to 38\% of that of \texttt{CBD}, while the average NC drops to only 11.52\% of that required by \texttt{CBD}.

\begin{figure}[t!]
    \centering
    \begin{subfigure}[b]{0.48\textwidth}  
        \includegraphics[clip,width=\textwidth]{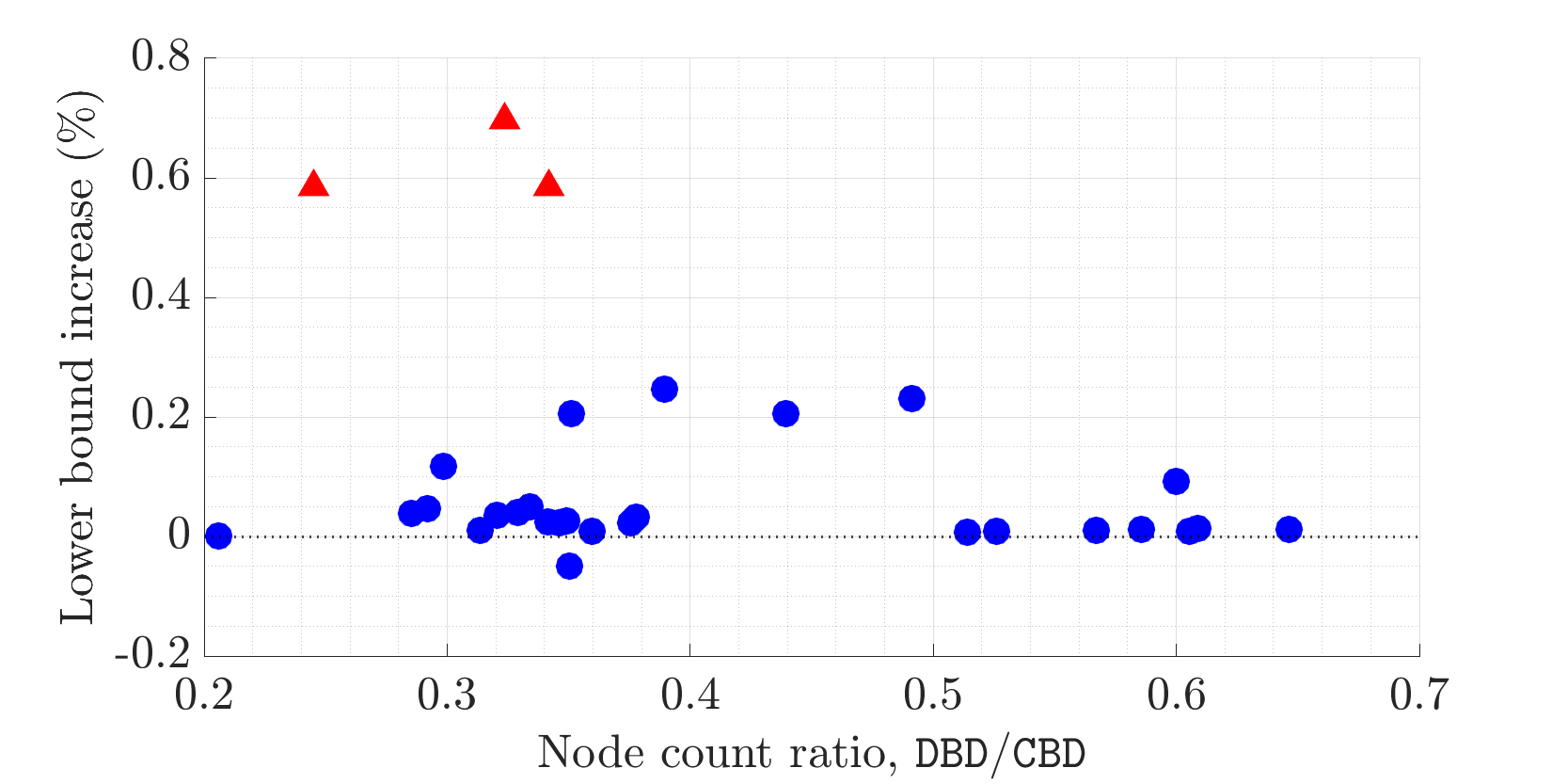}
        \caption{UFLP-500 instances}
        \label{fig:nc-500-cplex}
    \end{subfigure}
    \hspace{4mm}
    \begin{subfigure}[b]{0.48\textwidth}
        \includegraphics[clip,width=\textwidth]{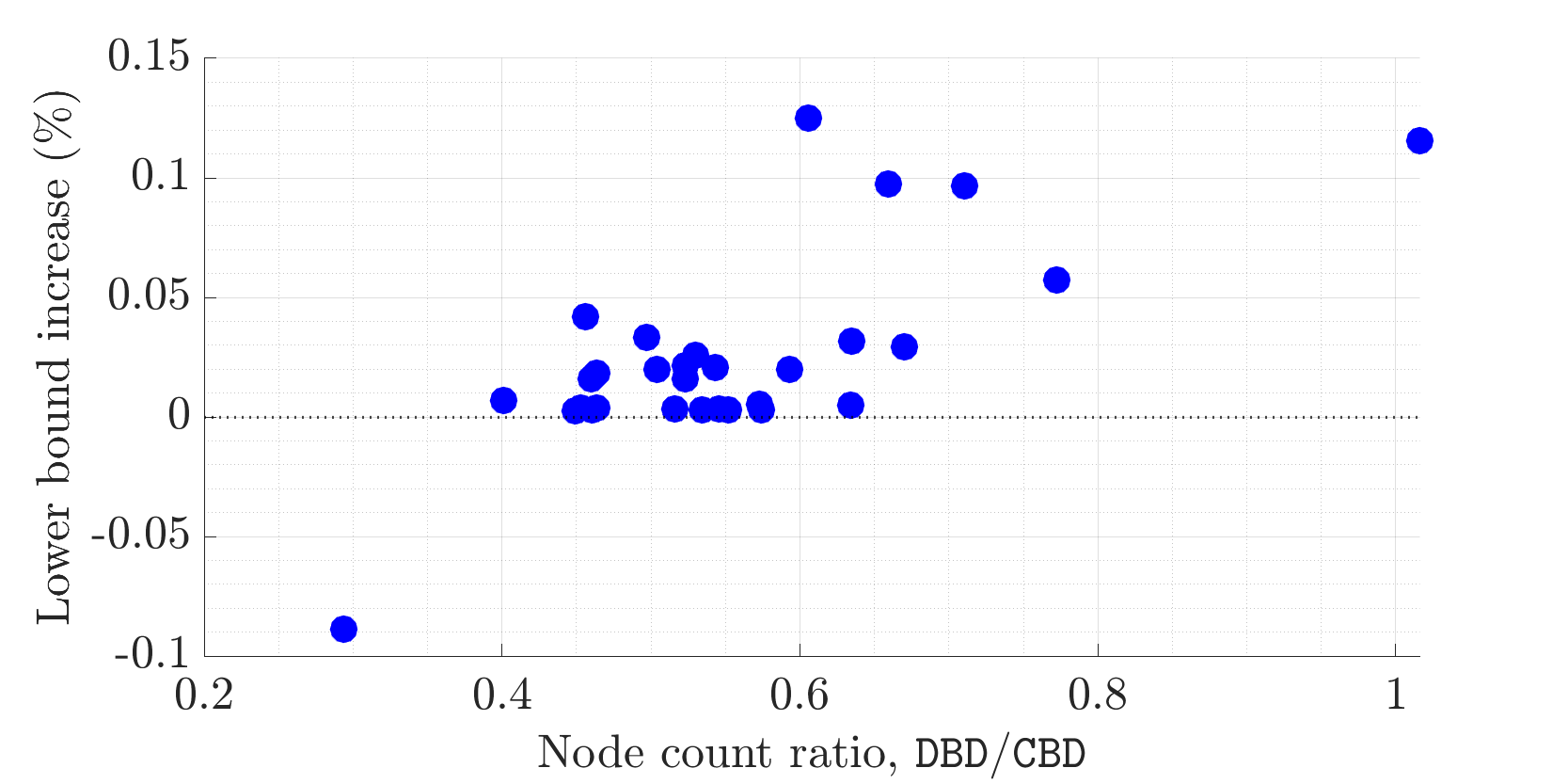}
        \caption{UFLP-750 instances}
        \label{fig:nc-750-cplex}
    \end{subfigure}
    \caption{Relationship between the NC ratio and the relative lower-bound increase. Triangles denote instances for which \texttt{CBD} reaches the time limit, whereas \texttt{DBD} solves to optimality.}
    \label{fig:nc-500-750-cplex}
\end{figure}

For the remaining 60 instances of sizes 500 and 750, \texttt{EXT} reached the 4-hour time limit and performed worst among the three methods. Consequently, we focus on comparing the lower bounds and the number of explored nodes of \texttt{CBD} and \texttt{DBD} in order to evaluate the strength of the generated cuts. We intentionally exclude upper-bound comparisons, as they are more sensitive to heuristic behavior and therefore less informative regarding cut quality. Figure~\ref{fig:nc-500-750-cplex} illustrates the relationship between the NC ratio (\texttt{DBD}/\texttt{CBD}) and the relative lower-bound improvement achieved by \texttt{DBD} over \texttt{CBD}. Across both the UFLP-500 and UFLP-750 instances, most points lie below an NC ratio of 1 and above a lower-bound improvement of zero, indicating that \texttt{DBD} almost always explores fewer nodes than \texttt{CBD} and produces stronger lower bounds within the same computational time budget. Moreover, the observed lower-bound improvements are not negligible: the average final optimality gap of \texttt{CBD} on these instances is only 0.65\%, within which even small improvements are difficult for \texttt{CBD} to achieve despite substantial computational effort. Notably, \texttt{DBD} solves three of these instances to optimality, whereas \texttt{CBD} fails to solve any within the time limit.
}


\subsection{SNIP}
For all \texttt{snipno3} and \texttt{snipno4} instances, \texttt{EXT} fails to reach optimality within the 4-hour time limit. However, both test sets are considered relatively easy under the BD framework, as \texttt{CBD} solves all 70 instances with an average runtime of 194.41 seconds. Only one case—Instance~1 in \texttt{snipno3} with budget~90—takes longer than 15 minutes, terminating after 1,274.76 seconds. Since \texttt{CBD} is already highly efficient on these instances, the additional cost of solving all 456 subproblems multiple times to generate disjunctive cuts is not offset. As a result, \texttt{DBD} incurs greater computational overhead: on average, it requires 2.30 times more runtime than \texttt{CBD} for \texttt{snipno3}, and 1.81 times more for \texttt{snipno4}, despite exploring only about 79\% as many nodes as \texttt{CBD}.

Nonetheless, our results show that \texttt{DBD} can outperform \texttt{CBD} in both runtime and node count particularly in instances where \texttt{CBD} struggles. In the most difficult case—Instance~1 of \texttt{snipno3} with budget~90—\texttt{DBD} solves the problem approximately 100 seconds faster than \texttt{CBD}, despite the overhead from subproblem solves. This advantage stems from the significant reduction in the final branch-and-bound tree: \texttt{DBD} explores only 16\% of the nodes required by \texttt{CBD}. To further investigate, we conducted additional experiments with budget values ranging from 100 to 150. We identified Instance~1 with budget~90, and Instances~0 and 1 with budget~140 in \texttt{snipno3} as the most challenging. For these instances, \texttt{CBD} requires an average of 1,764.14 seconds, while \texttt{DBD} completes them in just 79.2\% of that time. Moreover, the average node count ratio across these instances is only 10.1\%, highlighting the substantial improvement in node efficiency.

To summarize, our experiments on \textbf{SNIP} align with the conclusions of \cite{bodur2017strengthened}: as instances become easier, the cost of repeatedly solving subproblems may outweigh the benefits of stronger cuts. However, the results indicate that \texttt{DBD} remains particularly effective for difficult instances, demonstrating its scalability, and its performance can be further enhanced through parallel computing or scenario sampling strategies.

\section{Conclusions} \label{sec:conclusion}
This paper introduced disjunctive Benders decomposition (DBD), a new enhancement of classic BD, which integrates disjunctive programming theory to generate stronger cuts valid for the convex hull of the Benders reformulation. Unlike conventional BD—which typically produces cuts tight only for the continuous relaxation—our approach explicitly incorporates integrality information into the cut-generation process. Specifically, we developed both exact and approximate variants of a disjunctive cut-generating oracle for the Benders reformulation that can operate on top of any existing typical Benders oracles. The proposed disjunctive oracles construct valid inequalities via iterative separation in an extended space, without introducing additional integer variables or constraints into the subproblem. 
We also extended standard strengthening and lifting techniques—originally developed for disjunctive cuts based on simple binary-variable splits—into the BD framework. These enhancements improve cut quality and reduce the computational burden of the proposed oracles while ensuring their validity throughout the entire branch-and-bound tree.
{\color{black}In addition, we introduced a unified normalization framework for cut-generating programs that encompasses widely used norm-based and reverse polar normalization schemes. The framework recovers their geometric interpretations in a unified manner and provides a tool for designing new normalization schemes.
}
Computational experiments on large-scale uncapacitated facility location (UFLP) and stochastic network interdiction (SNIP) problems demonstrated the practical advantages of DBD. In the UFLP setting, DBD consistently outperforms both the extensive formulation and conventional BD, solving more instances to optimality and reducing node counts by up to a factor of ten. 
Notably, DBD retains much of the strength of original formulations while exploiting the computational benefits of decomposition, especially on challenging UFLP instances. For SNIP, while subproblem complexity impacts overall runtime, DBD achieves significantly better node efficiency on harder instances, highlighting its promise for solving large-scale mixed-integer programming problems.
In summary, DBD offers a scalable and general-purpose framework for high-performance mixed-integer optimization. Future work includes accelerating oracle execution through bundle-type methods, leveraging parallelism more effectively, and developing adaptive strategies for disjunction selection to further enhance performance.

%

\bibliographystyle{plain}
\bibliography{reference}

\begin{appendix}

\section{Omitted proofs}
{\color{black}
\subsection{Proof of Proposition \ref{prop:reverse}}
\begin{align*}
\overline{\operatorname{conv}}(\mathcal S) & = \bigcap_{(\boldsymbol{a_x},a_0) \in \mathcal S^*}\{\boldsymbol{x}:\alpha_0 - \boldsymbol{\alpha_x}^\top x \le 0\}\\
& = \{\boldsymbol{x}: (-\boldsymbol{x}, 1)^\top (\boldsymbol{\alpha_x}, \alpha_0) \le 0, \ \forall (\boldsymbol{a_x},a_0) \in \mathcal S^*\}\\
& = \{\boldsymbol{x}: (-\boldsymbol{x}, 1) \in (\mathcal S^*)^\circ\},
\end{align*}
where the first equality follows from the separating hyperplane theorem \cite[Corollary 11.5.1]{rockafellar1997convex}.\qed

\subsection{Proof of Theorem \ref{theo:dcgp}}\label{proof:theo:dcgp}

\noindent\textbf{Proof.}
Let $\boldsymbol{p}:= (-\boldsymbol{\hat x}, 1)$ and $\boldsymbol{\alpha} = (\boldsymbol{\alpha_x}, \alpha_0)$. Since \eqref{prob:cgp} is a convex optimization problem, $\boldsymbol{0} \in \mathcal S^*$, and $\sigma_{\mathcal C}(\boldsymbol{0})=0<1$, Slater's condition holds. Therefore, it has strong duality. 

Using $\mathcal C' := \mathcal C \times \{0\}$, we have $\sigma_{\mathcal C'}(\boldsymbol{\alpha}) = \sigma_{\mathcal C}(\boldsymbol{\alpha}_x)$. The dual can thus be written as
\begin{align}
\eqref{prob:cgp}
&= \inf_{\tau \ge 0} \sup_{\boldsymbol{\alpha} \in \mathcal{S}^*}
\boldsymbol{p}^\top \boldsymbol{\alpha} + \tau - \tau \sigma_{\mathcal C'}(\boldsymbol{\alpha}) \nonumber\\
&= \inf_{\tau \ge 0} \tau + \sup_{\boldsymbol{\alpha} \in \mathcal{S}^*}
\boldsymbol{p}^\top \boldsymbol{\alpha} - \tau \sigma_{\mathcal C'}(\boldsymbol{\alpha}) \nonumber\\
&= \inf_{\tau \ge 0,\boldsymbol{s}} \tau + (\tau \sigma_{\mathcal C'})^*(\boldsymbol s)
: \boldsymbol{p}-\boldsymbol{s} \in (\mathcal{S}^*)^\circ,
\label{eq:cgp-dual-1}
\end{align}
where the last equality follows from the conjugate identity
\(
(f_1+f_2)^* = f_1^*\square f_2^*
\)
with $f_1=\tau \sigma_{\mathcal C'}$ and $f_2=\delta_{\mathcal S^*}$ \cite[Theorem 16.4]{rockafellar1997convex}. This identity applies without closure because $\mathcal C'$ is compact, hence $\sigma_{\mathcal C'}$ is real-valued. Since $\mathcal S^*$ is a convex cone,
\(
(\delta_{\mathcal S^*})^*=\delta_{(\mathcal S^*)^\circ}.
\)

Moreover, for $\tau>0$,
\(
(\tau \sigma_{\mathcal C'})^*(\boldsymbol s)=\tau \sigma_{\mathcal C'}^*\!\left(\frac{\boldsymbol s}{\tau}\right),
\)
while for $\tau=0$,
\(
(\tau \sigma_{\mathcal C'})^*=\delta_{\{\boldsymbol 0\}}.
\)
Therefore,
\begin{equation}
\eqref{eq:cgp-dual-1}
=
\inf_{\tau \ge 0,\boldsymbol{s}}
((1+\sigma_{\mathcal C'}^*)\tau)(\boldsymbol s)
:
\boldsymbol{p}-\boldsymbol{s} \in (\mathcal S^*)^\circ,
\label{eq:cgp-dual:2}
\end{equation}
where
\[
((1+\sigma_{\mathcal C'}^*)\tau)(\boldsymbol{s}) :=
\begin{cases}
\tau (1+\sigma_{\mathcal C'}^*)\!\left(\frac{\boldsymbol{s}}{\tau}\right), & \tau > 0,\\[2mm]
\delta_{\{\boldsymbol{0}\}}(\boldsymbol{s}), & \tau =0.
\end{cases}
\]
Note that 
\[
\inf_{\tau \ge 0} ((1+\sigma_{\mathcal C'}^*)\tau)(\boldsymbol{s})
=
\inf\{\tau \ge 0 : \boldsymbol{s} \in \tau \mathcal C'\}=\inf\{\tau \ge 0 : \boldsymbol{s_x}\in \tau \mathcal C,\ s_0=0\} =: \gamma_{\mathcal{C}}(\boldsymbol{s_x}).
\]
where the first equality follows from the standard gauge identity \cite[p.~35]{rockafellar1997convex} and the second equality is from \(\mathcal C'=\mathcal C\times\{0\}\).

Substituting this into \eqref{eq:cgp-dual:2} yields \eqref{prob:dcgp}, together with the relation
$$(-\boldsymbol{\hat x},1)-(\boldsymbol{s_x},0) = (\boldsymbol{x},z).$$ 
For notational simplicity, we henceforth let $\boldsymbol{s}$ denote $\boldsymbol{s_x}$.

The final equivalence follows from Proposition~\ref{prop:reverse}.
\hfill $\qed$

\subsection{Proof of Proposition \ref{prop:supporting}}\label{proof:prop:supporting}
\noindent\textbf{Proof.}
By dualizing \eqref{eq:dcgp:b}--\eqref{eq:dcgp:c} and rearranging terms, we obtain
\begin{equation}
\sup_{\boldsymbol\alpha} \ \alpha_0 - \boldsymbol{\alpha}_x^\top \boldsymbol{\hat{x}}
+ \inf_{\boldsymbol s} \left\{ \gamma_{\mathcal C}(\boldsymbol{s}) - \boldsymbol{\alpha}_x^\top \boldsymbol{s} \right\}
+ \inf_{(-\boldsymbol{x}, z) \in (\mathcal{S}^*)^{\circ}} \left\{ \boldsymbol{\alpha}_x^\top \boldsymbol{x} - \alpha_0 z \right\}.
\label{eq:dcgp:lag-dual}
\end{equation}

By the KKT stationarity condition with respect to $\boldsymbol{s}$, $\boldsymbol{\hat\alpha}_x$ is a subgradient of $\gamma_{\mathcal C}$ at $\boldsymbol{\hat s}$. Hence, by the Fenchel-Young inequality,
\[
\boldsymbol{\hat\alpha_x}^\top \boldsymbol{\hat s}
= \gamma_{\mathcal C}(\boldsymbol{\hat s}) + \gamma_{\mathcal C}^*(\boldsymbol{\hat \alpha_x}).
\]
Since $\gamma_{\mathcal C}^* = \delta_{\mathcal C^\circ}$, it follows that
\[
\boldsymbol{\hat\alpha_x} \in \mathcal C^\circ
\quad \text{and} \quad
\boldsymbol{\hat\alpha_x}^\top \boldsymbol{\hat s} = \gamma_{\mathcal C}(\boldsymbol{\hat s}) = \hat \tau.
\]
Using $\boldsymbol{\hat s} = \boldsymbol{x'} - \boldsymbol{\hat{x}}$, we obtain
\[
\boldsymbol{\hat\alpha_x}^\top \boldsymbol{x'} = \hat\alpha_0.
\]

The third term in \eqref{eq:dcgp:lag-dual} equals $-\sigma_{(\mathcal S^*)^\circ}(\boldsymbol\alpha) = -\delta_{\mathcal S^*}(\boldsymbol\alpha)$, since $(\mathcal S^*)^\circ$ is a closed convex cone. Thus, $\boldsymbol{\hat\alpha} \in \mathcal S^*$, implying that it defines a valid inequality for $\overline{\operatorname{conv}}(\mathcal S)$.

Therefore, the inequality is valid and tight at $\boldsymbol{x'}$, and hence defines a supporting hyperplane.

Finally, if $\hat\tau = 0$, then $\boldsymbol{\hat s} = \boldsymbol{0}$, implying $\boldsymbol{x'} = \boldsymbol{\hat{x}} \in \overline{\operatorname{conv}}(\mathcal{S})$. \qed
}
\subsection{Proof of Proposition \ref{prop:algorithm-cut}} \label{proof:prop:algorithm-cut}
{\color{black}Note that, by the definition of homogenization,
$$
(\boldsymbol{x},t, z) \in \mathcal L^\# \Leftrightarrow 
\begin{cases}
(\boldsymbol{x}/z, t/z) \in \mathcal L, \ & \mbox{ if }z>0,\\
(\boldsymbol{x},t) \in \operatorname{recc}(\mathcal L), \ & \mbox{ if }z = 0.
\end{cases}$$
Therefore, if $(\boldsymbol{x},t, z) \in \mathcal L^\#$,
$$\beta_0 z-\boldsymbol{\beta_x}^\top\boldsymbol{x} -\beta_t t \le 0, \ \forall \boldsymbol{\beta} \in \mathcal L^*.$$
Thus, any valid inequality for $\mathcal L$ induces a valid inequality for $\mathcal L^\#$. Now suppose that $\boldsymbol{\beta} \in \mathcal L^*$ separates $(\boldsymbol{x}',t') \notin \mathcal L$, i.e.,
\[
\boldsymbol{\beta}_x^\top \boldsymbol{x}' + \beta_t t' < \beta_0.
\]
Multiplying by $\hat\omega_0 > 0$ and using $(\boldsymbol{x}',t') = (\boldsymbol{\hat\omega_x}/\hat\omega_0, \hat\omega_t/\hat\omega_0)$, we obtain
\[
\boldsymbol{\beta}_x^\top \boldsymbol{\hat\omega_x} + \beta_t \hat\omega_t < \beta_0 \hat\omega_0,
\]
which shows that $(\boldsymbol{\hat\omega_x}, \hat\omega_t, \hat\omega_0)$ violates the inequality and is therefore separated from $\mathcal L^\#$.

Conversely, if $(\boldsymbol{x'}, t') \in \mathcal{L}$, then by definition of homogenization,
\[
(\boldsymbol{\hat\omega_x}, \hat\omega_t, \hat\omega_0)
= \hat\omega_0 (\boldsymbol{x'}, t', 1) \in \mathcal{L}^\#.
\]\qed
}
   
\subsection{Proof of Proposition \ref{prop:inexact}}\label{proof:prop:inexact}
\texttt{R}$^{(\rho)}$ corresponds to \eqref{prob:dcglp} with the polyhedral set $\mathcal L$ replaced by its polyhedral relaxation $\mathcal L^{(\rho)}$. Thus, the statements follow directly from Corollary \ref{coro:dcglp}.\qed

\subsection{Proof of Theorem \ref{theo:finite-convergence}}\label{proof:theo:finite-convergence}
By design, Lines \ref{line:special:disjunctive:begin}-\ref{line:special:disjunctive:end} in Algorithm \ref{algo:special} are executed only when $(\boldsymbol{\hat x}, \hat t) \in \arg\min \{\boldsymbol{c}^\top \boldsymbol{x} +t : (\boldsymbol{x},t) \in \mathcal L, \ \boldsymbol{x} \in \mathcal X_{LP}, \ (\boldsymbol{x},t) \in \mathcal D\}$, where $\mathcal D$ is the polyhedron defined by the disjunctive cuts found up to that point.  
Therefore, Algorithm \ref{algo:special} can be interpreted exactly as the specialized, finitely convergent lift-and-project algorithm for solving $\min \{\boldsymbol{c}^\top \boldsymbol{x} + t: (\boldsymbol{x},t) \in \mathcal L, \ \boldsymbol{x} \in \mathcal X\}$, which was proposed in \cite{balas1993lift}, with the continuous relaxation problems between consecutive disjunctive cut separations are solved using standard BD. Specifically, Algorithm \ref{algo:special} initially solves the continuous relaxation of \eqref{prob:benders-reformulation} using standard BD, and each time a new disjunctive cut is added to $\mathcal D$ in Line \ref{line:special:disjunctive:add}, Algorithm \ref{algo:special} solves the updated linear program $\min \{\boldsymbol{c}^\top \boldsymbol{x} +t : (\boldsymbol{x},t) \in \mathcal L, \ \boldsymbol{x} \in \mathcal X_{LP}, \ (\boldsymbol{x},t) \in \mathcal D\}$ again using standard BD. 

From Theorem 3.1 in \cite{balas1993lift} (Theorem 5.24 in \cite{conforti2014integer}), it follows that a finite number of disjunctive cuts are required for termination.
In addition, for each addition of a new disjunctive cut, Benders decomposition is finitely convergent for solving the updated linear program. Therefore, the overall algorithm is guaranteed to converge in a finite number of iterations. \qed

\subsection{Proof of Proposition \ref{prop:strengthen}}\label{proof:prop:strengthen}
Consider a disjunction:
$$x_i - m^\top x \le 0 \ \lor \  x_i - m^\top x \ge 1,$$
for some integral vector $m \in \mathbb Z^{n_x}$ with $m_j = 0$ for $j \not\in \mathcal I$. 
Note that since $\underline{x}_j =0$ for $j \in \mathcal I$, $\hat\nu^r_j$ serves as slack variables in \eqref{eq:cglp:c} for $j \in \mathcal I$. Therefore, we can construct a valid cut $\hat\alpha_0' - (\boldsymbol{\hat\alpha_x'})^\top \boldsymbol{x}- \hat\alpha_t' t \le 0$ for $\mathcal P_{\mathcal R}$ with
\[\begin{cases}
\hat\alpha_0'  = \hat\alpha_0,\\
(\boldsymbol{\hat\alpha_{x}'})_j  = \max\{a^1_j + \hat \kappa^{1} m_j, a^2_j - \hat \kappa^{2} m_j\}, \ \forall j \in \mathcal I,\\
(\boldsymbol{\hat\alpha_{x}'})_j  = (\boldsymbol{\hat\alpha_{x}})_j, \ \forall j \in [n_x]\setminus \mathcal I,\\
\hat\alpha_t'  = \hat \alpha_t.
\end{cases}\]
 
\begin{figure}
\centering
\includegraphics[width=0.25\textwidth]{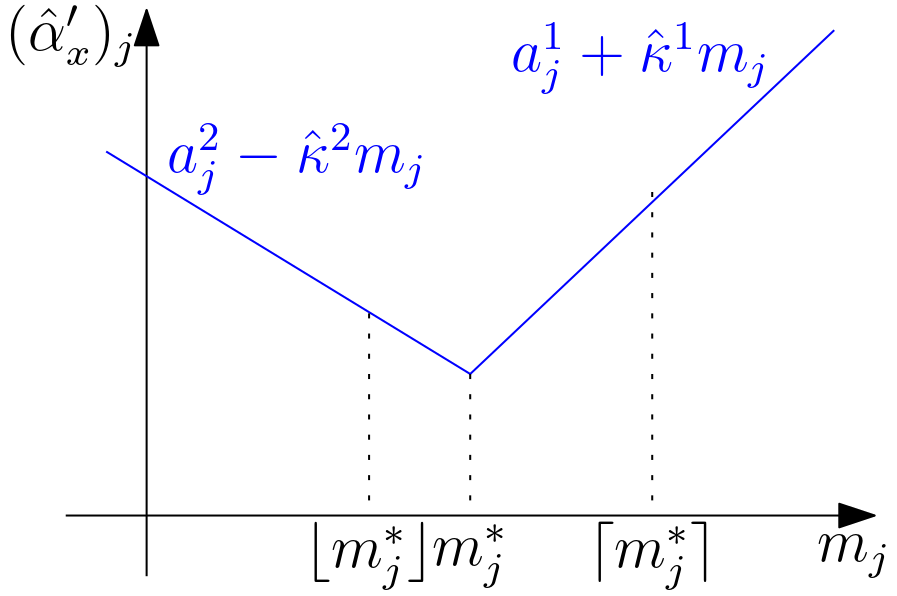}
\caption{The graph of $(\boldsymbol{\hat\alpha_{x}'})_j$ as a function of $m_j$.}
\label{fig:strengthen}
\end{figure}

\noindent Note that $(\boldsymbol{\hat\alpha_{x}'})_j$ is a function of $m_j$, as illustrated in Figure \ref{fig:strengthen}. Since $x_j \ge 0$ for $j \in \mathcal I$, the lower the value of the corresponding $(\boldsymbol{\hat\alpha_{x}'})_j$ is, the tighter the cut $\hat\alpha_0' - (\boldsymbol{\hat\alpha_x'})^\top \boldsymbol{x}- \hat\alpha_t' t \le 0$ is. Thus, we choose the value of $m_j$ so that $(\boldsymbol{\hat\alpha_{x}'})_j=\max\{{a}^1_j + \hat \kappa^{1} m_j, a^2_j - \hat \kappa^{2} m_j\}$ is minimized for each $j\in \mathcal I$. Since both $\hat\kappa^{r}$'s are negative, the $m_j$-value that gives the minimum value of $(\boldsymbol{\hat\alpha_{x}'})_j$ is the point where the two affine functions meet (see Figure \ref{fig:strengthen}), i.e., 
\begin{align*}
a_j^1 + \hat \kappa^{1}  m_j^* = a_j^2 - \hat \kappa^{2}  m_j^* \Rightarrow  m_j^* =\frac{a_j^{2} - a_j^{1}}{\hat \kappa^{1} + \hat \kappa^{2}}
\end{align*}
Since $m_j$ should be an integer, we can choose either $\lceil  m_j^* \rceil$ or $\lfloor  m_j^* \rfloor$ whichever gives a smaller value of $(\boldsymbol{\hat\alpha_{x}'})_j$ (see Figure \ref{fig:strengthen}), i.e., 
\begin{equation*}(\boldsymbol{\hat\alpha_{x}'})_j = \min\left\{a^1_j + \hat \kappa^{1} \left\lfloor \frac{a_j^2 - a_j^1}{\hat \kappa^{1} + \hat \kappa^{2}} \right\rfloor, a^2_j -\hat \kappa^{2} \left\lceil  \frac{a_j^2 - a_j^1}{\hat \kappa^{1} + \hat \kappa^{2}} \right\rceil\right\}, \ \forall j \in \mathcal I.
\end{equation*}


\qed
\subsection{Proof of Proposition \ref{prop:lifting}}\label{proof:prop:lifting}

We can verify $\boldsymbol{\widehat{\alpha}} \in \mathcal P_i$ directly by checking whether it satisfies \eqref{eqs:valid-cut}. To do so, we express $\mathcal P_i^\n$ and $\mathcal P_i$ as follows, using some appropriate matrices and vectors $A^r, \boldsymbol{b^r}, \boldsymbol{e}$:\\
\begin{minipage}{0.48\textwidth}\scriptsize
\centering
\begin{align*}
& \mathcal P_i^\n:\\
& [\boldsymbol{\psi^1}] \ A^1 \boldsymbol{x} + \boldsymbol{e}t \ge \boldsymbol{b^1} &&   \ [\boldsymbol{\psi^2}] \ A^2 \boldsymbol{x} + \boldsymbol{e}t \ge \boldsymbol{b^2},\\
& [\xi^1_j] \ x_j \ge 1, \forall j \in \mathcal I_{\n,1} && \ [\xi^2_j] \ x_j \ge 1, \forall j \in \mathcal I_{\n,1}\\
& [\zeta^1_j] \ -x_j \ge 0, \forall j \in \mathcal I_{\n,0} && \ [\zeta^2_j] \ -x_j \ge 0, \forall j \in \mathcal I_{\n,0}\\
& [\nu^1_j] \ x_j \ge 0, \forall j \in \mathcal I && \ [\nu^2_j] \ x_j \ge 0, \forall j \in \mathcal I\\
& [\eta^1_j] \ -x_j \ge -1, \forall j \in \mathcal I && \ [\eta^2_j] \ -x_j \ge -1, \forall j \in \mathcal I\\
\end{align*}
    \end{minipage}
    \hfill
    \begin{minipage}{0.48\textwidth}\scriptsize
\centering
\begin{align*}
& \mathcal P_{i}:\\
& [\boldsymbol{\psi^1}] \ A^1 \boldsymbol{x} + \boldsymbol{e}t \ge \boldsymbol{b^1} &&   \ [\boldsymbol{\psi^2}] \ A^2 \boldsymbol{x} + \boldsymbol{e}t \ge \boldsymbol{b^2},\\
& [\nu^1_j] \ x_j \ge 0, \forall j \in \mathcal I && \ [\nu^2_j] \ x_j \ge 0, \forall j \in \mathcal I\\
& [\eta^1_j] \ -x_j \ge -1, \forall j \in \mathcal I && \ [\eta^2_j] \ -x_j \ge -1, \forall j \in \mathcal I\\
\end{align*}
    \end{minipage}
The characterization of valid inequalities for $\mathcal P_i^\n$ and $\mathcal P_i$ are as follows \cite[Theorem 1.2]{balas2018disjunctive}:\\
\begin{minipage}{0.48\textwidth}\scriptsize
\centering
\begin{align*}
\alpha_t & = \boldsymbol{e}^\top \boldsymbol{\psi^1} = \boldsymbol{e}^\top \boldsymbol{\psi^2},\\
(\boldsymbol{\alpha_{x}})_j  & = (A^1_{\cdot j})^\top \boldsymbol{\psi^1} = (A^2_{\cdot j})^\top \boldsymbol{\psi^2}, \ \forall j \in [n_x] \setminus \mathcal I,\\
(\boldsymbol{\alpha_{x}})_j & = (A^1_{\cdot j})^\top \boldsymbol{\psi^1} + \nu^1_j - \eta^1_j \\
& = (A^2_{\cdot j})^\top \boldsymbol{\psi^2} + \nu^2_j - \eta^2_j, \forall j \in \mathcal I \setminus \mathcal I_{\n,0}\cup \mathcal I_{\n,1},\\
(\boldsymbol{\alpha_{x}})_j & = (A^1_{\cdot j})^\top \boldsymbol{\psi^1} + \nu^1_j - \eta^1_j - \zeta_j^1 \\
& = (A^2_{\cdot j})^\top \boldsymbol{\psi^2} + \nu^2_j - \eta^2_j- \zeta_j^2, \forall j \in \mathcal I_{\n,0},\\
(\boldsymbol{\alpha_{x}})_j & = (A^1_{\cdot j})^\top \boldsymbol{\psi^1} + \nu^1_j - \eta^1_j + \xi_j^1 \\
& = (A^2_{\cdot j})^\top \boldsymbol{\psi^2} + \nu^2_j - \eta^2_j+ \xi_j^2, \forall j \in \mathcal I_{\n,1},\\
\alpha_{0} & \le (\boldsymbol{b^1})^\top \boldsymbol{\psi^1} - \sum_{j \in \mathcal I}\eta_j^1 + \sum_{j \in \mathcal I_{\n,1}} \xi_j^1,\\
\alpha_{0} & \le (\boldsymbol{b^2})^\top \boldsymbol{\psi^2} - \sum_{j \in \mathcal I}\eta_j^2 + \sum_{j \in \mathcal I_{\n,1}} \xi_j^2.
\end{align*}
    \end{minipage}
    \hfill
    \begin{minipage}{0.48\textwidth}\scriptsize
\centering
\begin{align*}
\alpha_t & = \boldsymbol{e}^\top \boldsymbol{\psi^1} = \boldsymbol{e}^\top \boldsymbol{\psi^2},\\
(\boldsymbol{\alpha_{x}})_j  & = (A^1_{\cdot j})^\top \boldsymbol{\psi^1} = (A^2_{\cdot j})^\top \boldsymbol{\psi^2}, \ \forall j \in [n_x] \setminus \mathcal I,\\
(\boldsymbol{\alpha_{x}})_j &  = (A^1_{\cdot j})^\top \boldsymbol{\psi^1} + \nu^1_j - \eta^1_j \\
& = (A^2_{\cdot j})^\top \boldsymbol{\psi^2} + \nu^2_j - \eta^2_j, \forall j \in \mathcal I \setminus \mathcal I_{\n,0}\cup \mathcal I_{\n,1},\\
 (\boldsymbol{\alpha_{x}})_j &= (A^1_{\cdot j})^\top \boldsymbol{\psi^1} + \nu^1_j - \eta^1_j \\
&= (A^2_{\cdot j})^\top \boldsymbol{\psi^2} + \nu^2_j - \eta^2_j, \forall j \in \mathcal I_{\n,0},\\
(\boldsymbol{\alpha_{x}})_j & = (A^1_{\cdot j})^\top \boldsymbol{\psi^1} + \nu^1_j - \eta^1_j \\
&= (A^2_{\cdot j})^\top \boldsymbol{\psi^2} + \nu^2_j - \eta^2_j, \forall j \in \mathcal I_{\n,1},\\
 \alpha_{0} & \le (\boldsymbol{b^1})^\top \boldsymbol{\psi^1} - \sum_{j \in \mathcal I}\eta_j^1,\\
 \alpha_{0} & \le (\boldsymbol{b^2})^\top \boldsymbol{\psi^2} - \sum_{j \in \mathcal I}\eta_j^2,
\end{align*}
    \end{minipage}
Thus, $(\boldsymbol{\hat\alpha}, \boldsymbol{\hat\chi}, (\boldsymbol{\hat\zeta^{r}},\boldsymbol{\hat\xi^{r}})_{r =1,2})$ satisfies the system on the left. Consequently, the modification in \eqref{eq:lifted-cut} results in:
{\scriptsize\begin{align*}
\alpha_t & = e^\top \boldsymbol{\psi^1} = e^\top \boldsymbol{\psi^2},\\
(\boldsymbol{\alpha_{x}})_j  & = (A^1_{\cdot j})^\top \boldsymbol{\psi^1} = (A^2_{\cdot j})^\top \boldsymbol{\psi^2}, \ \forall j \in [n_x] \setminus \mathcal I,\\
(\boldsymbol{\alpha_{x}})_j & = (A^1_{\cdot j})^\top \boldsymbol{\psi^1} + \nu^1_j - \eta^1_j \\
& = (A^2_{\cdot j})^\top \boldsymbol{\psi^2} + \nu^2_j - \eta^2_j, \forall j \in \mathcal I \setminus \mathcal I_{\n,0}\cup \mathcal I_{\n,1},\\
{\color{black}(\boldsymbol{\widehat{\alpha_{x}}})_j} & = (A^1_{\cdot j})^\top \boldsymbol{\psi^1} + \nu^1_j - \eta^1_j - \zeta_j^1{\color{black} +\max\{\zeta_j^1, \zeta_j^2\}} \\
& = (A^2_{\cdot j})^\top \boldsymbol{\psi^2} + \nu^2_j - \eta^2_j- \zeta_j^2{\color{black} +\max\{\zeta_j^1, \zeta_j^2\}}, \forall j \in \mathcal I_{\n,0},\\
{\color{black}(\boldsymbol{\widehat{\alpha_{x}}})_j} & = (A^1_{\cdot j})^\top \boldsymbol{\psi^1} + \nu^1_j - \eta^1_j + \xi_j^1 {\color{black} -\max\{\xi^1, \xi_j^2\}}\\
& = (A^2_{\cdot j})^\top \boldsymbol{\psi^2} + \nu^2_j - \eta^2_j+ \xi_j^2{\color{black} -\max\{\xi^1, \xi_j^2\}}, \forall j \in \mathcal I_{\n,1},\\
{\color{black}\widehat{\alpha_{0}}} & \le (\boldsymbol{b^1})^\top \boldsymbol{\psi^1} - \sum_{j \in \mathcal I}\eta_j^1 + \sum_{j \in \mathcal I_{\n,1}} (\xi_j^1{\color{black} -\max\{\xi^1, \xi_j^2\}}),\\
{\color{black}\widehat{\alpha_{0}}} & \le (\boldsymbol{b^2})^\top \boldsymbol{\psi^2} - \sum_{j \in \mathcal I}\eta_j^2 + \sum_{j \in \mathcal I_{\n,1}} (\xi_j^2{\color{black} -\max\{\xi^1, \xi_j^2\}}),
\end{align*}
}
which show clearly $\boldsymbol{\widehat{\alpha}}$, together with $\boldsymbol{\widehat{\chi}}$ in \eqref{eq:lifting:farkas}, satisfies the valid inequality characterization for $\mathcal P_i$.

Morevoer, we have $\widehat{\alpha_0} - \boldsymbol{\widehat{\alpha_x}}^\top  \boldsymbol{\hat x} - \widehat{\alpha_t}\hat t =
\hat\alpha_0 - \boldsymbol{\hat\alpha_x}^\top \boldsymbol{\hat x} - \hat\alpha_t \hat t>0$.
\qed

\subsection{Proof of Proposition \ref{prop:approx-oracle}}
Upon full convergence of the while loop in Line \ref{oracle-approx:while}, $(\tau, \boldsymbol\alpha)$ is an optimal solution to (\texttt{DCGLP}$_{\mathcal{I}}$) where (\texttt{DCGLP}$_{\mathcal{I}}$) is a restriction of \eqref{prob:dcglp} incorporating the additional constraints from Line \ref{oracle:approx:const}.
Therefore, $\hat\tau \geq v_{\eqref{prob:dcglp}} = v_{\eqref{prob:cgp}}$. 
From Proposition \ref{prop:lifting}, we have $\widehat{\alpha} \in \mathcal S^* := \mathcal P_{\mathcal R}^*$. 
Since $\mathcal S^*$ is conic and $\sigma_{\mathcal C}$ is real-valued and positively homogeneous, $\boldsymbol{\widehat{\alpha}}$ can be made feasible to \eqref{prob:cgp} by normalizing it with a positive scalar $\frac{1}{\max\{\sigma_{\mathcal C}(\boldsymbol{\widehat{\alpha_x}}),1\}}$. 
The resultant feasible solution evaluates the objective function of \eqref{prob:cgp} at $ \frac{\hat\tau}{\max\{\sigma_{\mathcal C}(\boldsymbol{\widehat{\alpha_x}}), 1\}}$, and thus $v_\eqref{prob:cgp} \geq \frac{\hat\tau}{\max\{\sigma_{\mathcal C}(\boldsymbol{\widehat{\alpha_x}}), 1\}}$.\qed

\section{Formulation of benchmark problems}
\subsection{Uncapacitated Facility Location Problem}\label{formu:uflp}

The Uncapacitated Facility Location Problem (UFLP) is formulated as a MILP to minimize total costs, including facility opening costs and transportation costs. Given $I$ as the total number of potential facility locations and $J$ as the total number of customer nodes, let $c_{i,j}$ be the transportation cost from facility $i$ to customer $j$. Binary variables $x_i$ indicate whether facility $i$ is opened, and continuous variables $y_{i,j}$ represent the fraction of customer $j$'s demand satisfied by facility $i$. The classical UFLP model is expressed as follows:

\begin{equation*}
\begin{aligned}
& \min_{x, y} \quad && \sum_{i \in [I]} f_i x_i + \sum_{i \in [I]} \sum_{j \in [J]} c_{i,j}y_{i,j} \\
& \text{s.t.} \quad && \sum_{i \in [I]} y_{i,j} = 1, && \forall j \in [J] \\
& && y_{i,j} \leq x_i, && \forall i \in [I], j \in [J] \\
& && x_i \in \{0,1\}, \,\, y_{i,j} \geq 0, && \forall i \in [I], j \in [J]
\end{aligned}
\end{equation*}

Since the UFLP has no coupling constraint, the subproblem becomes separable with respect to each customer $j$. In addition to the classical UFLP model, we consider $\sum_{i \in [I]}x_i \geq 2$ in the master problem of BD. Since we can easily compute the optimal objective value for the case when $\sum_{i \in [I]}x_i = 1$, we exclude it in the BD framework as proposed by \cite{fischetti2017redesigning}.





\subsection{Stochastic Network Interdiction Problem}\label{formu:snip}

The Stochastic Network Interdiction Problem (SNIP) aims to minimize the expected probability of an intruder reaching the destination across scenarios by deploying sensors on eligible arcs. Let $N$ denote the total number of nodes and $A$ set of arcs, with $D \subseteq A$ as the subset of arcs where sensors can be placed. Let $K$ denote the total number of scenarios. In the first-stage problem, given a budget $b$, the interdictor places sensors on arcs, knowing a priori probability of an intruder avoiding detection with and without a sensor on arc $(i,j)$ denoted as $r_{i,j}$ and $q_{i,j}$ respectively. The interdictor makes the decision without knowing the origin ($o_k$) and the destination ($d_k$) of the intruder for each scenario $k \in [K]$. In the second-stage problem, the intruder maximizes the likelihood of evading surveillance by selecting a maximum reliability path. When no sensors are placed, $\psi_{j,k}$ denotes the value of a maximum reliability path from node $j \in [N]$ to the destination $d_k$, computed by solving a shortest path problem \cite{pan2008minimizing}, and is used to adjust the evasion probability based on the installation. Binary variables $x_{i,j}$, used in the first-stage, indicate whether a sensor is deployed on arc $(i,j) \in D$, while continuous variables $y_{i,k}$, used in the second-stage, represent the probability of the intruder traversing from node $i \in [N]$ to destination $d_k$ without detection. The extensive formulation of SNIP is as follows:

\begin{equation*}
\begin{aligned}
& \min_{x,y} \quad && \sum_{k \in [K]} p_k y_{o_k,k} \\
& \text{s.t.} \quad && \sum_{(i,j) \in D} x_{ij} \leq b\\
& && y_{d_k,k} = 1, && \forall k \in [K] \\
& && y_{i,k} - r_{i,j}y_{j,k} \geq 0, && \forall (i,j) \in A \backslash D, k \in [K] \\
& && y_{i,k} - r_{i,j}y_{j,k} \geq -(r_{i,j}-q_{i,j})\psi_{j,k}x_{i,j}, && \forall (i,j) \in D, k \in [K] \\
& && y_{i,k} - q_{i,j}y_{j,k} \geq 0, && \forall (i,j) \in D, k \in [K] \\
& && x_{i,j} \in \{0,1\}^{|D|}, \,\, y_{i,k} \geq 0, && \forall i \in [N], k \in [K]
\end{aligned}
\end{equation*}

{\color{black}
\section{Additional experiments for \textbf{UFLP} with Gurobi} \label{apx:gurobi}

\begin{figure}[b!]
    \centering
    \begin{subfigure}[b]{0.45\textwidth}
        \includegraphics[clip,width=\textwidth]{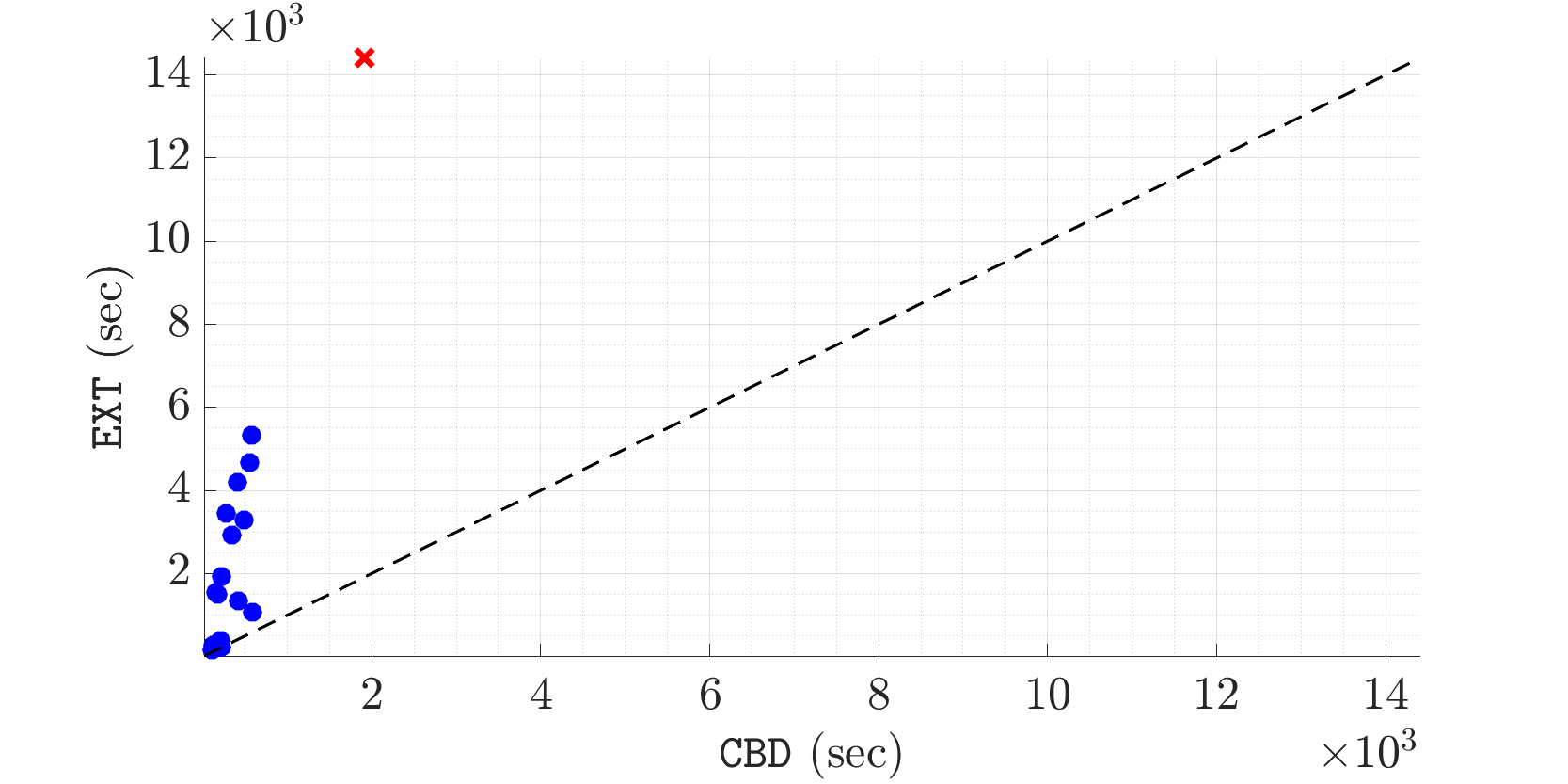}
        \caption{Time, \texttt{CBD} and \texttt{EXT}}
        \label{fig:time-ext-250-gurobi}
    \end{subfigure}
 \begin{subfigure}[b]{0.45\textwidth}
        \includegraphics[clip,width=\textwidth]{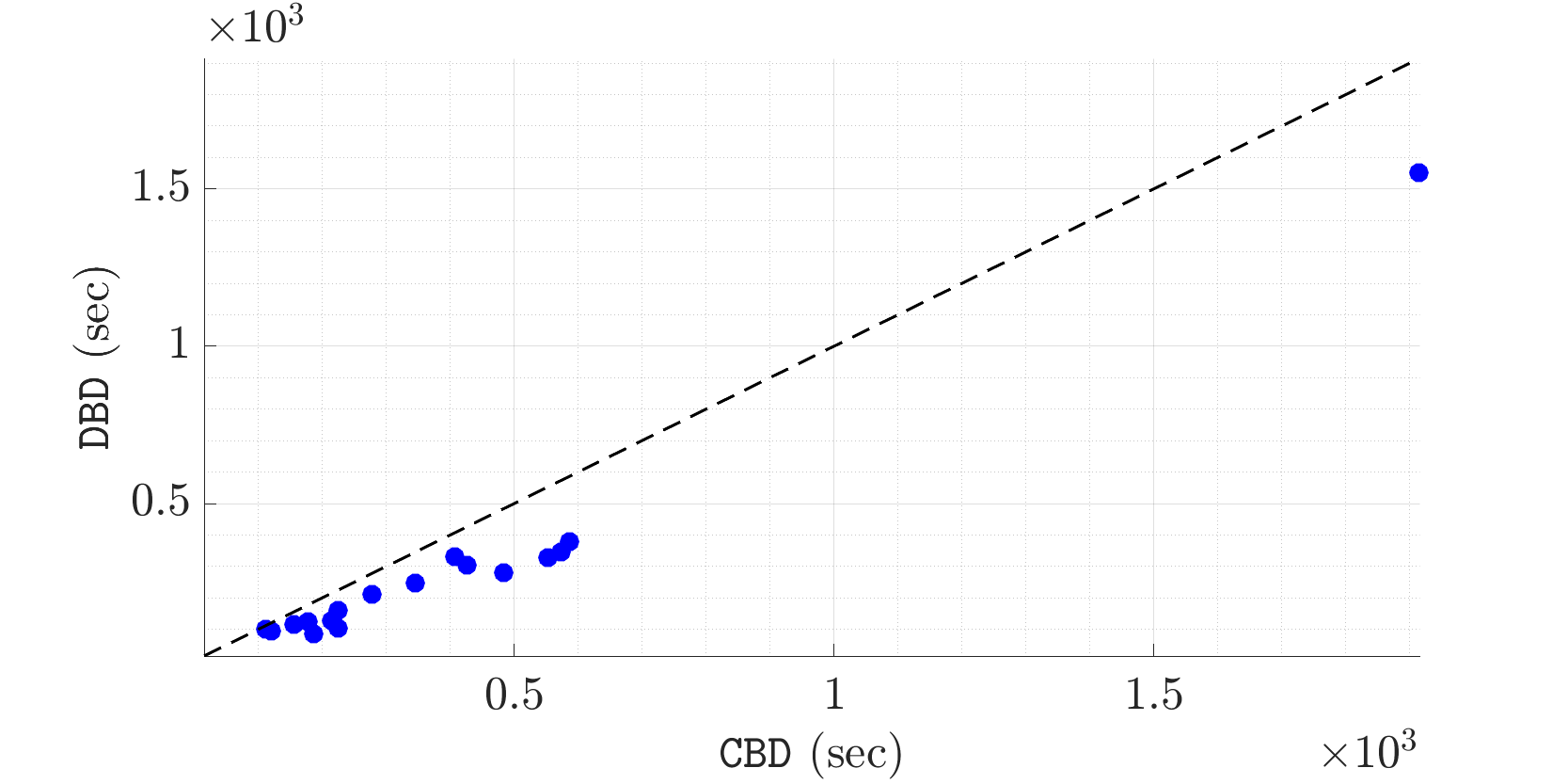}
        \caption{Time, \texttt{CBD} and \texttt{DBD}}
        \label{fig:time-dbd-250-gurobi}
    \end{subfigure}
    \begin{subfigure}[b]{0.45\textwidth}  
        \includegraphics[clip,width=\textwidth]{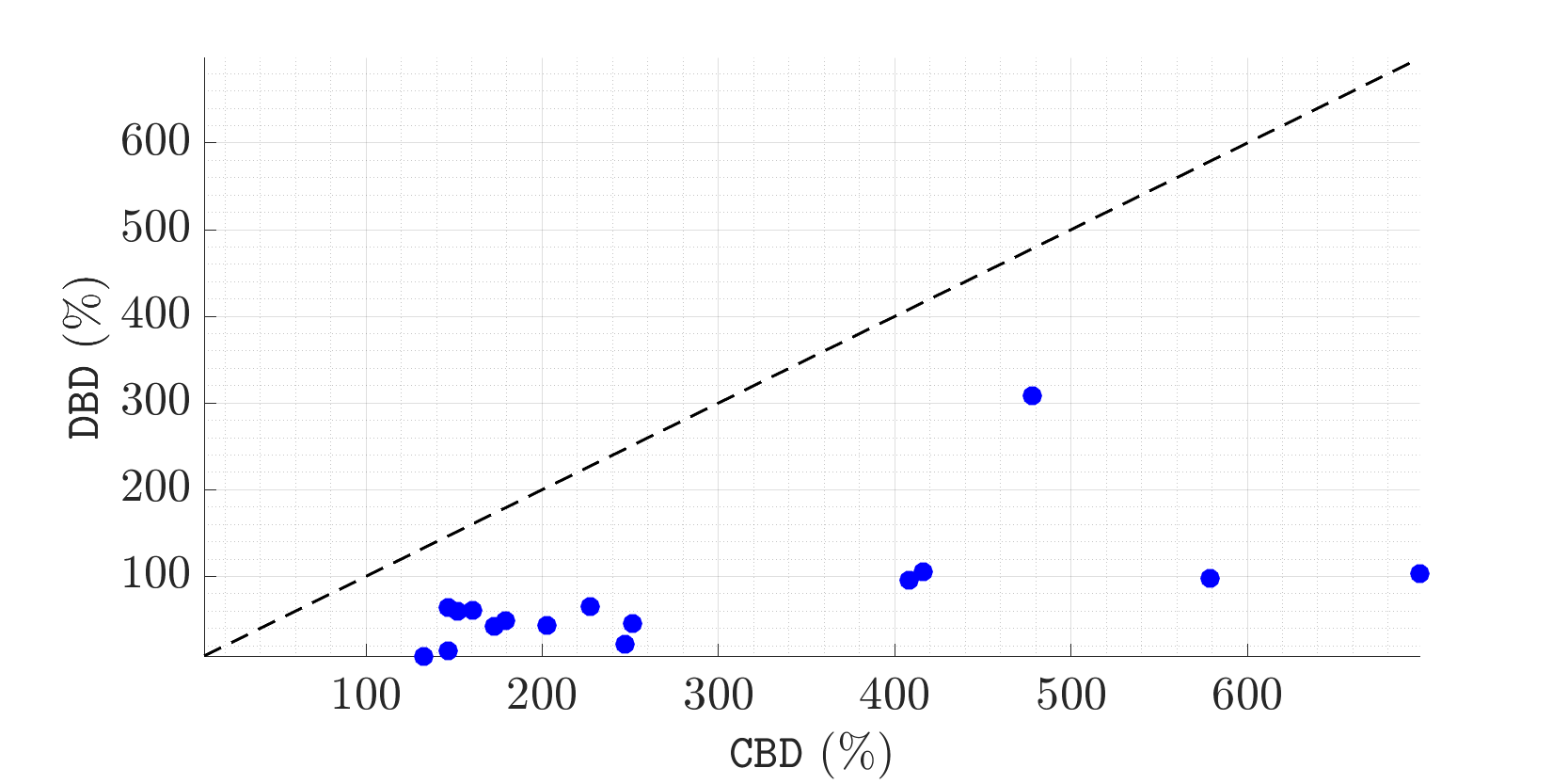}
        \caption{Relative increase in NC with respect to \texttt{EXT} on instances solved to optimality by \texttt{EXT}}
        \label{fig:nc-rel-250-gurobi}
    \end{subfigure}
    \begin{subfigure}[b]{0.45\textwidth}
        \includegraphics[clip,width=\textwidth]{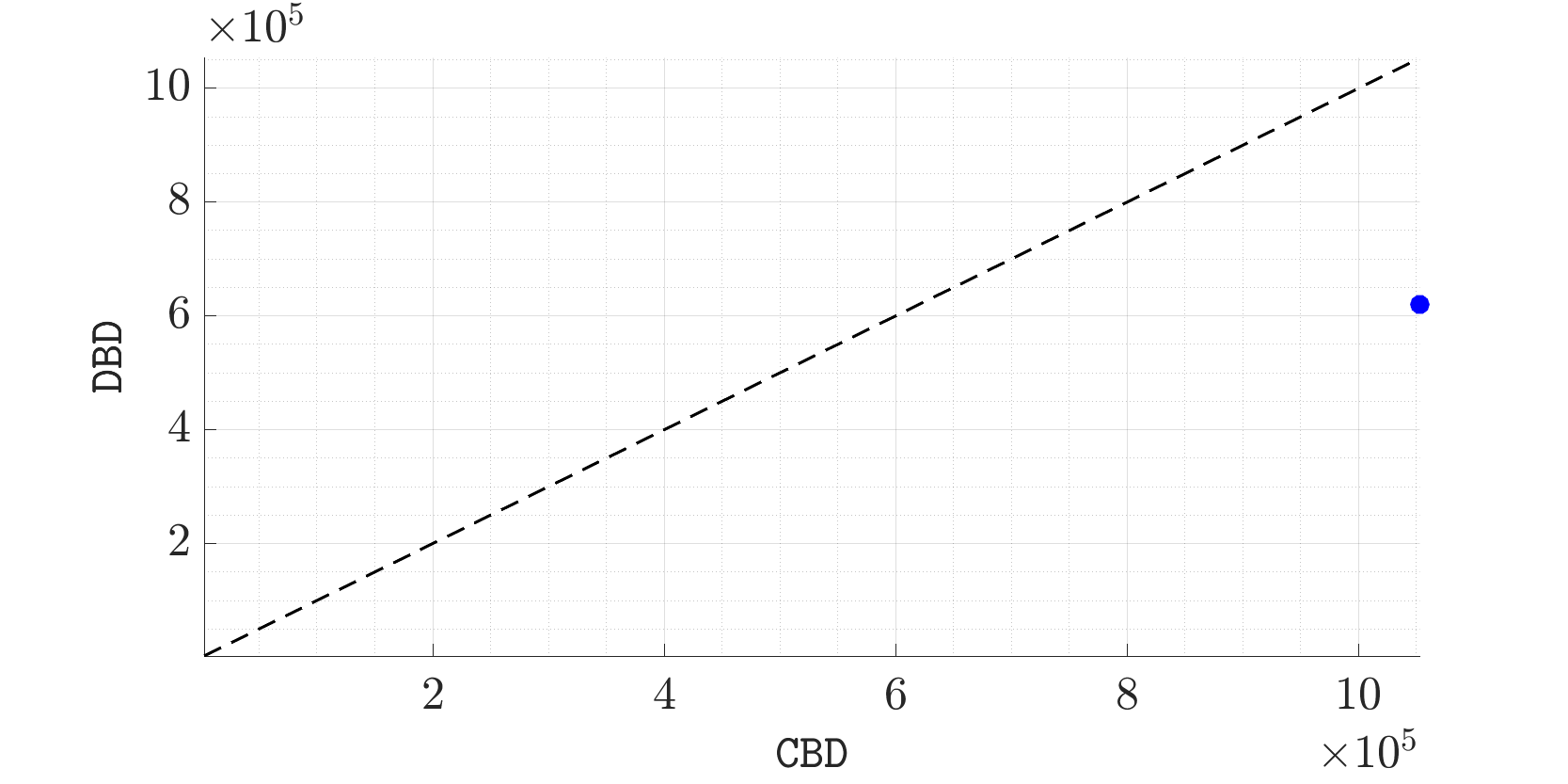}
        \caption{NC on instances for which \texttt{EXT} reached the time limit.}
        \label{fig:nc-abs-250-gurobi}
    \end{subfigure}
    \caption{UFLP-250 instances; the \texttt{x} marker denote the instances for which \texttt{EXT} reached the 4-hour time limit.}
    \label{fig:250-gurobi}
\end{figure}

Gurobi exhibits trends similar to those discussed in Section~\ref{sec:experiment:uflp}, although the performance of all three methods improves noticeably. The results suggest that Gurobi internally generates strong cutting planes from the Benders cuts added during the solution process. Consequently, the performance gap between \texttt{CBD} and \texttt{DBD} becomes smaller. Nevertheless, \texttt{DBD} still outperforms \texttt{CBD} on the majority of instances.

\begin{figure}[t!]
    \centering
    \begin{subfigure}[b]{0.45\textwidth}  
        \includegraphics[clip,width=\textwidth]{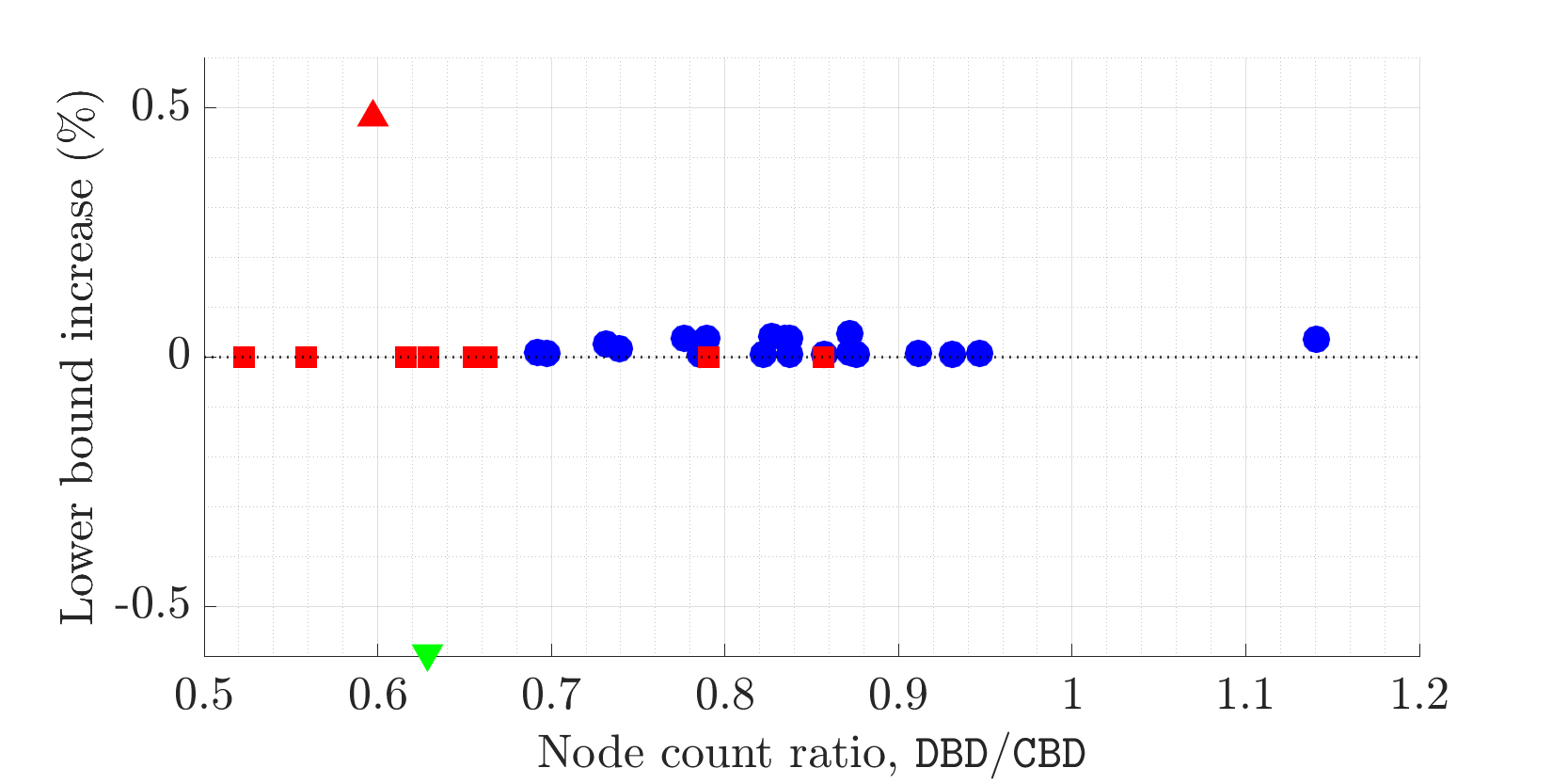}
        \caption{UFLP-500 instances}
        \label{fig:nc-500-gurobi}
    \end{subfigure}
    \hspace{4mm}
    \begin{subfigure}[b]{0.45\textwidth}
        \includegraphics[clip,width=\textwidth]{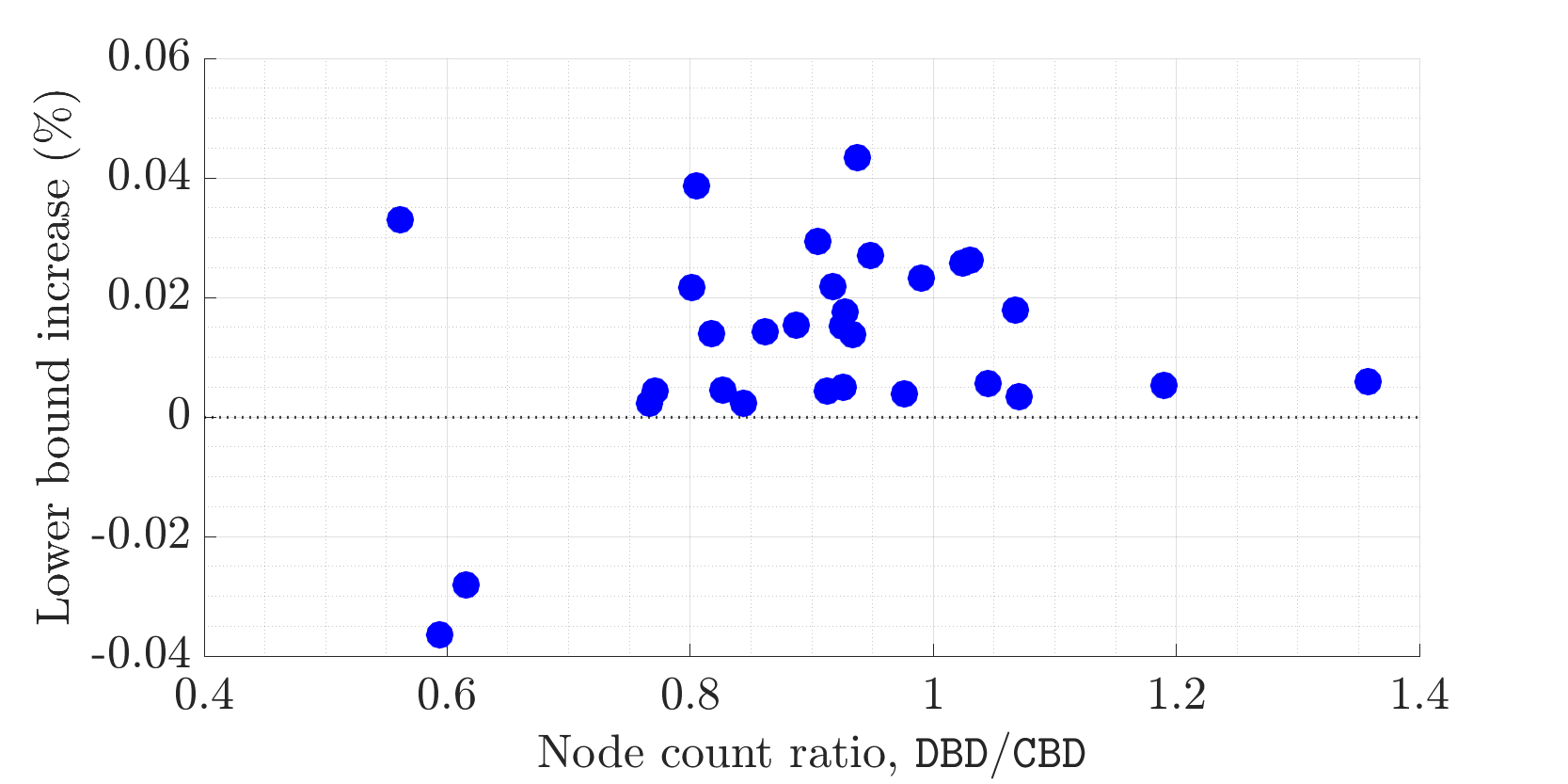}
        \caption{UFLP-750 instances}
        \label{fig:nc-750-gurobi}
    \end{subfigure}
   \caption{Relationship between the NC ratio and the relative lower-bound improvement. Triangles (resp.\ upside-down triangles) denote instances for which only \texttt{CBD} (resp.\ only \texttt{DBD}) reaches the time limit. Squares denote instances for which both methods solve to optimality.}
    \label{fig:nc-500-750-gurobi}
\end{figure}
}

\end{appendix}

\end{document}